\documentclass[11pt]{amsart}
\usepackage{amssymb}
\usepackage{amscd}

\newtheorem{theorem}{Theorem}[section]
\newtheorem{corollary}[theorem]{Corollary}

\newtheorem{lemma}[theorem]{Lemma}
\newtheorem{proposition}[theorem]{Proposition}

\newtheorem{remark}{Remark}[section]

\newtheorem{definition}{Definition}[section]

\theoremstyle{definition}
\theoremstyle{remark}
\numberwithin{equation}{section}

\begin{document}
\author{S.A. Argyros and V. Felouzis}
\title{Interpolating Hereditarily Indecomposable Banach Spaces}
\date{December, 1997}

\begin{abstract}
It is shown that every Banach space either contains $\ell ^1$ or it has an
infinite dimensional closed subspace which is a quotient of a H.I. Banach
space.Further on, $L^p(\lambda )$, $1<p<\infty $,  is a quotient of a H.I
Banach space.
\end{abstract}
\maketitle

\textbf{Introduction} .A Banach space $X$ is said to be Hereditarily
Indecomposable (H.I.) if for every pair of closed subspaces $Y$, $Z$ of
$X$ with $%
Y\cap Z=\left\{ 0\right\} $, $Y+Z$ is not a closed subspace.(By
``subspace'', in the sequel, we mean closed infinite dimensional subspace of 
$X$ ). The H.I spaces are a new and, as we believe, a fundamental class of
Banach spaces. The celebrated example of a Banach space with no
unconditional basic sequence, due to W. Gowers and B. Maurey (\cite{GM}), is
the first construction of a H.I. space. It is easy to see that every H.I.
space does not contain any unconditional basic sequence. Actually, the
concept of H.I. spaces came after W.Johnson's observation that this was a
property of Gowers - Maurey example. To describe even further the peculiar
structure of a H.I. space, we recall an alternative definition of such a
space. So, a Banach space $X$ is a H.I. space if and only if for every pair
of subspaces $Y$, $Z$ and $\varepsilon >0$ there exist $y\in Y$, $z\in Z$
with $||y||=||z||=1$ and $||y-z||<\varepsilon $. Thus, H.I. spaces are
structurally irrelevant to classical Banach spaces, in particular to Hilbert
spaces. Other constructions of H.I. spaces already exist. We mention Argyros
and Deliyanni construction of H.I spaces which are asymptotic $\ell ^1$
spaces (\cite{AD2}), V.Ferenczi's example of a uniformly convex H.I. space (%
\cite{F2}) and H.I. modified asymptotic $\ell ^1$ spaces contained in \cite
{ADKM}.

The construction of such a space is hard and builds upon two fundamental
ideas. The first is Tsirelson's recursive definition of saturated norms (%
\cite{Ts}) and the second is Maurey - Rosenthal's construction of weakly
null sequences without unconditional basic subsequence (\cite{MR}). It is
natural to expect that H.I. spaces share special and interesting properties
not located in the previously known Banach spaces. Indeed, the following
theorem is included in \cite{GM}: Every bounded linear operator from a H.I.
space $X$ to itself is of the form $\lambda I+S$, where $I$ denotes the
identity operator and $S$ is a strictly singular operator. We recall that an
operator is strictly singular if its restriction to any subspace is not an
isomorphism. As a consequence of this theorem, every H.I. space is not
isomorphic to its hyperplanes. Later on, Gowers proved a new dichotomy for
Banach spaces (\cite{G}). The result is that every Banach space either
contains an unconditional basic sequence or has a subspace which is a H.I.
space. We will use this result to prove our dichotomy related to quotients
of H.I. spaces. Also N.Tomczak-Jaegermann solved the distortion problem for
H.I. spaces by showing that every H.I. space is arbitrarily distortable (%
\cite{To}) .

In the present paper we demonstrate that in spite of the fact that the
structure of H.I. spaces is irrelevant to that of spaces with an
unconditional basis, still, there are ways to connect these two classes.
Thus we show that the class of Banach spaces which are quotients of H.I.
spaces is extensive and, further on, large classes of operators between
Banach spaces are factorized through H.I. spaces.

The structure of the quotients of H.I. spaces has also been  studied by
Ferenczi in \cite{F1}, who, showed that every quotient of Gowers - Maurey
space is a H.I. space and, at the same time, there exists a quotient of a
H.I. space which is not a H.I. space. The answer to the question whether
every separable Banach space is a quotient of a H.I. space is negative, due
to the lifting property of $\ell ^1\left( \mathbb{N}\right) $. Indeed, as it is
well known, every Banach space that has $\ell ^1\left( \mathbb{N}\right) $ as a
quotient it has a complemented subspace isomorphic to $\ell ^1\left( \mathbb{N}%
\right) $ and hence is not a H.I. space. The following result shows that $%
\ell ^1\left( \mathbb{N}\right) $ is somehow the only exception. To be more
precise we prove the following dichotomy:

\begin{theorem}
Every Banach space $X$ either contains a subspace isomorphic to $\ell ^1(%
\mathbb{N})$ or it has a subspace which is a quotient of a H.I. Banach space.
\end{theorem}

Further on, every member in the family of the, so called, classical Banach
spaces not containing $\ell ^1(\mathbb{N})$ is actually a quotient of a H.I.
space. For example, we show that separable Hilbert spaces or, even, $%
L^p\left( \lambda \right) $, $1<p<\infty $, $c_0\left( \mathbb{N}\right) $, are
quotients of H.I. Banach spaces. We point out, by standard duality
arguments, that Theorem 0.1. gives that the dual of a H.I. space could
contain isomorphically $L^p\left( \lambda \right) $, $1<p<\infty $, or $\ell
^1\left( \mathbb{N}\right) $.

The proof of these results is quite long and is based upon the following two
ingrendients. The first is a general construction of interpolation H.I.
spaces and the second is the geometric concept of a thin or an $\mathbf{a}$%
-thin norming set. We will say more about these shortly. We now give the
definition of an $\mathbf{a}$-thin set, suggested to us by B. Maurey, in
order to present the factorization results of this paper.

For a null sequence $\mathbf{a=}\left( a_n\right) _{n\in \mathbb{N}}$ of
positive real numbers and a bounded convex symmetric subset $W$ of a Banach
space $X$ we say that $W$ is an $\mathbf{a}$-thin set if the sequence $%
\left( ||\;||_n\right) _n$ of the equivalent norms on $X$, defined by
Minkowski's gauges of the sets $2^nW+a_nB_X$, is not uniformly bounded on
the unit ball of every subspace of $X$.

The second result of this paper concerns factorization of operators and it
is the following:

\begin{theorem}
Let $T:X\rightarrow Y$ be a bounded linear operator between Banach spaces
such that $T[B_X]$ is an $\mathbf{a}$-thin set. Then there exists a H.I.
space $Z$ and bounded linear operators $F_1:X\rightarrow Z,$ $%
F_2:Z\rightarrow Y$ such that $T=F_2\circ F_1$. (i.e. $T$ is factorized
through a H.I. space.)
\end{theorem}

As a consequence of this theorem we show that every $T\in \mathcal{L}\left(
\ell ^p,\ell ^q\right) $, where $p\neq q$, $p,$ $q\in [1,\infty )$ is
factorized through a H.I. space. Additionally, the identity map $I:L^\infty
\left( \lambda \right) \rightarrow L^1\left( \lambda \right) $ is also
factorized through a H.I. space and so is for every strictly singular
operator $T\in \mathcal{L}\left( \ell ^p,\ell ^p\right) .$

The next result of the paper refers to the structure of the H.I. spaces. As
we have already mentioned at the beginning, Gowers and Maurey have shown
that for $X$ a H.I. space, every $T$ in $\mathcal{L}\left( X,X\right) $, is
of the form $\lambda I+S$. It is not known whether a Banach space $X$ such
that every $T$ in $\mathcal{L}\left( X,X\right) $ is of the form $\lambda I+K
$, where $K$ is a compact operator, does exist. However, the construction of
strictly singular but not compact operators on a H.I. space does not seem
easy. What is already known is a construction, due to Gowers, of an operator 
$T$ from a subspace $Y$ of Gowers - Maurey space $X$ to the whole space
which is strictly singular and not compact. It is shown here that H.I.
spaces with many strictly singular non compact operators do exist. More
precisely we have the following result:

\begin{theorem}
There exists a H.I. space $X$ with the property that for every $Y$ subspace
of $X$ there exists a strictly singular non-compact operator $T\in \mathcal{L%
}\left( X,X\right) $ with the range of $T$ contained in $Y$.
\end{theorem}

The last result of the paper that we would like to mention concerns $\ell
^p- $saturated Banach spaces. We recall that a Banach space $X$ is said to
be $\ell ^p$-saturated, for $p\in [1,\infty )$, if every subspace of $X$
contains a further subspace which is isomorphic to $\ell ^p$. (The class of $%
c_0$-saturated Banach spaces is similarly defined).

\begin{theorem}
Every reflexive Banach space $X$ with an unconditional basis contains a
subspace $Y$ satisfying the property that for every $p\in [1,\infty )$,
there exists an $\ell _p$-saturated space $Z$ which has as quotient the
space $Y$. The same holds for some $c_0$-saturated space $Z$.
\end{theorem}

We also show that when $X$ is equal to some $L^q(\lambda )$, $1<q<\infty $,
the subspace $Y$ coincides with the whole space. Probably, the most
interesting case of the above theorem is that of the quotients of $c_0$%
-saturated spaces and we mention that D.H. Leung (\cite{Le}) has shown that
separable Hilbert spaces are quotients of $c_0$-saturated spaces. Our
methods are different from his.

We now move on to describe how the paper is organized and also to explain
briefly the basic ideas of the proofs of the theorems mentioned above.

As is indicated in the title of the paper, the general scheme that we follow
is interpolation methods. It is convenient for us to use the Davis, Figiel,
Johnson and Pelczynski's method (\cite{DFJP}) which we extend for cases
where the external norm is not necessarily unconditional. Thus in the first
section we introduce the $d$-product norm of the space $\Omega _{00}=\left(
\prod\limits_{n=1}^\infty X_n\right) _{00}$, where $(X_n)_n$ is a sequence
of Banach spaces. Then for a sequence $(||\;||_n)_n$ of equivalent norms on
a given Banach space $X$ and a $d$-product norm on $\Omega _{00}=\left(
\prod\limits_{n=1}^\infty (X,||\;||_n)\right) _{00}$ we consider the
diagonal space $\Delta \tilde{X}$, in the completion of $\Omega _{00}$,
consisting of the vectors of the form $\overline{x}=(x,x,...)$ such that $%
x\in X$ and $\overline{x}$ belongs to the completion of $\Omega _{00}$. That
describes the interpolation spaces that we are going to use.

In the second section we define and study the key notions of thin and $%
\mathbf{a}$-thin sets. The concept of a thin set was introduced in Neidinger
Ph.D. thesis (\cite{N}) and it is defined as follows.

A convex bounded symmetric set $W$ in a space $X$ is called a thin set if
there exists $\varepsilon >0$ such that for every subspace $Z$ of $X$ and
every real $\lambda $ we have that $B_Z$ is not contained in the set $%
\lambda W+\varepsilon B_X$. In a previous version of the paper we
exclusively used the notion of a thin set. Then B. Maurey made the important
observation that everything works if we use the concept of $\mathbf{a}$-thin
set, introduced above, instead of thin set. To explain the difference
between these two notions we notice that if $W$ is a thin set then it is an $%
\mathbf{a}$-thin set for every null sequence $\mathbf{a}=(a_n)_n$ of
positive real numbers. We would like to point out that the use of $\mathbf{a-%
}$thin sets has many advantages which allow us to arrive to the present form
of the paper. Further on, B. Maurey showed us how to prove that $B_{L^\infty
\left( \lambda \right) }$ is an $\mathbf{a}$-thin set for an appropriate
sequence $\mathbf{a}$. Notice that Rosenthal had observed that $B_{L^\infty
\left( \lambda \right) }$ is not a thin set (\cite{N},\cite{N2}). Therefore $%
B_{L^\infty \left( \lambda \right) }$ separates these two classes. We use
the $\mathbf{a}$-thin sets in the following manner: Start with a set $W$
which for a sequence $\mathbf{a}=(a_n)_n$ is an $\mathbf{a}$-thin set.
Denote also by $||\;||_n$ the equivalent norm defined by the set $2^nW+a_nB_X
$ and suppose that there exists a $d$-product norm on $\Omega _{00}=\left(
\prod\limits_{n=1}^\infty (X,||\;||_n)\right) _{00}$ which is block H.I.
(Definition 2.3).Then the diagonal space $\Delta \tilde{X}$ is a H.I. space.
Therefore the use of thin or $\mathbf{a}$-thin sets is the essential tool to
construct H.I. interpolation spaces.

The next three sections are devoted to the construction of $\mathbf{a}$-thin
(or thin) norming sets. We recall that a bounded set $W$ in a Banach space
norms a subspace $Y$ of $X^{*\text{ }}$if there exists a constant $C>0$ such
that for every $x^{*}$ in $Y$, $||x^{*}||\leq C\sup \left\{ |x^{*}\left(
w\right) |:w\in W\right\} $. Our goal is, given $A$, a member of a certain
class of separable Banach spaces, to construct a corresponding space $X_A$
and a set $W$, subset of $X_A$, which is $\mathbf{a}$-thin (or thin) and
norms a subspace $Y$ of $X_A^{*}$ isometric and $w^{*}$-homeomorphic to $%
A^{*}$. It does not seem obvious how to construct such a $W$ even for
concrete Banach spaces. However there were indications in the literature for
the existence of such sets. We mention a result due to J.Bourgain in
Radon-Nikodym theory that describes a similar phenomenon, from which we got
the initial idea of constructing such sets. It is proven in \cite{Bo2} that
if $K$ is a bounded convex set in a Banach space and $\varepsilon >0$ such
that the convex combinations of slices of $K$ have diameter greater than $%
\varepsilon $, then the set $K$ norms a subspace of $X^{*\text{ }}$
isomorphic to $\ell ^1\left( \mathbb{N}\right) $. It was somehow expected that
one can construct such sets $K$, with the additional property to be thin
sets. Our construction concerns spaces $A$ with an unconditional basis. For
such a space $A$ we consider an auxiliary space $X_A$ which is of the form $%
\left( \sum \bigoplus \ell ^1\left( k_n\right) \right) _A$ and we ``spread''
the positive part of the unit ball of $A$ on the branches of an appropriate
tree in the ball of $X_A$. Thus after defining such a set $K$, the set $W$
is the set $\overline{co}\left( K\cup -K\right) .$ It easily follows as a
result of this construction that $W$ norms a subspace of $X_A^{*}$ isometric
to $A^{*}$. The difficult part is to show that $W$ is an $\mathbf{a}$-thin
set. This is not always true. For example if $A$ contains $\ell ^1\left( 
\mathbb{N}\right) $ then $W$ is not an $\mathbf{a}$-thin set.

In the third section of the paper we prove that $W$ is an $\mathbf{a}$- thin
set for some special cases of reflexive Banach spaces with an unconditional
basis. The case of classical spaces $L^p\left( \lambda \right) $, $%
1<p<\infty $, is included.

In the fourth section we show that every reflexive space $A$ with an
unconditional basis contains a block subspace $B$ such that the set $W$ in $%
X_B$ is an $\mathbf{a}$-thin set. This is the most general result that
follows for our constructions. The proof depends on some earlier results
from \cite{AMT} on the summability of weakly null sequences.

The fifth section contains the proof that the set $W$, defined for the space 
$c_0\left( \mathbb{N}\right) $, is a thin set.

To explain the role of $\mathbf{a}$-thin (or thin ) norming sets let us
assume that for a sequence $\left( X_n\right) _{n\in \mathbb{N}}$ of separable
Banach spaces there exists a $d-$product norm on $\Omega _{00}=\left(
\prod\limits_{n=1}^\infty X_n\right) _{00}$ which is block-H.I.. In that
case, as we mentioned earlier, the diagonal space $\Delta \tilde{X}_A$
defined by the $\mathbf{a}$-thin set $W$ in $X_A$ and the block H.I. $d$%
-product norm is a H.I. space. Further on, since $W$ norms a subspace of $%
X_A^{*}$ isometric and $w^{*}$-homeomorphic to $A^{*}$ we get that $A^{*}$
is isomorphic and $w^{*}$-homeomorphic to a subspace of $\left( \Delta 
\tilde{X}_A\right) ^{*}$. This implies that $A$ is a quotient of $\Delta 
\tilde{X}_A$.

The sixth section contains the general construction of block-H.I. $d-$%
product norms. In particular we prove that for every sequence $\left(
X_n\right) _{n\in \mathbb{N}}$ of separable Banach spaces there exists a $d$%
-product norm on the space $\Omega _{00}=\left( \prod_{n=1}^\infty
X_n\right) _{00}$ which is block-H.I.. This means that for every block
normalized sequence $\left( x_n\right) _{n\in \mathbb{N}}$ in $\Omega _{00}$
the Banach space $\overline{span}[\left( x_n\right) _{n\in \mathbb{N}}]$ is a
H.I. space. We follow the scheme presented in the construction of H.I.
asymptotic $\ell ^1$ spaces in \cite{AD2}. Here we use a simpler sequence of
compact families instead of the Schreier families and hence the proofs are
not so complicated. For somebody experienced in such constructions it will
be clear that it is possible to use either the Gowers-Maurey scheme or the
asymptotic $\ell ^1$ scheme and to obtain variations of block-H.I. $d$%
-product norms. The choice of the particular one has been made because it
explains in a simpler setting the approach through the analysis of
functionals introduced in \cite{AD2}.

In the seventh section the final results, including the results mentioned
above and some problems related to this work, are presented.

We expect our results to contribute towards the better understanding of the
structure of a H.I. space and possibly help solve certain problems of the
general theory.

\textbf{Acknowledgments:}We would like to thank Bernard Maurey for his
valuable help which allowed us to make an essential step forward and to
arrive at the present form of the paper. We also thank Nicole
Tomczak-Jaegermann for the discussions we had during her visit at Herakleion
and Athens. Her suggestion to use quantities similar to Definition \ref
{Def4.01} was really important to us. Finally, we thank Antonis Manoussakis,
Irene Deliyanni and Apostolos Giannopoulos for their help during the
preparation of the paper.

\section{\textbf{DIRECT PRODUCTS-DIAGONAL SPACES.}}

In this section we present a generalization of the Davis-Figiel-
Johnson-Pelczynski's interpolation method, by extending this to cases where
the connecting external norm is not necessarily unconditional. We will use
this in the next sections to construct interpolation spaces which are
Hereditarily Indecomposable.

\textbf{Notation. }Let $(X_n,||\,\,||_n)_{n\in \mathbb{N}}$ be a sequence of
Banach spaces and $\Omega =\overset{\infty }{\underset{\,n=1}{\Pi }}X_n$
their Cartesian product.\newline
The \textbf{support }of a vector, denoted by $supp(\tilde{x})$, is the set
of all $n\in \mathbb{N}$ such that $\tilde{x}(n)\neq 0$. $\Omega _{\,00}=\{%
\tilde{x}\in \Omega \,:supp(\tilde{x})$is a finite set$\}$.\newline
The \textbf{range} of a vector $\tilde{x}\in \Omega $, written $rang\tilde{x}
$ is the interval of integers $[\min $\textit{supp}$\tilde{x},$ $\max $%
\textit{supp}$\tilde{x}]$.\newline
For any $A\subseteq \mathbb{N}$ we define a linear transformation $P_A:\Omega
\rightarrow \Omega $ by $P_A(\tilde{x})(n)=\left\{ 
\begin{array}{c}
\tilde{x}(n)\,\,\,\text{if \thinspace }n\in A \\ 
0\,\,\,\,\,\,\,\,\,\,\,\,\,\,\,\,\text{if \thinspace }n\notin A
\end{array}
\right. $

$P_n=P_{\{1,...,n\}},$ $\,\,\,\,\,\,P_{n^c}=P_{\{n+1,....\}},\,\,\,\,\,\,\,$ 
$\pi _n=P_{\{n\}}.$\newline
For any $k\in \mathbb{N}$ we denote by $i_k:X_k\rightarrow \Omega _{00}$ the
natural embedding defined by

$i_k(x)(n)=\left\{ 
\begin{array}{c}
x\,\,\,\,\,\text{if\thinspace \thinspace \thinspace \thinspace \thinspace }%
n=k \\ 
0\,\,\,\,\,\text{if\thinspace \thinspace \thinspace \thinspace }n\neq k
\end{array}
\right. $

\begin{definition}
A Banach space $(\tilde{X},||\,\,||)$ is called {\bf d-product} of the
sequence $(X_n,||\,\,||_n)_{n\in \mathbb{N}}$ if the following conditions are
satisfied.\newline
1.\thinspace \thinspace \thinspace \thinspace \thinspace $\Omega
_{00}\subseteq \tilde{X}\subseteq \Omega $ and \thinspace $\Omega _{00}$ is
a dense subset of $\tilde{X}$.\newline
2.\thinspace \thinspace \thinspace \thinspace $||x||_n=||i_n(x)||$, for
every $n\in \mathbb{N}$ and $x\in X_n$.\newline
3.\thinspace \thinspace \thinspace \thinspace \thinspace $P_n:\tilde{X}%
\rightarrow \tilde{X}$ is a bounded transformation for any $n\in \mathbb{N}$
and for any $\tilde{x}\in \tilde{X},\,\,\,\tilde{x}=\underset{n\rightarrow
\infty }{lim}P_n\tilde{x}$.
\end{definition}

If $(X_n,\,||\,\,||_n)_{n\in \mathbb{N}}$ is a sequence of Banach spaces and $A$
is a Banach space with 1-unconditional basis $(e_n)_{n\in \mathbb{N}}$ then%
\newline
$\left( \underset{n=1}{\overset{\infty }{\prod }}X_n\right) _A=\left\{
(x_n)_{n\in \mathbb{N}}\in \underset{n=1}{\overset{\infty }{\prod }}X_n:\left\| 
\underset{n=1}{\overset{\infty }{\sum }}||x_n||_ne_n\right\| _A<\infty
\right\} $ with norm\newline
$||(x_n)_{n\in \mathbb{N}}||=\left\| \underset{n=1}{\overset{\infty }{\sum }}%
||x_n||_ne_n\right\| _A$ \thinspace \thinspace is an example of d-product of
the sequence $(X_n)_{n\in \mathbb{N}}$.

By the uniform boundedness principle there exists a $C>0$ such that $||P_n%
\tilde{x}||\leq C||\tilde{x}||$, for any $\tilde{x}\in \tilde{X}$ and any $%
n\in \mathbb{N}$. Therefore the sequence $(X_n)_{n\in \mathbb{N}}$ is a Schauder
decomposition of the space $\tilde{X}$ (for a detailed study of this notion
see \cite{LT} p. 47-52)\newline
The family $(P_n)_{n=1}^\infty $satisfies the following properties: \newline
$\mathbf{(P1)}$ $P_n\circ P_m=P_{\min \{m,n\}}$,\newline
$\mathbf{(P2)\;}$ $\underset{n\in \mathbf{N}}{\sup }||P_n||<\infty $, 
\newline
$\mathbf{(P3)}$ \ \ For every $x\in X$ $x=\underset{n\rightarrow \infty }{%
\lim }P_nx$ or equivalently,\newline
$\mathbf{(P3)}^{\prime }$ $\;X=\overline{\underset{n=1}{\overset{\infty }{%
\bigcup }}P_n(X)}$.\newline
Conversely if in a Banach space $X$ a family $(P_n:X\rightarrow
X)_{n=1}^\infty $ of linear operators satisfies $\mathbf{(P1),(P2),(P3)}$
the space $X$ is the d-product of the family $\left( (P_n-P_{n-1})(X)\right)
_{n=1}^\infty $ where $P_0=0$.\newline
The norm $||\,\,||$ of the space $\tilde{X}$ is called:\newline
(a) \thinspace \thinspace {\bf C-bimonotone} if there is a $C>0$ such
that $||P_A(\tilde{x})||\leq C||\tilde{x}||$ for any interval $A\subseteq 
\mathbb{N}$ and any $\tilde{x}\in \tilde{X}$.\newline
(b)\thinspace \thinspace \thinspace \thinspace {\bf unconditional} if
there is a $C>0$ such that $||P_A(\tilde{x})||\leq C||\tilde{x}||$ for any $%
A\subseteq \mathbb{N}$ and any $\tilde{x}\in \tilde{X}$.\newline
(c)\thinspace \thinspace \thinspace \thinspace {\bf boundedly complete }%
if for any $\tilde{x}\in \Omega $,\thinspace \thinspace if $\underset{n\in 
\mathbf{N}}{\sup }||P_n\tilde{x}||<\infty $ then $\tilde{x}\in \tilde{X}$.%
\newline
(d)\thinspace \thinspace \thinspace \thinspace {\bf shrinking} if for any 
$\tilde{x}^{*}\in \tilde{X}^{*}$, $\tilde{x}^{*}=\underset{n\rightarrow
\infty }{\lim }\,\,\,\,\,P_n^{*}\tilde{x}^{*}$, where

\begin{center}
$\,P_n^{*}\tilde{x}^{*}(x)=\underset{i=1}{\overset{n}{\sum }}\tilde{x}%
^{*}(\pi _i(\tilde{x}))=\tilde{x}^{*}(P_n\tilde{x})$
\end{center}

Since $(P_n^{*})_{n=1}^\infty $ satisfies $\mathbf{(P1)}$ and $\mathbf{(P2)}$
the norm of $\tilde{X}$ is shrinking if and only if $\tilde{X}^{*}$ is a
d-product of the sequence $(X_n^{*})_{n=1}^\infty $.\newline
In the sequel we shall assume that $||\,\,||$ is always 1-bimonote and
therefore\newline
\begin{equation*}
||\tilde{x}||_\infty \leq ||\tilde{x}||\leq ||\tilde{x}||_1
\end{equation*}
where , $||\tilde{x}||_\infty =\underset{n\in \mathbf{N}}{\sup }||\tilde{x}%
(n)||_n,\,\,||\tilde{x}||_1=\underset{n=1}{\overset{\infty }{\sum }}||\tilde{%
x}(n)||_n.$\newline
In general $\;\;\;\;\;\;\;\;\;\;\;\;\;\;\,||\tilde{x}||_{\ell _p}=\left( 
\underset{n=1}{\overset{\infty }{\sum }}||\tilde{x}(n)||^p\right) ^{\frac
1p}.$

\begin{definition}
Let $(X_n,||\,\,||_n)_{n\in \mathbb{N}}$ be a sequence of Banach spaces such
that $X_n=X_1$ and $||\,\,||_n$ is equivalent to $\|\,\,\|_1$ for any $n\in 
\mathbb{N}$. Let $\tilde{X}$ be any d-product of the sequence $%
(X_n,||\,\,||_n)_{n\in \mathbb{N}}$.\newline
The {\bf diagonal space} $\Delta \tilde{X}$ of $\tilde{X}$ is the
(closed) subspace of $\tilde{X}$ consisting of all $\tilde{x}\in \tilde{X}$
such that $\tilde{x}(n)=\tilde{x}(1)$, for any $n\in \mathbb{N}$.\newline
Define $J:\Delta \tilde{X}\rightarrow X_1,\,\,\,I:\Delta \tilde{X}%
\rightarrow \tilde{X},\,\,\pi _1:\tilde{X}\rightarrow X_1$ by $J(\tilde{x})=%
\tilde{x}(1),\,\,I(\tilde{x})=\tilde{x},\,\,\pi _1(\tilde{x})=\tilde{x}(1)$.
\end{definition}

Then $J=\pi _1\circ I$ and $J$ is a $1-1$, linear, continuous transformation.%
\newline
The operator $J$ is not always an isomorphism, i.e. $J(\Delta \tilde{X})$ is
not necessarily a closed subspace of $X_1$.\newline
If the norm of $\tilde{X}$ is boundedly complete and shrinking then $J$ has
an interesting property, it is a Tauberian operator.

\begin{definition}
\cite{KW}  An operator $T:X\rightarrow Y$ is a {\bf Tauberian operator }%
iff $(T^{**})^{-1}(Y)\subseteq X$
\end{definition}

(for a study and characterizations of Tauberian operators see \cite{N}, \cite
{N2}, \cite{NR}).

The dual $\tilde{X}^{*}\,\,$of $\tilde{X}$ is identified with a subspace of $%
\underset{n=1}{\overset{\infty }{\prod }}X_n^{*}$ by the $1-1$ linear map $%
\tau ^{*}:\tilde{X}^{*}\rightarrow \underset{n=1}{\overset{\infty }{\prod }}%
X_n^{*}$ \thinspace \thinspace \thinspace \thinspace \thinspace \thinspace
\thinspace \thinspace where$\;\;\;$%
\begin{equation*}
\tau ^{*}(\tilde{x}^{*})=(x_n^{*})_{n\in \mathbb{N}}
\end{equation*}
\thinspace \thinspace \thinspace \thinspace \thinspace \thinspace \thinspace
\thinspace \thinspace \thinspace \thinspace and\thinspace \thinspace
\thinspace \thinspace \thinspace \thinspace \thinspace \thinspace $%
x_n^{*}(x)=\tilde{x}^{*}(i_n(x)),\,\,\,$for any $x\in X_n$.

The space $\tilde{X}^{*}$ is not always a $d$-product of the sequence $%
(X_n^{*})_{n\in \mathbb{N}}$. This happens if and only if the norm of $\tilde{X}
$ is shrinking.\newline
We also define $\tau ^{**}:\tilde{X}^{**}\rightarrow \underset{n=1}{\overset{%
\infty }{\prod }}X_n^{**}$ by

\begin{center}
$\tau ^{**}(\tilde{x}^{**})=(x_n^{**})_{n\in \mathbb{N}}$
\end{center}

where $x_n^{**}(x^{*})=\tilde{x}^{**}(i_n^{*}x^{*})$, for any $x^{*}\in
X_n^{*}$ and $n\in \mathbb{N}$.

The map $\tau ^{**}$ is linear but not always $1-1$. The latter condition is
satisfied if the norm $||\,\,||$ of $\tilde{X}$ is shrinking. We notice that
the space $\tilde{X}$ is reflexive if and only if each $X_n$ is reflexive
and the norm of $\tilde{X}$ is boundedly complete and shrinking.(The proof
is similar to that of the corresponding theorem of James for spaces with
Schauder basis \cite{LT})

\begin{proposition}
Assume that the norm of $\tilde{X}$ is shrinking and boundedly complete then 
$J:\Delta \tilde{X}\rightarrow X_1$ is a $1-1$ Tauberian operator.
\end{proposition}

\proof%
We denote by $J_n:\Delta \tilde{X}\rightarrow X_n$ the operator $\pi _n\circ
I$. Since $\tilde{X}$ is shrinking each $\tilde{x}^{*}$ in $\tilde{X}^{*}$
is represented as $(x_n^{*})_{n\in \mathbb{N}}$ and hence

\begin{center}
\begin{equation}
I^{*}(x_n^{*})_{n\in \mathbb{N}}=(J_n^{*}x_n^{*})_{n\in \mathbb{N}}  \tag{1}
\label{1}
\end{equation}
\end{center}

Also, since $\tilde{X}$ is boundedly complete for any $x^{**}$ in $(\Delta 
\tilde{X})^{**}$ we have by ($\ref{1})$

\begin{center}
\begin{equation}
I^{**}(\tilde{x}^{**})=(J_n^{**}x^{**})_{n\in \mathbb{N}}  \tag{2}  \label{2}
\end{equation}
\end{center}

$I^{**}$ is an isometry, therefore $J^{**}$ is $1-1$.\newline
Suppose now that $J^{**}(\tilde{x}^{**})\in X_1$ then $(J_n^{**}x^{**})_{n%
\in \mathbb{N}}\in \tilde{X}$ and by (\ref{2}), $I^{**}(\tilde{x}^{**})\in
\Delta \tilde{X}$, hence $\tilde{x}^{**}\in \Delta \tilde{X}$.\newline
Therefore $(J^{**})^{-1}(X_1)\subseteq \Delta \tilde{X}$.$\,\,\,$%
\endproof%

\begin{definition}
Let $W$ be a bounded convex symmetric non-empty subset of a Banach space $%
(X,|||\,|||)$, $\mathbf{a}=(a_n)_{n=1}^\infty $ a null sequence of positive
numbers and $W_n=2^nW+a_nB_X$. \newline
Let $||\,||_n$ be the Minkowski functional of $W_n$ which is a norm
equivalent of $\,|||\,|||$, $(X,||\,||_n)_{n\in \mathbb{N}}$ be the above
defined sequence of Banach spaces and $\tilde{X}$ be any d-product of $%
(X,||\,||_n)_{n\in \mathbb{N}}.$ The diagonal space $\Delta \tilde{X}$ is
called the $\mathbf{(a},\tilde{X},W\mathbf{)}${\bf -diagonal space,}or
simply, a{\bf \ diagonal space}.
\end{definition}

\begin{remark}
If $a_n=2^{-n}$ for every $n$ and $\tilde{X}=\left( \underset{n=1}{\overset{%
\infty }{\sum }}\oplus X_n\right) _{\ell ^2}$ then the $\mathbf{(a},\tilde{X}%
,W\mathbf{)}$-diagonal space is the space introduced in \cite{DFJP}.
\end{remark}

\begin{proposition}
\label{A}Let $\Delta \tilde{X}$ be any diagonal space of the sequence 
\newline
$(X,||\,||_n)_{n\in \mathbb{N}}$. Then the following hold:\newline
(a) $W\subseteq JB_{\Delta \tilde{X}}$\newline
(b) For every $\varepsilon >0$ there exists $\lambda >0$ such that\newline
\begin{equation*}
JB_{\Delta \tilde{X}}\subseteq \lambda W+\varepsilon B_X.
\end{equation*}
\end{proposition}

\proof%
: (a) If $w\in W$ then $||w||_n\leq \frac 1{2^n}$ and therefore $%
||(w,w,...)||\leq \underset{n=1}{\overset{\infty }{\sum }}\frac 1{2^n}\leq 1$%
\newline
(b) If $x\in JB_{\Delta \tilde{X}}$ then $||x||_n\leq 1$ and so $x\in 
\overline{2^nW+a_nB_X}$ for every $n\in \mathbb{N}$. Choose $n$ such that $%
a_n<\varepsilon $. Then $x\in 2^nW+\varepsilon B_X$ so $JB_{\Delta \tilde{X}%
}\subseteq 2^nW+\varepsilon B_X$.

\begin{definition}
Let $(X,||\,||)$ be a Banach space, $A,B\subseteq X$ and $\varepsilon >0$.%
\newline
We say that\newline
(i)The set $A$ $\mathbf{\varepsilon }${\bf -absorbs} $B$ if there exists
a $\lambda >0$ such that $B\subseteq \lambda A+\varepsilon B_X$.\newline
(ii)The set $A$ {\bf almost absorbs }$B$ if $A$ $\varepsilon $-absorbs $B$
for any $\varepsilon >0$.
\end{definition}

The following Lemma is due to A.Grothendieck (\cite{Gr}) and for a proof see 
\cite{D}

\begin{lemma}
\label{Gr}Let $K$ be a weakly closed subset of the Banach space $X$.\newline
Suppose for each $\varepsilon >0$ there exists a weakly compact set $%
K_{\varepsilon _{}}$ in $X$ such that $K\subset K_{\varepsilon
_{}}+\varepsilon B_X$\newline
Then $K$ is weakly compact.%
\endproof%
\end{lemma}

\begin{proposition}
Let $\Delta \tilde{X}$ be a $(\tilde{X},W)$-diagonal space such that the
norm of $\tilde{X}$ is boundedly complete and shrinking. Then $\Delta \tilde{%
X}$ is reflexive if and only if $W$ is relatively weakly compact subset of $X
$.
\end{proposition}

\proof%
. If $\Delta \tilde{X}$ is reflexive then $W$ is a relatively weakly compact
subset of $J\Delta \tilde{X}$ since $J^{-1}(W)$ is a closed convex subset of 
$B_{\Delta \tilde{X}}$.\newline
For the converse we observe that by Proposition \ref{A} (b) $W$ almost
absorbs the set $JB_{\Delta \tilde{X}}$ and hence from Lemma \ref{Gr} $%
JB_{\Delta \tilde{X}}$ is relatively weakly compact set . It follows that $%
\overline{JB_{\Delta \tilde{X}}}^{w^{*}}\subseteq B_X$ hence $J^{**}(\Delta 
\tilde{X}^{**})\subseteq X$. Since $J$ is a Tauberian operator we get that $%
\Delta \tilde{X}=\Delta \tilde{X}^{**}$. 
\endproof%

\section{{\bf THIN SETS}}

In this section we start the study of thin and $\mathbf{a}$-thin sets that
we will continue in the next three sections. The notion of a thin set was
introduced in \cite{N},\cite{N2} and that of $\mathbf{a-}$thin set ,which is
weaker, was suggested to us by B.Maurey who also showed us how to prove that
the unit ball of $L^\infty \left( \lambda \right) $ is an $\mathbf{a}$-thin
subset of $L^1\left( \lambda \right) $. It was observed by Rosenthal (\cite
{N},\cite{N2}) that this set is not a thin subset of $L^1\left( \lambda
\right) .$\newline
Thin or $\mathbf{a}$-thin sets are important for our results since the
interpolation diagonal space defined by an $\mathbf{a}$-thin set (and hence
by a thin set) and an appropriate external norm is, as we show in this
section (Proposition \ref{Prop2.3}), an H.I-space.

\begin{definition}
a) A bounded convex symmetric non-empty subset $W$ of a Banach space $X$ is
said to be a {\bf thin set} if $W$ does not almost-absorb the unit ball $%
B_Y$ of any infinite dimensional closed subspace $Y$ of $X$. (i.e. for any $Y
$, infinite dimensional closed subspace of $X$, there exists $\varepsilon >0$
such that for all $\lambda >0$, $B_Y\nsubseteq \lambda W+\varepsilon B_X$).%
\newline
b) An operator $T:Z\rightarrow X$ is called {\bf thin} if $TB_Z$ is a
thin subset of $X$.
\end{definition}

\begin{remark}
\label{Rem2.1}(a) There are several examples of thin sets . The closed
convex hull of the basis $\left( e_n\right) _{n=1}^\infty $ in any $\ell
^p\left( \mathbb{N}\right) $,for $1<p<\infty ,$ or $c_0\left( \mathbb{N}\right) $
is a thin set.\newline
(b) Every thin operator is always a strictly singular operator. The converse
is not true. Indeed,as it is shown in \cite{N} and we mentioned above the
set $B_{L^\infty \left( \lambda \right) }$ in $L^1\left( \lambda \right) $
is not a thin set and the identity operator $I:L^\infty \left( \lambda
\right) \rightarrow L^1\left( \lambda \right) $ is a strictly singular
operator.\newline
(c) It follows easily from the definition that if $W$ is a thin subset in a
Banach space $X$ and $\mathbf{a=}\left( a_n\right) _{n\in \mathbb{N}}$ is a
null sequence of positive numbers then the sequence of equivalent norms$%
\left( ||\;||_n\right) _{n\in \mathbb{N}}$ on $X$ defined by the Minkowski
gauges $\left( 2^nW+a_nB_X\right) _{n\in \mathbb{N}}$ are not uniformly bounded
on the $B_Z$ , for every $Z$ infinite dimensional closed subspace of $X$ .
\end{remark}

\begin{definition}
a) Let $\mathbf{a}=(a_n)_{n=1}^\infty $ be a null sequence of positive
numbers. $W$ is said to be a $\mathbf{a}${\bf -thin} subset of $X$ if for
any $Y$, infinite dimensional closed subspace of $X$, $\underset{n\in %
\mathbf{N}}{\sup }\left\{ \underset{y\in B_Y}{\sup }||y||_n\right\} =\infty $%
, where $||$ $||_n$ is the gauge of $W_n=2^nW+a_nB_X$.\newline
b) An operator $T:Z\rightarrow X$ is called $\mathbf{a}${\bf -thin} if $%
TB_Z$ is an $\mathbf{a}$-thin subset of $X$.
\end{definition}

\begin{remark}
Every thin subset of $X$ is $\mathbf{a}$-thin for every null positive
sequence $\mathbf{a}$ and every thin operator is $\mathbf{a}$-thin for every
null positive sequence $\mathbf{a}$. An $\mathbf{a}$-thin operator is
strictly singular. In particular if $W$ is a $\mathbf{a}$-thin set then the
operator $J$ is strictly singular.\newline
The two notions are not equivalent. More precisely we will show \newline
(Lemma\ref{M1},Proposition \ref{M} ) that the unit ball $B_{L^\infty \left(
[0,1]\right) }$ of $L^\infty \left( \left[ 0,1\right] \right) $ is an $%
\mathbf{a-}$thin subset of $L^1\left( \left[ 0,1\right] \right) $ for an
appropriate positive null sequence $\mathbf{a}$. Therefore ,as follows from
Remark \ref{Rem2.1} (b) , the identity operator $I:L^\infty \left( \left[
0,1\right] \right) \rightarrow L^1\left( \left[ 0,1\right] \right) $ is
strictly singular but not thin.\newline
A natural question that arises is if a strictly singular operator is $%
\mathbf{a}$-thin for some null positive sequence $\mathbf{a.}$
\end{remark}

\begin{lemma}
\label{M1}Let $Z$ be an infinite dimensional subspace of $L_2[0,1]$ such
that $Z\subset L_p[0,1]$ for every $p<\infty $ and let 
\begin{equation*}
C_p(Z)=\sup \left\{ \left\| f\right\| _p:f\in Z\text{ and }\left\| f\right\|
_2\leq 1\right\} 
\end{equation*}
then 
\begin{equation*}
\underset{p\rightarrow +\infty }{\liminf }\frac{C_p(Z)}{\sqrt{p}}\geq
e^{-\frac 12}
\end{equation*}
\end{lemma}

\proof%
Let us choose an orthonormal basis $\mathbf{f}=(f_1,...,f_n,...)$ for our
subspace $Z$ of $L_2[0,1]$ and let us fix $p>2$.By assumption,we have $%
Z=span[f_j]_{j=1}^\infty \subset L_p[0,1]$ and $C_p(Z)$ is equal to the
smallest constant $K$ such that 
\begin{equation*}
\left\| \underset{j=1}{\overset{n}{\sum }}c_jf_j\right\| _p\leq K\left\| 
\underset{j=1}{\overset{n}{\sum }}c_jf_j\right\| _2=K\left( \underset{j=1}{%
\overset{n}{\sum }}\left| c_j\right| ^2\right) ^{\frac 12}
\end{equation*}
\newline
for all $n\geq 1$ and for all real scalars $(c_j)$.\newline
For a subspace $G$ of $L_2(\Omega ,P)$spanned by a Gaussian system $\mathbf{g%
}=(g_1,...g_n,...)$ of independent random variables with the common
distribution $(2\pi )^{-\frac 12}e^{-x^2}dx$ ,we know that all linear
combinations $\sum_jc_jg_j$ such that $\sum_j\left| c_j\right| ^2=1$ have
the same distribution ,namely that of $g_1,$therefore 
\begin{equation}
\left\| \sum_{j=1}^nc_jg_j\right\| _p=C_p(G)\underset{j=1}{\overset{n}{\sum }%
}\left| c_j\right| ^2
\end{equation}
\newline
for all real scalars and 
\begin{eqnarray}
C_p(G) &=&\left\| g_1\right\| _p=\left( \int\limits_{-\infty }^{+\infty
}\left| x\right| ^pe^{-\frac{x^2}2}\frac{dx}{2\pi }\right) ^{\frac 1p}= \\
&=&\left( \frac{2^{\frac p2}}{\sqrt{\pi }}\Gamma \left( \frac{p+1}2\right)
\right) ^{\frac 1p}\sim \sqrt{\frac pe}
\end{eqnarray}
\newline
as $p\rightarrow \infty ,$using Stirling's formula..\newline
Let us consider $\overset{n}{\sum\limits_{j=1}}f_j(s)g_j(t).$By (1) we get
that 
\begin{equation*}
\int \left| \overset{n}{\sum\limits_{j=1}}f_j(s)g_j(t)\right| ^pdt=\left(
C_p(G)\right) ^p\left( \underset{j=1}{\overset{n}{\sum }}\left|
f_j(s)\right| ^2\right) ^{\frac p2}
\end{equation*}
\newline
therefore if $p>2$%
\begin{equation*}
\overset{1}{\underset{0}{\int }}\int \left| \overset{n}{\sum\limits_{j=1}}%
f_j(s)g_j(t)\right| ^pds\;dt=\left( C_p(G)\right) ^p\underset{0}{\overset{1}{%
\int }}\left( \underset{j=1}{\overset{n}{\sum }}\left| f_j(s)\right|
^2\right) ^{\frac p2}ds\geq
\end{equation*}
\newline
\begin{equation*}
\geq \left( C_p(G)\right) ^p\left( \underset{0}{\overset{1}{\int }}\left( 
\underset{j=1}{\overset{n}{\sum }}\left| f_j(s)\right| ^2\right) ds\right)
^{\frac p2}=n^{\frac p2}\left( C_p(G)\right) ^p
\end{equation*}
\newline
On the other hand, 
\begin{equation*}
\overset{1}{\underset{0}{\int }}\left| \overset{n}{\sum\limits_{j=1}}%
f_j(s)g_j(t)\right| ^pds\leq \left( C_p(Z)\right) ^p\left( \underset{j=1}{%
\overset{n}{\sum }}\left| g_j(t)\right| ^2\right) ^{\frac p2}
\end{equation*}
\newline
\begin{equation*}
\int \overset{1}{\underset{0}{\int }}\left| \overset{n}{\sum\limits_{j=1}}%
f_j(s)g_j(t)\right| ^pdsdt\leq \left( C_p(Z)\right) ^p\underset{}{\overset{}{%
\int }}\left( \underset{j=1}{\overset{n}{\sum }}\left| g_j(t)\right|
^2\right) ^{\frac p2}dt
\end{equation*}
\newline
and 
\begin{equation*}
\varphi (p,n)=\int \left( \underset{j=1}{\overset{n}{\sum }}\left|
g_j(t)\right| ^2\right) ^{\frac p2}dt=\int\nolimits_{\mathbb{R}^n}\left\|
x\right\| ^pd\gamma _n(x)=\mathcal{\upsilon }_{n-1}\int\limits_0^{+\infty
}r^pe^{-\frac{r^2}2}r^{n-1}\frac{dr}{\left( 2\pi \right) ^{\frac n2}}
\end{equation*}
\newline
where $\mathcal{\upsilon }_{n-1}$ is the surface of the unit sphere in $\mathbb{%
R}^n$.By change of variable , 
\begin{equation*}
\varphi (p,n)=\mathcal{\upsilon }_{n-1}\dfrac{2^{\frac p2+\frac n2-1}}{%
\left( 2\pi \right) ^{\frac n2}}\int\limits_0^{+\infty }u^{\frac p2+\frac
n2-1}e^{-u}du=\mathcal{\upsilon }_{n-1}\frac{2^{\frac p2-1}}{\pi ^{\frac n2}}%
\Gamma \left( \frac{n+p}2\right)
\end{equation*}
\newline
Taking $p=2$ gives 
\begin{equation*}
n=\varphi (2,n)=\frac{\mathcal{\upsilon }_{n-1}}{\pi ^{\frac n2}}\Gamma
\left( \frac n2+1\right)
\end{equation*}
\newline
\begin{equation*}
\frac{\varphi (p,n)}n=\frac{\varphi (p,n)}{\varphi (2,n)}=2^{\frac p2-1}%
\frac{\Gamma \left( \frac n2+\frac p2\right) }{\Gamma \left( \frac
n2+1\right) }
\end{equation*}
\newline
so that\newline
\begin{equation*}
\varphi (p,n)=2^{\frac p2}\frac{\Gamma \left( \frac n2+\frac p2\right) }{%
\Gamma \left( \frac n2\right) }
\end{equation*}
\newline
Summing up our informations, 
\begin{equation*}
n^{\frac 12}C_p(G)\leq C_p(Z)(\varphi (p,n))^{\frac 1p}
\end{equation*}
\newline
or 
\begin{equation*}
C_p(G)\leq \sqrt{\frac 2n}\left( \frac{\Gamma \left( \frac n2+\frac
p2\right) }{\Gamma \left( \frac n2\right) }\right) ^{\frac 1p}C_p(Z)
\end{equation*}
\newline
We know by Stirling's formula that for any fixed $x>0$ we have that $\Gamma
(t+x)\sim t^x\Gamma (t)$ when $t$ tends to $+\infty .$Applying this with $%
x=\frac p2$ we get \newline
\begin{equation*}
\underset{n\rightarrow +\infty }{\lim }\sqrt{\frac 2n}\left( \frac{\Gamma
\left( \frac n2+\frac p2\right) }{\Gamma \left( \frac n2\right) }\right)
^{^{\frac 1p}}=1
\end{equation*}
$\newline
$and it follows that\newline
\begin{equation*}
C_p(G)\leq C_p(Z)
\end{equation*}
\newline
This says that the Gaussian case is the best possible inequality between the 
$L_p$ and $L_2$ norms.By (2) we already know that $C_p(G)$ is equivalent to $%
\sqrt{\frac pe}$ and this concludes the proof.%
\endproof%

\begin{proposition}
\label{M}The unit ball $B_{L_\infty }$ is an $\mathbf{a}-$thin subset of $L_1
$ for $\mathbf{a=}\left( a_n\right) _n$ with $a_n=16^{-n}.$\newline
\proof%
Assume the contrary;then there exists an infinite dimensional subspace $Z$
of $L_1$ and $M>0$ such that 
\begin{equation*}
\frac 1MB_Z\subset 2^nB_{L_\infty }+a_nB_{L_1}
\end{equation*}
\newline
for every $n\in \mathbb{N}\mathbf{.}$\newline
So,every $z\in Z$ with $\left\| z\right\| _1\leq \frac 1M$ has a
decomposition $z=z_0+z_1$ with $\left\| z_0\right\| _{}\leq 2^n$ and $%
\left\| z_1\right\| _1\leq 16^{-n}$.This implies that\newline
$P\left( s\in [0,1]:\left| z(s)\right| >2^{n+1}\right) \leq P\left( s\in
[0,1]:\left| z_1(s)\right| >2^n\right) \leq 2^{-n}a_n$\newline
Let $\lambda >2$ and $n\geq 0$ such that $2^{n+1}<\lambda \leq 2^{n+2}.$We
get that\newline
$P\left( s\in [0,1]:\left| z(s)\right| >\lambda \right) \leq \frac 4\lambda
\exp \left( \frac{-\lambda ^4}{4^4}\right) $\newline
Then\newline
\begin{equation*}
\left\| z\right\| _p^p\leq 2^p+\int\limits_2^{+\infty }p\lambda ^{p-1}\frac
4\lambda \exp \left( \frac{-\lambda ^4}{4^4}\right) d\lambda 
\end{equation*}
\newline
We set $s=\left( \frac \lambda 4\right) ^4\,$and we get that \newline
\begin{equation*}
\left\| z\right\| _p^p\leq 2^p+p4^{p-1}\int\limits_0^{+\infty }s^{\frac{p-1}%
4-1}e^{-s}ds=2^p+p4^{p-1}\Gamma \left( \frac{p-1}4\right) 
\end{equation*}
\newline
Applying Stirling's formula we have that\newline
\begin{equation*}
\left\| z\right\| _p\leq {\large O}\left( {\large p}^{\frac 14}\right) 
\end{equation*}
\newline
when $p\rightarrow +\infty ,$ uniformly for $z$ in the unit ball of $Z$\
contradicting the Lemma \ref{M1}%
\endproof%

{\bf Notation }a) For $A,B$ finite non-empty subsets of $\mathbb{N}$ we
write $A<B$ iff $\max A<\min B$.\newline
b) If $\tilde{x},\tilde{y}\in \left( \underset{n=1}{\overset{\infty }{\prod }%
}X_n\right) _{00}$ we write $\tilde{x}<\tilde{y}$ iff supp$(\tilde{x})<$supp$%
(\tilde{y})$.\newline
c) A sequence $(\tilde{x}_n)_{n\in \mathbb{N}}$ of non-zero vectors of $\left( 
\underset{n=1}{\overset{\infty }{\prod }}X_n\right) _{00}$ is said to be a 
{\bf block sequence} if $\tilde{x}_n<\tilde{x}_{n+1}$ for any $n\in \mathbb{N%
}$.\newline
d) If $\tilde{X}$ is any $d$-product of a sequence $(X_n,||$ $||_n)_{n\in 
\mathbb{N}}$ of Banach spaces then a {\bf block subspace} of $\tilde{X}$ is
the closed linear span of a block sequence $(\tilde{x}_k)_{k\in \mathbb{N}}$.
\end{proposition}

\begin{proposition}
\label{2.3}Let $X$ be a Banach space, $\mathbf{a}=(a_n)_{n=1}^\infty $ a
null sequence of positive numbers and $W$ a $\mathbf{a}$-thin subset of $X$.%
\newline
We denote by $||$ $||_n$ the $n^{\text{th}}$ equivalent norm defined by the
set $2^nW+a_nB_X$ and $Y$ be the $(\mathbf{a},W,\tilde{X})$-diagonal space.%
\newline
Then for any $\varepsilon >0$ and any infinite dimensional closed subspace $Z
$ of $Y$ there exists a block subspace $\tilde{X}_1$ of $\tilde{X}$ and an
infinite dimensional closed subspace $Z_1$ of $Z$ which are $(1+\varepsilon )
$-isomorphic.
\end{proposition}

\proof%
. Given $0<\varepsilon <1$, we say that a normalized vector $\tilde{x}\in 
\tilde{X}$ has $\mathbf{(1-\varepsilon )}${\bf support} the interval $%
E=[n,m]$ if $\max \{||P_n\tilde{x}||,\,||P_{m^c}\tilde{x}||\}<\frac
\varepsilon 4$.\newline
We identify $Y$ with the linear space $\Delta \tilde{X}=\{y\in X:\tilde{y}%
=(y,y,...)\}$.\newline
We start with the following observations: \newline
(1) Any space $(X,\,||\,||_n)$ is isomorphic to $(X,\,||\,||_1)$, so\newline
For any $\varepsilon >0$ and $n\in \mathbb{N}$ there exists $\delta >0$ such
that if $y\in Y$, $||y||_1<\delta $, then $||P_n\tilde{y}||<\varepsilon $,
where $\tilde{y}=(y,y,...)$.\newline
2) Since $W$ is a thin set and $W$ almost absorbs the set $JY$ we get that $%
J $ is a thin operator and hence a strictly singular operator therefore any $%
\delta >0$ and any $Y^{\prime }$, infinite dimensional closed subspace of $Y$%
, there is $y\in Y^{\prime }$ such that $||y||_Y=1$ and $||y||_1<\delta $.%
\newline
Observations (1) and (2) permit us to use a ``traveling hump
argument'' to construct in $Y^{\prime }$ a normalized basic
sequence $(y_i)_{i\in \mathbb{N}}$ and a sequence $(E_i)_{i\in \mathbb{N}}$ of
intervals of $\mathbb{N}$ such that $E_i<E_{i+1}$ and any $y_i$ has $%
(1-\varepsilon _i)$-support the interval $E_i$, where $(\varepsilon
_i)_{i=1}^\infty $ is a given sequence of positive numbers. It is easy to
see that $(y_i)_{i\in \mathbb{N}}$ is $(1+\varepsilon )$ equivalent to $(E_i%
\tilde{y}_i)_{i\in \mathbb{N}}$, provided $\underset{i=1}{\overset{\infty }{%
\sum }}\varepsilon _i<\min \left\{ \frac \varepsilon 4,\frac 14\right\} $. 
\endproof%

\begin{definition}
\label{Def2.3}a) A Banach space $X$ is called a {\bf hereditarily
indecomposable} space or $\mathbf{H.I}${\bf -space} if $X$ has no
subspace that can be decomposed as a topological direct sum of two infinite
dimensional subspaces, or equivalently for any pair of infinite dimensional
subspaces $Y$,$Z$ of $X$ and any $\varepsilon >0$ there exist $y\in Y$, $%
z\in Z$ such that $||y-z||<\varepsilon ||y+z||$.\newline
b) A $d$-product $\tilde{X}$ is called a {\bf block-}$\mathbf{H.I}$ space
if any block subspace of $\tilde{X}$ is an $H.I$ space.
\end{definition}

\begin{proposition}
\label{Prop2.3}Let $\Delta \tilde{X}$ be the $(W,\tilde{X})$-diagonal space
produced by a $\mathbf{a}$-thin set $W$ of a Banach space $X$ and such that $%
\tilde{X}$ is a block-$H.I$ space. Then $\Delta \tilde{X}$ is a $H.I$-space
\end{proposition}

\proof%
Suppose that $\Delta \tilde{X}$ is not a $H.I$-space. Then there exist $%
\tilde{Y}_1,\tilde{Y}_2$ infinite dimensional closed subspaces of $\Delta 
\tilde{X}$ and $\varepsilon >0$ such that $||\tilde{y}-\tilde{z}||\geq
\varepsilon ||\tilde{y}+\tilde{z}||$, for any $\tilde{y}\in \tilde{Y}_1$, $%
\tilde{z}\in \tilde{Y}_2$. We can construct a sequence of vectors $\tilde{y}%
_1,\tilde{z}_1,\tilde{y}_2,\tilde{z}_2,....$ and a sequence of intervals $%
E_1,F_1,E_2,F_2,...$ such that $\tilde{y}_i\in \tilde{Y}_1,\,\,\tilde{z}%
_i\in \tilde{Y}_2,$ $E_1<F_1<E_2<F_2<...$, each $\tilde{y}_i$ has $(1-\frac
1{10^i})$-support the $E_i$ and $\tilde{z}_i$ has $(1-\frac 1{10^i})$%
-support the interval $F_i$.\newline
We set $\tilde{Y}_i^{\prime }=\overline{\left\langle E_i\tilde{y}%
_i\right\rangle _{i=1}^\infty }$, $\tilde{Y}_2^{\prime }=\overline{%
\left\langle F_i\tilde{z}_i\right\rangle _{i=1}^\infty }$ and $\tilde{Z}%
^{\prime }=\overline{\left\langle E_i\tilde{y}_i,F_i\tilde{z}_i\right\rangle
_{i=1}^\infty }$. Then $\tilde{Y}_1^{\prime }$, $\tilde{Y}_2^{\prime }$ are
block subspaces of $\tilde{Z}^{\prime }$ and there exists a constant $%
\varepsilon ^{\prime }>0$ such that for every $\tilde{y}\in \tilde{Y}%
_1^{\prime }$ and $\tilde{z}\in Y_2^{\prime }$, $||\tilde{y}-\tilde{z}%
||>\varepsilon ^{\prime }||\tilde{y}+\tilde{z}||$ therefore $\tilde{Z}%
^{\prime }$ is not an $H.I$-space, a contradiction. 
\endproof%

\begin{definition}
\label{24}A bounded subset $W$ of a Banach space $X$ is called a $\mathbf{C}$%
{\bf -norming} set for a subspace $Z$ of $X^{*}$,where $C$ is a positive
number,if\newline
$||z||\leq C\underset{w\in W}{\sup }|z(w)|$, for any $z\in Z$.\newline
The set $W$ is a {\bf norming} set for the subspace $Z$ if it is $C$%
-norming for $Z$ for some $C>0$.
\end{definition}

\begin{proposition}
\label{Prop2.4}Let $Y$ be any diagonal space produced by a subset $W$ of a
Banach space $X$, which is a $C$-norming set of a subspace $Z$ of $X^{*}$.
Then $Z$ is isomorphic to a closed subspace of $Y^{*}$.
\end{proposition}

\proof%
Denote by $J^{*}$ the adjoint of the operator $J$. From our assumption for
any $z\in Z$ there exists $w_z\in W$ such that $||z||<2C|z(w_z)|$.\newline
Since $W\subseteq B_Y,||J^{*}z||\geq |J^{*}(J^{-1}(w_z))|=|z(w_z)|\geq \frac
1{2C}||z||$. Therefore $J^{*}|_Z$ is an isomorphism onto a subspace of $%
Y^{*} $%
\endproof%

\begin{proposition}
\label{Prop2.5}Let $X,Y,Z,W$ as before and $A$ a Banach space such that $%
A^{*}$ is $w^{*}$-isomorphic to $Z$. Then $A$ is a quotient of $Y$.
\end{proposition}

\begin{proof}
Let $S:A\rightarrow Z$ be a $w^{*}$-continuous isomorphism. Then $J^{*}\circ
S$ is also a $w^{*}$-continuous isomorphism hence $(J^{*}\circ S)^{*}$ maps
the space $Y$ onto the space $A$.
\end{proof}

\section{{\bf THIN NORMING SETS I}}

\begin{center}
{\bf (special reflexive cases)}
\end{center}

This is the first section devoted to the $\mathbf{a-}$thin (thin) norming
sets.\newline
An $\mathbf{a-}$thin subset $W$ of a Banach space $X$ is norming for a
subspace $Y$ of $X^{*}$ if it defines an equivalent norm on $Y$ (Definition $%
\ref{24}$ ).\newline
In this section, for a reflexive Banach space $A$ with an unconditional
basis we construct a space $X_A$ which is of the form $\left(
\sum\limits_{n=1}^\infty \bigoplus \ell ^1\left( k_n\right) \right) _A$ and
a symmetric convex closed subset $W$ of the unit ball of $X_A$ which, in
certain cases, is an $\mathbf{a-}$thin subset of $X_A$ and norms a subspace
of $X_A^{*}$ isometric to $A^{*}.$\newline
The set $W$ is of the form $\overline{co}\left( K\cup -K\right) $ where $K$
forms a tree in the branches of which we have ``spread'' in a regular way a
dense subset of the positive part of the unit ball of $A$.\newline
The fact that such a $W$ norms $A^{*}$ is quite easy and it is true for all
reflexive Banach spaces $A$ with an unconditional basis.\newline
The difficult part to be shown is that $W$ is an $\mathbf{a-}$thin set. For
this we need certain measure-theoretical combinatorial results.\newline
The first, Proposition \ref{AC}, ``the finite version'', will be used in the
next section and the second, Proposition \ref{AE}, ``the infinite version'',
will be used in this section.\newline
The property (P) given in Definition \ref{DP} is the main tool we use to
show that for certain classes of reflexive spaces the corresponding set $W$
is an $\mathbf{a-}$thin set. The formulation of this, in a slightly
different form, was suggested to us by B.Maurey.\newline
We show that if $A$ satisfies the property (P) then the set $W\,$ in the
space $X_A$ is an $\mathbf{a-}$thin set for an appropriate sequence $\mathbf{%
a}$. Thus we are able to prove that for spaces like $L^p\left( \lambda
\right) ,$ $\ell ^p\left( \mathbb{N}\right) ,$ $1<p<\infty $, and others, the
corresponding set $W$ is an $\mathbf{a-}$thin set.\newline
Finally we prove a corresponding result for spaces with many complemented
subspaces, suggested to us by N. Tomczak-Jaegermann.

(a) {\bf Definition of the space }$\mathbf{X}_A$\newline
For every reflexive Banach spaces with 1-unconditional basis we shall
construct a Banach space $X_A$ and a subset $W$ of $X_A$ which is norming
for a subspace $Z$ of $X_A^{*}$ isometric to $A^{*}$ and it is $\mathbf{a}$%
-thin, for an appropriate sequence $\mathbf{a}$ if $A$ satisfies some
general further condition.

Consider a tree $\frak{D}$ of height $\omega $ with a least element $\varrho 
$ such that every $\delta $ in $\frak{D}$ with height equal to $n$ has $%
3^{4(n+1)}+1$ immediate successors. We can define $\frak{D}$ as the set of
all finite sequences $(k_1,...k_n)$, $n\in \mathbb{N}$ of natural numbers such
that $0\leq k_i\leq 3^{4i}$ and with least element the empty sequence.
Consider the partial order $\prec $ in $\frak{D}$: $(k_1,...k_n)\prec
(k_1^{\prime },...,k_\ell ^{\prime })$ iff $(k_1,...,k_n)$ is an initial
segment of $(k_1^{\prime },...,k_\ell ^{\prime })$, i.e. $n\leq \ell $ and $%
k_i=k_i^{\prime }$ for any $i=1,2,...,n$. We write $s\sqsubset \gamma $ if $%
s $ is an initial segment of $\gamma $.\newline
A {\bf segment }is a subset of $\frak{D}$ of the form $d=[\delta
_1,\delta _2]=\left\{ \delta \in \frak{D}:\delta _1\prec \delta \prec \delta
_2\right\} .$\newline
A segment $d$ of the form $d=[\rho ,\delta ]$ is called an {\bf initial
segment }of $\frak{D}$ .\newline
For every segment $d=[\delta _1,\delta _2]$we define $ext(d)=[\rho ,\delta
_2].$\newline
A {\bf branch } is maximal linearly ordered subset of $\frak{D}$ .{\bf %
\ }\newline
A branch is identified to an infinite sequence $(k_i)_{i\in \mathbb{N}}$ of
natural numbers, where $0\leq k_i\leq 3^{4i}$ \newline
The {\bf height} of $\delta =(k_1,...,k_n)$ is denote $|\delta |$ and it
is equal to $n$.We set $\left| \rho \right| =0.$ \newline
We denote by $c_{00}(\frak{D})$ the linear space of all functions $f:\frak{D}%
\rightarrow \mathbb{R}$ such that $\mathit{supp}(f)=\{\delta \in \frak{D}%
:f(\delta )\neq 0\}$ is a finite set.\newline
Also,denote by $e_{\delta _{}}$ the characteristic function of $\{\delta \}$%
, $\delta \in \frak{D}$. The vectors $(e_{\delta _{}})_{\delta \in \frak{D}}$%
, form a Hamel basis for $c_{00}(\frak{D})$.\newline
If $A,B$ are finite subsets of $\frak{D}$ we write\newline
$A<B$ iff $\max \{|\alpha |:\alpha \in A\}<\min \{|\beta |:\beta \in B\}$,
and\newline
$A\prec B$ iff $\alpha \prec \beta $ for any $\alpha \in A$ and any $\beta
\in B$.\newline
If $f\in c_{00}(\frak{D})$, we set $range\left( f\right) =\{\delta \in \frak{%
D}:|\alpha |\leq |\delta |\leq |\beta |$ for some $\alpha ,\beta \in $%
\textit{supp}$(f)\}.$We call the set $range\left( f\right) $ the {\bf %
range} of $f$.\newline
If $f,g\in c_{00}(\frak{D})$ we write $f<g$ iff \textit{supp}$(f)<$\textit{%
supp}$(g)$ and $f\prec g$ iff \textit{supp}$(f)\prec $\textit{supp}$(g)$.%
\newline
For a given space $A$ with an unconditional basis $(e_n)_{n\in \mathbb{N}}$ we
define the following norm of $c_{00}(\frak{D})$:\newline
For $f=\underset{\delta \in \frak{D}}{\sum }\lambda _{\delta _{}}e_{\delta
_{}}\in c_{00}(\frak{D})$ we set\newline
\begin{equation*}
||f||=\left\| \underset{n=0}{\overset{\infty }{\sum }}\left( \underset{%
|\delta |=n}{\sum }|\lambda _{\delta _{}}|\right) e_n\right\| _A
\end{equation*}
\newline
We denote by $X_A$ the completion of $c_{00}(\frak{D})$ with the above
defined norm.\newline
It is clear that $X_A$ is isometric to $\left( \underset{n=0}{\overset{%
\infty }{\sum }}\bigoplus \ell ^1(k_n)\right) _A$ for some sequence $%
(k_n)_{n=1}^\infty $ and if $A$ has a shrinking basis $X_A^{*}=\left( 
\underset{n=0}{\overset{\infty }{\sum }}\bigoplus \ell ^\infty (k_n)\right)
_{A^{*}}$. If $A$ is reflexive Banach space then $X_A$ is also reflexive.%
\newline
Set $y_n^{*}=\underset{|\delta |=n}{\sum }e_{\delta _{}}^{*}$ and $Y^{*}=%
\overline{\left\langle (y_n^{*})_{n\in \mathbb{N}}\right\rangle }$. Then $%
||y_n^{*}||=1$ and it is easy to see that $Y^{*}$ is isometric to $A^{*}$.
In the sequel we consider reflexive spaces $A$ which the basis $(e_n)_{n\in 
\mathbb{N}}$ of $A$ is $1$-unconditional.\newline
$_{}$\newline
$_{}$\newline
b) {\bf Definition of the set }$\mathbf{W}$\newline
For every $\delta =(k_1,...,k_n)$we set $a_{\delta _{}}=\frac{k_n}{3^{4n}}$%
\newline
For any infinite branch $\gamma $ of $\frak{D}$ we denote by $x_{\gamma _{}}$%
the typical series 
\begin{equation*}
x_{\gamma _{}}=\underset{\delta \in \gamma }{\sum }a_{\delta _{}}e_{\delta
_{}}
\end{equation*}
where $\delta _n=(k_1,...,k_n)\in \frak{D}$ the initial segment of $\gamma $%
. \newline
Also we set 
\begin{equation*}
|||x_{\gamma _{}}|||=\underset{n=0}{\overset{\infty }{\sum }}\frac{k_n}{%
3^{4n}}e_{\delta _n}
\end{equation*}
\newline
It follows from the definitions that if $|||x_{\gamma _{}}|||<\infty $ then $%
x_{\gamma _{}}\in X_A$ and $||x_{\gamma _{}}||=|||x_{\gamma _{}}|||$.\newline
We set\newline
\begin{equation*}
K=\{x_{\gamma _{}}:||x_{\gamma _{}}||\leq 1\}\text{ and }W=\overline{co}%
(K\cup -K)
\end{equation*}
\newline
The set $W$ is a closed bounded symmetric convex subset of $B_{X_A}$.

\begin{lemma}
The set $K$ (and so $W$) is a $\frac 14$-norming set for the space $%
Y^{*}\subseteq X_A^{*}$.
\end{lemma}

\proof%
Let $y^{*}=\underset{k=1}{\overset{\infty }{\sum }}\lambda _ky_k^{*}\in
Y^{*} $ such that $||y^{*}||=\left\| \underset{k=1}{\overset{\infty }{\sum }}%
\lambda _ke_k^{*}\right\| _{A^{*}}=1$.\newline
Define $z^{*}=\underset{k=1}{\overset{\infty }{\sum }}\lambda _ke_k^{*}\in
A^{*}$. Then $||z^{*}||_{A^{*}}=||y^{*}||=1$\newline
Let $z=\underset{k=1}{\overset{\infty }{\sum }}\mu _ke_k\in A$ be such that $%
|z^{*}(z)|=\left| \underset{k=1}{\overset{\infty }{\sum }}\lambda _k\mu
_k\right| =1$. Since $(e_n)_{n\in \mathbb{N}}$ is $1$-unconditional and
bimonotone\newline
\begin{equation*}
1=\left| \underset{k=1}{\overset{\infty }{\sum }}\lambda _k\mu _k\right|
\leq \left| \underset{\mu _k>0}{\sum }\lambda _k\mu _k\right| +\left| 
\underset{\mu _k<0}{\sum }\lambda _k\mu _k\right|
\end{equation*}
\newline
So, we can suppose that $\left| \underset{\mu _k>0}{\sum }\lambda _k\mu
_k\right| \geq \frac 12$.\newline
For any $k\in \mathbb{N}$ we define $n_k\in \mathbb{N}$ as follows.\newline
(i) \thinspace If $\mu _k\leq 0$ then $n_k=0$\newline
(ii) If $\mu _k>0$ then $n_k$ is the unique natural number with the property%
\newline
$\frac{n_k}{3^{4k}}\leq \mu _k<\frac{n_k+1}{3^{4k}}$.\newline
Since $(e_n)_{n\in \mathbb{N}}$ is bimonotone, $\mu _k\leq 1$ for any $k\in 
\mathbb{N}$ and therefore $n_k\leq 3^{4k}$ for any $k\in \mathbb{N}$.\newline
Let $\gamma =(n_1,n_2,...)$, an infinite branch of $\frak{D}$. Then

\begin{center}
$||x_{\gamma _{}}||=\left\| \underset{k=1}{\overset{\infty }{\sum }}\frac{n_k%
}{3^{4k}}e_k\right\| _A\leq \left\| \underset{k=1}{\overset{\infty }{\sum }}%
\mu _ke_k\right\| _A\leq 1$, \thinspace \thinspace i.e. \thinspace $%
x_{\gamma _{}}\in K$.
\end{center}

And,\newline
$|y^{*}(x_{\gamma _{}})|=\left| \underset{k=1}{\overset{\infty }{\sum }}%
\lambda _k\frac{n_k}{3^{4k}}\right| =\left| \underset{\mu _k>0}{\sum }%
\lambda _k\frac{n_k}{3^{4k}}\right| \geq \left| \underset{\mu _k>0}{\sum }%
\lambda _k\mu _k\right| -\left| \underset{\mu _k>0}{\sum }\lambda _k\left(
\mu _k-\frac{n_k}{3^{4k}}\right) \right| \geq $

$\geq \frac 12-\underset{k=1}{\overset{\infty }{\sum }}\frac 1{3^{4k}}>\frac
14=\frac 14||y^{*}||$\newline
and the proof is completed. 
\endproof%

In the following lemma $\phi $ is a function on $X_A$ which is either

i) $\phi (x)=y^{*}(x)$ for some non-negative functional $y^{*}$ or

ii) $\phi (x)=||x||$.

In fact we may consider any real valued bounded sublinear function $\phi $
on $X_A$ such that for any segments $d_1$, $d_2$ of $\frak{D}$ with $%
d_1\subset d_2$, $0\leq \phi (x_{d_1})\leq \phi (x_{d_2})$.\newline
It is more convenient to consider the set $W_0=co\left( K\cup -K\right) $
instead of $W=\overline{W_0}.$

\begin{lemma}
$\label{AA}$Let $E$ be a subset of $\frak{D}$ of the form $E=\{\delta \in 
\frak{D}:m\leq |\delta |\leq M\}$, $\phi $ a function on $X_A$ as before and 
$\varepsilon >0$.\newline
There exists a decomposition of $E$ into two disjoint subsets $E^{\prime }$, 
$E^{\prime \prime }$ such that\newline
(1) $|\phi (E^{\prime }w)|<\varepsilon $ for every $w\in W_0$\newline
(2) If $d=[\delta _1,\delta _2]$ is a segment of the tree such that $|\delta
_1|\leq m$, $|\delta _2|\geq M$ and $d\cap E^{\prime \prime }\neq \emptyset $
then $|\phi (Ex_d)|\geq \varepsilon $.
\end{lemma}

\proof%
For every $\delta \in E$ let $ext_E(\delta )=ext(\delta )\cap E=[\delta
^E,\delta ]$, where $\delta ^E$ is the unique element of $E$ such that $%
\delta ^E<\delta $ and $|\delta ^E|=m$. We set\newline
\begin{equation*}
E^{\prime \prime }=\{\delta \in E:||x_{ext(\delta )}||\leq 1\text{ and }\phi
(x_{ext_E(\delta )})\geq \varepsilon \newline
\}
\end{equation*}
\begin{equation*}
E^{\prime }=E\backslash E^{\prime \prime }\text{.}
\end{equation*}
\newline
Let $w\in W_0$, then $Ew$ can be represented as

\begin{equation*}
Ew=\underset{d\in L}{\sum }\lambda _dx_d
\end{equation*}
\newline
where $L$ is a subset of segments of $E$ such that $||x_{ext(d)}||\leq 1$
for every $d\in L$ and $\underset{d\in L}{\sum }|\lambda _d|\leq 1$. For
every $d\in L$ the set $d^{\prime }=d\cap E^{\prime }$ is either empty or a
segment of $E$ and $\phi (x_{d^{\prime }})<\varepsilon $. So $|\phi
(E^{\prime }w)|\leq \underset{d\in L}{\sum }|\lambda _d|$\thinspace $\,|\phi
(x_{d^{\prime }})|<\varepsilon $, which proves property (1).\newline
Property (2) of $E^{\prime \prime }$ is obvious%
\endproof%

The following two Propositions \ref{AC} , \ref{AE} ,are of combinatorial
nature and we will use them to show that the set $W$ is an $\mathbf{a}$-thin
set. Both are related to the existence of incomparably supported elements of
the set $W$ .The present form which is a variation of our original approach
suggested by B.Maurey.\newline
The following lemma is well known and the sake of completeness we give a
proof of it.

\begin{lemma}
$\label{AB}$Let $(\Omega ,\mu )$ be a measure space with $\mu (\Omega )\leq 1
$ and let $B_1,...,B_N$ be measurable subsets of $\Omega $ which satisfy $%
\mu (B_i)\geq \varepsilon $ for $i=1,...,N$. If $k<N\varepsilon $ there
exist $B_{i_1}...B_{i_k}$ $1\leq i_1\leq ...\leq i_k\leq n$ such that $%
\underset{j=1}{\overset{k}{\bigcap }}B_{i_j}\neq \emptyset $.
\end{lemma}

\proof%
If not then the function $\underset{i=1}{\overset{N}{\sum }}x_{B_i}$ is
bounded by $k$ on $\Omega $ and therefore 
\begin{equation*}
N\varepsilon \leq \underset{\Omega }{\int }\left( \overset{N}{\underset{i=1}{%
\sum }}x_{B_i}\right) d\mu \leq k
\end{equation*}

a contradiction.%
\endproof%

Next we prove the first incomparability result released to finite families
of elements of the set $W$.

\begin{proposition}
\label{AC}Let $(E_i)_{i=1}^n$ be a finite sequence of subsets of $\frak{D}$
of the form\newline
$E_i=\{\delta \in \frak{D}:m_i\leq |\delta |\leq M_i\},\,\,m_i,M_i\in \mathbb{N}%
,\,\,m_i<M_{i+1}$ for $i=1,...,n-1$.\newline
Then for every $\varepsilon >0$, $x^{*}\in B_{X_A}^{+}$ and $(w_i)_{i=1}^n$
a finite sequence of vectors of $W_0$ there exists a partition $D_1,...,D_N$
of $\{1,2,...,n\}$ and $F_i\subset E_i$ for $i=1,...,n$ such that\newline
(1) If $r,s$ belong to the same set $D_i$ and $r<s$ then $ext(\delta )\cap
F_r=\emptyset $ for every $\delta \in F_s$ (i.e. the sets $(F_r)_{r\in D_i}$
are pairwise incomparable for every $i=1,...,N$)\newline
(2) $x^{*}\left( (E_i\backslash F_i)w_i\right) <\varepsilon $ for every $%
i=1,...,n$\newline
(3) The number $N$ of the members of the partition is less than $\frac
5{\varepsilon ^2}$.
\end{proposition}

\proof%
Applying lemma \ref{AA} we find for every $i=1,...,n$ a decomposition of $%
E_i $ into two disjoint subsets $E_i^{\prime }$, $E_i^{\prime \prime }$ such
that\newline
(a) $|x^{*}(E_i^{\prime }w_i)|<\frac \varepsilon 2$\newline
(b) If $d=[\delta _1,\delta _2]$ is a segment such that $|\delta _1|\leq m_i$%
, $|\delta _2|\geq M_i$ and $d\cap E_i^{\prime \prime }\neq \emptyset $ then 
$x^{*}(x_d)>\frac \varepsilon 2$.\newline
For every $i=1,2,...,n$ the vector $E_iw_i$ has a representation

\begin{equation*}
E_iw_i=\underset{d\in L_i}{\sum }\lambda _dx_d
\end{equation*}
\newline
where $L_i$ is a set of segments of $E_i$ with $||x_{ext(d)}||\leq 1$ for
every $d\in L_i$ and $\underset{d\in L_i}{\sum }\lambda _d\leq 1$.\newline
This representation defines a positive measure $\mu _i$ on the set 
\begin{equation*}
seg(E_i)=\{d:d\text{ is a segment of }E_i\text{ and }||x_{ext(d)}||\leq 1\}
\end{equation*}
,with $||\mu _i||\leq 1,$ as follows: 
\begin{equation*}
\mu _i(A)=\underset{d\in A\cap L_i}{\sum }|\lambda _d|\,\,\,\text{for every
\thinspace }A\subset seg(E_i)\text{.}
\end{equation*}
\newline
For every $S\subset \{1,...,n\}$ and every $i\in \{1,...,n\}$ we set\newline
\begin{equation*}
A_S^i=\left\{ d\in seg(E_i):ext(d)\cap \left( \underset{s\in S\backslash
\{i\}}{\bigcup }E_s^{\prime \prime }\right) \neq \emptyset \right\} \text{.}
\end{equation*}
\newline
Given any non-empty subset $J$ of $\{1,...,n\}$ we inductively define a
non-empty subset $D)J)$ of $J$ by the following manner:\newline
Suppose that $J=\{i_1,...,i_m\}$ then\newline
(1) $\min J=j_1\in D(J)$.\newline
(2) $i_2\in D(J)$ if and only if $\mu _{i_2}\left( A_{\{i_1\}}^{i_2}\right)
<\frac \varepsilon 2$.\newline
(3) If $k\leq m$ and we have selected the elements $\{j_1,...,j_r\}$ of $%
D(J) $ up to $k-1$ then the $i_k$ will belong to $D(J)$ if and only if%
\newline
\begin{equation*}
\mu _{i_k}\left( A_{\{j_1,...,j_r\}}^{i_k}\right) <\frac \varepsilon 2\text{.%
}
\end{equation*}
\newline
This complete the inductive definition of $D(J)$.\newline
We set also 
\begin{equation*}
F_j=\left\{ \delta \in E_j^{\prime \prime }:ext(\delta )\cap \bigcup \left\{
E_s^{\prime \prime }:s<j,s\in D(J)\right\} =\emptyset \right\}
\end{equation*}
\newline
for every $j\in D(J)$.\newline
The sets $(F_j)_{j\in D(J)}$ are pairwise incomparable by definition.\newline
We show that $|x^{*}(E_j\backslash F_j)w_j|<\varepsilon $ for every $j\in
D(J)$.\newline
For $d\in L_j$ we set $d^{\prime }=E_j^{\prime }\cap d$, $d_a^{\prime \prime
}=F_j\cap d$ and $d_b^{\prime \prime }=(E_j^{\prime \prime }\backslash
F_j)\cap d$.\newline
It is clear that for every $d\in L_j$\thinspace \thinspace $d_b^{\prime
\prime }\in A_{D(J)}^j$ and since $\mu _j\left( A_{D(J)}^j\right) <\frac
\varepsilon 2$ we get that\newline
\begin{equation*}
|x^{*}((E_j\backslash F_j)w_j)|\leq |x^{*}(E_j^{\prime }w_j)|+\left|
x^{*}\left( \underset{d\in L_j}{\sum }\lambda _dx_{d_b^{\prime \prime
}}\right) \right| \leq \frac \varepsilon 2+\underset{d\in A_{D(J)}^j}{\sum }%
|\lambda _d|<\varepsilon \text{.}
\end{equation*}
\newline
We define the partition $D_1,...,D_N$ of $\{1,2,...,n\}$ as follows: $%
J_1=\{1,...,n\}$, $D_1=D(J_1)$, $J_2=J_1\backslash D_1$,$D_2=D(J_2)$ and we
continue until the set $J_{N+1}$ will be empty.\newline
It remains to show that $N<\frac 5{\varepsilon ^2}$.\newline
Indeed, since $J_N\neq \emptyset $ let $r\in J_N$. Then $\mu _r\left(
A_{D_j}^r\right) \geq \frac \varepsilon 2$ for $j=1,...,N-1$.\newline
If $N>\frac 5{\varepsilon ^2}$ then $\left( N-1\right) \frac \varepsilon
2>\frac 2\varepsilon $. By lemma \ref{AB} we can select $D_{j_1},...,D_{j_k}$
with $j_1<...<j_k$ and $\underset{s=1}{\overset{k}{\bigcap }}%
A_{D_{j_s}}^r\neq \emptyset $ for some $k$ such that $\frac \varepsilon 2k>1$%
.\newline
Let $d\in \underset{s=1}{\overset{k}{\bigcap }}A_{D_{j_s}}^r$ then for every 
$s=1,...,k$ there exists $r_s\in D_{j_s}$ such that $d\cap E_{r_s}^{\prime
\prime }\neq \emptyset $ so $x^{*}(E_{r_s}d)\geq \frac \varepsilon 2$ and $%
x^{*}(x_{ext(d)})\geq \frac \varepsilon 2k>1$ which is a contradiction.%
\endproof%

Before we prove the infinite analogous of the above proposition we give the
statement and the proof of a known auxiliary result.

\begin{lemma}
\label{AD}Let $(\Omega ,\mu )$ be a measure space with $\mu (\Omega )<\infty 
$, $\varepsilon >0$ and $(B_n)_{n=1}^\infty $ be a sequence of measurable
subsets of $\Omega $ such that $\mu (B_n)\geq \varepsilon $ for every $n\in 
\mathbb{N}$. There exists an infinite subset $M$ of $\mathbb{N}$ such that $%
\underset{n\in M}{\bigcap }B_n\neq \emptyset $.
\end{lemma}

\proof%
Since $\mu (\Omega )<\infty $, it follows that 
\begin{equation*}
\mu (\lim \sup B_n)=\mu \left( \underset{n=1}{\overset{\infty }{\bigcap }}%
\underset{m=n}{\overset{\infty }{\bigcup }}B_n\right) \geq \lim \inf \mu
(B_n)\geq \varepsilon \text{.}
\end{equation*}
\newline
Hence $\lim \sup B_n\neq \emptyset $ and for $\omega \in \lim \sup B_n$
there exists $M$ infinite subset of $\mathbb{N}$ such that $\omega \in $ $%
\underset{n\in M}{\bigcap }B_n$.%
\endproof%

\begin{proposition}
\label{AE}Let $(E_n)_{n=1}^\infty $ be a sequence of subsets of $\frak{D}$
of the form\newline
$E_n=\{\delta \in \frak{D}:m_n\leq |\delta |\leq M_n\}$, $m_n,M_n\in \mathbb{N}$
and $m_{n+1}>M_n$ $\forall n\in \mathbb{N}$.\newline
For every sequence $(w_n)_{n=1}^\infty $ of vectors of $W_0$ and $%
\varepsilon >0$ there exists an infinite subset $I$ of $\mathbb{N}$ and $%
F_n\subset E_n$ for every $n\in I$ such that\newline
1. The sets $(F_n)_{n\in I}$ are pairwise incomparable\newline
2. For every $n\in I$, $||(E_n\backslash F_n)w_n||<\varepsilon $.
\end{proposition}

\proof%
We apply Lemma \ref{AA} to find a decomposition of each $E_n$ into two
disjoint sets $E_n^{\prime },E_n^{\prime \prime }$ such that\newline
(1) $||E_n^{\prime }w_n||<\frac \varepsilon 2$\newline
(2) If $d=[\delta _1,\delta _2]$ is a segment with $|\delta _1|\leq m_n$, $%
|\delta _2|\geq M_n$ and $d\cap E_n^{\prime \prime }\neq \emptyset $ then $%
||E_nx_d||\geq \frac \varepsilon 2$.\newline
For every $n\in \mathbb{N}$ the vector $E_nw_n$ has a representation as 
\begin{equation*}
E_nw_n=\underset{d\in L_n}{\sum }\lambda _dx_d
\end{equation*}
\newline
where $L_n\subseteq seg(E_n)=\{d:d\subset E_n$ and $||x_{ext(d)}||\leq 1\}$
and \newline
$\underset{d\in L_n}{\sum }|\lambda _d|\leq 1$.\newline
This representation defines a positive measure $\mu _n$ on $seg(E_n)$ as in
the previous proposition: 
\begin{equation*}
\mu _n(A)=\underset{d\in A\cap L_n}{\sum }|\lambda _d|\text{, for every }%
A\subset seg(E_n)\text{.}
\end{equation*}
\newline
Consider any probability measure $\nu $ on the compact metrisable space $%
\Gamma $ of all infinite branches of $\frak{D}$ such that $\nu (W_d)\neq
\emptyset $ for every $d$ segment of the tree, where $W_d$ is the basic
clopen subset of $\Gamma $ consisting of all branches that contain $d$.%
\newline
For instance we can consider the natural measure $\nu $ defined by 
\begin{equation*}
\nu (W_{ext(\delta )})=\frac 1{k_{|\delta |}}
\end{equation*}
\newline
for every $\delta \in \frak{D}$, where $k_n=\#\{\delta \in \frak{D}:|\delta
|=n\}$.\newline
We define a measure $\mu $ on $\Gamma $ as follows:\newline
For every clopen subset $B$ of $\Gamma $ 
\begin{equation*}
\mu (B)=\underset{n\rightarrow \frak{U}}{\lim }\underset{d\in L_n}{\sum }%
|\lambda _d|\frac{\nu (W_d\cap B)}{\nu (W_d)}
\end{equation*}
\newline
where the limit is taken with respect some non-trivial ultrafilter $\frak{U}$
on $\mathbb{N}$.\newline
Using a diagonal argument we may actually assume that this limit is an
ordinary limit ,by passing to some subsequence and we shall consider that
this subsequence is the whole sequence of natural numbers.\newline
For every $n,s\in \mathbb{N}$ with $n<s$ we define 
\begin{equation*}
B_n=\{\gamma \in \Gamma :\gamma \cap E_n^{\prime \prime }\neq \emptyset \}
\end{equation*}
\newline
\begin{equation*}
A_n^s=\{d\in seg(E_s):ext(d)\cap E_n^{\prime \prime }\neq \emptyset \}
\end{equation*}
\newline
and observe that if $d\in A_n^s$ then $W_d\subseteq B_n$ and if $d\in
L_s\backslash A_n^s$ then $W_d\cap B_n=\emptyset $.\newline
Therefore 
\begin{equation*}
\underset{d\in L_s}{\sum }|\lambda _d|\frac{\nu (B_n\cap W_d)}{\nu (W_d)}=%
\underset{d\in L_s\cap A_n^s}{\sum }|\lambda _d|=\mu _s(A_n^s)
\end{equation*}
\newline
and so for every $n\in \mathbb{N}$\newline
\begin{equation}
\mu (B_n)=\underset{s\rightarrow \infty }{\lim }\mu _s(A_n^s)  \tag{1}
\label{!}
\end{equation}

{\bf Claim:} For every infinite subset $I$ of $\mathbb{N}$ and $\theta >0$
the set $I_{\theta _{}}=\{n\in I:\mu (B_n)<\theta \}$ is infinite.\newline
Indeed, otherwise there exists an infinite subset $J$ of $I$ such that $\mu
(B_n)\geq \theta $ for every $n\in J$ and by Lemma \ref{AD} there exists an
infinite subset $J^{\prime }$ of $J$ and $\gamma \in \underset{n\in
J^{\prime }}{\bigcap }B_n$. For every $n\in J^{\prime }$, $\gamma \cap
E_n^{\prime \prime }\neq \emptyset $ and so $||E_n^{\prime \prime }x_{\gamma
_{}}||\geq \frac \varepsilon 2$ for infinite values of $n$. Select $\delta
_n\in E_n^{\prime \prime }\cap \gamma $ then $||x_{ext(\delta _n)}||\leq 1$
so $||x_{\gamma _{}}||\leq 1$ which contradicts the fact the basis of $X_A$
is boundedly complete.\newline
We define 
\begin{eqnarray*}
\tilde{I}_1 &=&\left\{ n\in \mathbb{N}:\mu (B_n)<\frac \varepsilon
{2^2}\right\} ,\;\;\;\;\;_{}\;\;\;\;_{}\;\;\;\;\;n_1=\min \tilde{I}_1\text{,}
\\
I_1 &=&\left\{ s\in \tilde{I}_1:\mu _s(A_{n_1}^s)<\frac \varepsilon
{2^2}\right\} \,\,\,\,\,\,\,\,\,\,\,\,\,\,_{}\text{and} \\
F_{n_1} &=&E_{n_1}^{\prime \prime }
\end{eqnarray*}
\newline
and inductively for $k\geq 1$\newline
\begin{eqnarray*}
\tilde{I}_{k+1} &=&\left\{ n\in I_k:\mu (B_n)<\frac \varepsilon
{2^{k+2}}\right\} ,\;\;\;\;\;_{}\;\;\;\;_{}\;\;\;\;\;n_{k+1}=\min \tilde{I}%
_{k+1}\text{,} \\
I_{k+1} &=&\left\{ s\in \tilde{I}_{k+1}:\mu _{n_{k+1}}(A_{n_{k+1}}^s)<\frac
\varepsilon {2^{k+2}}\right\} \,\,\,\,\,\,\,\,\,\,\,\,\,\,_{}\text{and} \\
F_{n_{k+1}} &=&\left\{ \delta \in E_{n_{k+1}}^{\prime \prime }:ext(\delta
)\cap \left( \overset{}{\overset{k}{\underset{i=1}{\bigcup }}E_{n_i}^{\prime
\prime }}\right) =\emptyset \right\} \text{.}
\end{eqnarray*}
\newline
By the previous claim and relation (\ref{!}) the sets $\tilde{I}_{k+1}$, $%
I_{k+1}$ are infinite subsets of $I_k$ and since $n_{k+1}\in I_k$%
\begin{equation}
\mu _{n_{k+1}}(A_{n_k}^{n_{k+1}})<\frac \varepsilon {2^{k+1}}\text{.} 
\tag{2}  \label{a}
\end{equation}
\newline
We set $I=\{n_1,n_2,...\}$ and we observe that the sets $(F_{n_j})_{j=1}^%
\infty $ are pairwise incomparable by their own definition. \newline
We show that $||(E_{n_k}\backslash F_{n_k})w_{n_k}||<\varepsilon $.\newline
If $k=1$ it follows for the fact that $F_{n_1}=E_{n_1}^{\prime \prime }$.%
\newline
Let $k\in \mathbb{N}$ and set $r=n_{k+1}$. Consider the set: 
\begin{equation*}
A_r=A_{n_1}^r\cup ...\cup A_{n_k}^r\text{.}
\end{equation*}
\newline
by (\ref{a}) we get that 
\begin{equation}
\mu _r(A_r)\leq \frac \varepsilon {2^2}+...+\frac \varepsilon
{2^{k+1}}<\frac \varepsilon 2\text{.}  \tag{3}
\end{equation}
\newline
Let $d\notin A_r$ then $ext(\delta )\cap E_{n_i}^{\prime \prime }=\emptyset $
for every $\delta \in d$ and $i=1,...,k$, so $d\subset F_r$. Conversely if $%
d\in A_r$ then $d\cap F_r=\emptyset $. So 
\begin{equation*}
||(E_r^{\prime \prime }\backslash F_r)w_r||\leq \underset{d\in L_s\cap A_r}{%
\sum }|\lambda _d|\,\,\,||(E_r^{\prime \prime }\backslash F_r)x_d||\leq \mu
_r(A_r)<\frac \varepsilon 2\text{.}
\end{equation*}
\newline
Finally we have 
\begin{equation*}
||(E_r\backslash F_r)w_r||\leq ||(E_r^{\prime }w_r||+||(E_r^{\prime \prime
}\backslash F_r)w_r||<\frac \varepsilon 2+\frac \varepsilon 2=\varepsilon 
\text{.}%
\end{equation*}
\endproof%

\begin{definition}
\label{DP}Let $(A,||$\thinspace $||_A)$ be a reflexive Banach space with an
unconditional basis $(e_n)_{n=1}^\infty $. We say that $A$ {\bf satisfies
the property} (P) if there exists a sequence $(C_k)_{k=1}^\infty $of
positive numbers, such that:\newline
For every block sequence $(x_n)_{n=1}^\infty $ of $(e_n)_{n=1}^\infty $ and
every $k\in \mathbb{N}$ there exists a normalized block sequence $%
(y_n)_{n=1}^\infty $ of $(x_n)_{n=1}^\infty $ satisfying the following
condition.\newline
For every infinite subset $M$ of $\mathbb{N}$ there exists $E\subset M$ and $%
(\lambda _n)_{n\in E}\subseteq \mathbb{R}^{+}$, $(y_n^{*})_{n\in E}\subseteq
B_{A^{*}}$ such that\newline
\begin{equation}
\underset{n\in E}{\sum }\lambda _n=1,\left\| \underset{n\in E}{\sum }\lambda
_ny_n\right\| _A<\frac 1k  \tag{i}  \label{i}
\end{equation}
\newline
\begin{equation}
\mathit{supp}y_n^{*}\subseteq \mathit{supp}y_n\;\text{and\ }%
y_n^{*}(y_n)>\frac 12\;\text{for every }n\in E  \tag{ii}  \label{ii}
\end{equation}
\newline
\begin{equation}
\text{For every\ }x\in B_A,\left\| \underset{n\in E}{\sum }%
y_n^{*}(x)y_n\right\| _A\leqslant C_k  \tag{iii}  \label{iii}
\end{equation}
\end{definition}

\begin{remark}
\label{Rem3.1}Condition (\ref{iii}) in the above definition is equivalent to
say that the norm of the operator $P:X\rightarrow <\overline{(y_n)_{n\in E}>}
$ defined by the relation $P(x)=\underset{n\in E}{\sum }y_n^{*}(x)y_n$ is
dominated by $C_k$.\newline
There other conditions which imply condition (\ref{iii}) and in some cases
is more convenient to be checked. We give two of them that we will use later:%
\newline
(P1):\ \ \ \ $\left\| \underset{n\in E}{\sum y_n^{*}}\right\| \leq C_k$%
\newline
(P2): \ \ \ \ $\left| E\right| \leq C_k$\newline
In the sequel by property (P$_i),i=1,2$ we will mean property (P) where the
condition (\ref{iii}) has been replaced by (P$_i)\,i=1,2.$\newline
(P2)$\Rightarrow $(P1) is obvious.\newline
The implication (P1)$\Rightarrow $(P) follows from the unconditionality of
the basis $(e_n)_n$ .\newline
Indeed,given$(y_n)_{n\in E},(y_n^{*})_{n\in E}$ such that the condition (P1)
is fulfilled with$\left\| \underset{n\in E}{\sum y_n^{*}}\right\| \leq C,$we
set $P(x)=\underset{n\in E}{\sum }y_n^{*}(x)y_n.$In order to see that $%
\left\| P\right\| \leq C$ we observe that for $x\in X,\left\| x\right\| =1\;$%
there exists $\tilde{x}\in X$ such that $\left\| \tilde{x}\right\| =1$ and $%
y_n^{*}(\tilde{x})=\left| y_n^{*}(x)\right| $ for every $n$.Therefore for $%
x\in X,\left\| x\right\| =1$ we have that\newline
\begin{equation*}
\left\| P(x)\right\| =\left\| \underset{n\in E}{\sum }y_n^{*}(x)y_n\right\|
\leq \underset{n\in E}{\sum }\left| y_n^{*}(x)\right| =\underset{n\in E}{%
\sum }y_n^{*}(\tilde{x})=\left( \underset{n\in E}{\sum }y_n^{*}\right) (%
\tilde{x})\leq C
\end{equation*}
\end{remark}

\begin{lemma}
\label{PP}If the Banach space $(A,||$\thinspace $||_A)$ satisfies the
property $(P)$ then the same holds for the space $(X_A,||$ $||)$ .
\end{lemma}

\proof%
For every $x=\underset{\delta \in \frak{D}}{\sum }\lambda _{\delta
_{}}e_{\delta _{}}$ we set \newline
\begin{equation*}
\overline{x}=\overset{\infty }{\underset{i=1}{\sum }}\left( \underset{%
|\delta |=i}{\sum }|\lambda _{\delta _{}}|\right) e_i
\end{equation*}
\newline
Let $(x_n)_{n=1}^\infty =\left( \underset{|\delta |\in E_n}{\sum }\lambda
_{\delta _{}}e_{\delta _{}}\right) _{n\in \mathbb{N}}$ be a block sequence in $%
X_A$, where $E_1<E_2<...$ are successive subsets of $\mathbb{N}$,then $(%
\overline{x}_n)_{n=1}^\infty $ is a block sequence in $A$ such that $%
||x_n||=||\overline{x}_n||_A$, for every $n$.\newline
Since $A$ satisfies property $(P)$ for every $k>0$ there exists a sequence $(%
\overline{y}_n)_{n=1}^\infty =\left( \underset{i\in F_n}{\sum }\overline{%
\lambda }_i\overline{x}_i\right) _{n=1}^\infty $ of normalized blocks of $(%
\overline{y}_n)_{n=1}^\infty $ such that conditions (i),(ii),(iii) of
Definition \ref{DP} are fulfilled. Notice that if $\overline{y}_n^{*}=%
\underset{i\in F_n}{\sum }\left( \underset{j\in E_i}{\sum }m_je_j^{*}\right) 
$then $\overline{y}_n^{*}(\overline{y}_n)=\underset{i\in F_n}{\sum }%
\overline{\lambda }_i\left( \underset{j\in E_i}{\sum }m_j\left( \underset{%
|\delta |=j}{\sum }|\lambda _{\delta _{}}|\right) \right) $\newline
and $\overline{P}(\overline{x})=\underset{n\in E}{\sum }\overline{y}_n^{*}(x)%
\overline{y}_n.$\newline
We define 
\begin{equation*}
y_n=\underset{i\in F_n}{\sum }\overline{\lambda }_ix_i=\underset{i\in F_n}{%
\sum }\overline{\lambda }_i\left( \underset{|\delta |\in E_i}{\sum }\lambda
_{\delta _{}}e_{\delta _{}}\right) ,y_n^{*}=\underset{i\in F_n}{\sum }\;%
\underset{j\in E_i}{\sum }\left( \underset{|\delta |=j}{\sum }(sign\lambda
_{\delta _{}})m_je_{\delta _{}}^{*}\right)
\end{equation*}
\newline
where $sign\lambda =\left\{ 
\begin{array}{lll}
1 & \text{if} & \lambda >0 \\ 
-1 & \text{if} & \lambda <0 \\ 
0 & \text{if} & \lambda =0
\end{array}
\right. $\newline
Then $||\overline{y}_n^{*}||=\left\| \underset{i\in F_n}{\sum }\left( 
\underset{j\in E_i}{\sum }m_je_j^{*}\right) \right\| =||y_n^{*}||\leqslant 1$
, also 
\begin{equation*}
\overline{y}_n^{*}(y_n)=\underset{i\in F_n}{\sum }\underset{j\in E_i}{\sum }%
\overline{\lambda }_im_j\left( \underset{|\delta |=j}{\sum }|\lambda
_{\delta _{}}|\right) =\overline{y}_n^{*}(\overline{y}_n)>\frac 12
\end{equation*}
and \newline
\begin{equation*}
\left| \left| \underset{n\in E}{\sum }\theta _ny_n\right| \right| =\left\| 
\underset{n\in E}{\sum }\theta _n\overline{y}_n\right\| _A<\frac 1k
\end{equation*}
\newline
It remains to prove that $||P||\leqslant C_k$ where $P(x)=\underset{n\in E}{%
\sum }y_n^{*}(x)y_n$ for every $x\in X_A$.\newline
Let $x=\underset{\delta \in \frak{D}}{\sum }t_{\delta _{}}e_{\delta _{}}\in
B_{X_A}$ and $\overline{x}=\underset{i\in F_n}{\sum }\;\underset{j\in E_i}{%
\sum }\left( \underset{|\delta |=j}{\sum }|t_{\delta _{}}|\right) e_j$. Then
assuming that $(e_n)_{n=1}^\infty $is 1-unconditional 
\begin{equation*}
||P(x)||=\left| \left| \underset{n\in E}{\sum }\;\underset{i\in F_n}{\sum }\;%
\underset{j\in E_i}{\sum }\left( \underset{|\delta |=j}{\sum }sign\lambda
_{\delta _{}}t_{\delta _{}}m_j\right) y_n\right| \right| _A=
\end{equation*}
\newline
\begin{equation*}
\left| \left| \underset{n\in E}{\sum }\;\underset{i\in F_n}{\sum }\;%
\underset{j\in E_i}{\sum }m_j\left( \underset{|\delta |=j}{\sum }|t_{\delta
_{}}\right) y_n\right| \right| _A=||\underset{n\in E}{\sum }\overline{y}%
_n^{*}(\overline{x})y_n||=
\end{equation*}
\begin{equation*}
||\underset{n\in E}{\sum }\overline{y}_n^{*}(\overline{x})\overline{y}%
_n||_A=\left| \left| \overline{P}(\overline{x})\right| \right| \leq C
\end{equation*}

The following result is well known \cite{N} and says that if a closed set
almost absorbs the unit ball $B_{X\text{ }}$of the whole space then it also
absorbs $B_X.$We will state it in the following form:

\begin{lemma}
\label{LA}Let $W$ be a subset of a Banach space $X$ such that\newline
$B_X\subseteq \lambda W+\varepsilon B_X$ , for some $\lambda >0$ and $%
\varepsilon <1,$then $B_X\subset \dfrac \lambda {1-\varepsilon }\overline{W}$
\end{lemma}

\proof%
By our assumption\newline
$B_X\subseteq \lambda W+\varepsilon B_X\subset \lambda W+\varepsilon
(\lambda W+\varepsilon B_X)=\lambda (1+\varepsilon )W+\varepsilon ^2B_X$%
\newline
Inductively we prove that $B_X\subset \bigcap\limits_{n=1}^\infty \left(
\dfrac \lambda {1-\varepsilon }W+\varepsilon ^nB_X\right) $\newline
and by the completeness of the space this last set is equal to $\dfrac
\lambda {1-\varepsilon }\overline{W}$%
\endproof%

\begin{theorem}
\label{Th3.9}Let $(A,||$ $||_A)$ be reflexive Banach space with
unconditional basis $(e_n)_{n=1}^\infty $ which satisfies property $(P)$ and 
$\mathbf{a}=(a_n)_{n=1}^\infty =\left( \frac 1{nC_{n2^n}}\right)
_{n=1}^\infty $ where $(C_n)_{n=1}^\infty $ is the sequence defined by the
property $(P)$. Then $W$ is a $\mathbf{a}$-thin subset of $X_A$.
\end{theorem}

\begin{proof}
Suppose that $W$ is not a $\mathbf{a}-$thin subset of $X_A$. Then there
exists a normalized sequence $(x_n)_{n=1}^\infty $ of successive blocks of $%
(e_{\delta _{}})_{\delta \in \frak{D}}$ and $\lambda >0$ such that for every 
$n\in \mathbb{N}$\newline
\begin{equation}
B_X\subseteq \lambda (2^nW+2a_nB_{X_A})  \label{oe}
\end{equation}
where $X=\overline{span}(x_n)_{n=1}^\infty $.\newline
Fix a $k\in \mathbb{N}\mathbf{\,\;}$with $k>16\lambda .$\newline
By Lemma \ref{PP} $X_A$ also satisfies property $(P)$ for the same sequence $%
(C_k)_k$, hence there exists a normalized block sequence $(y_n)_{n=1}^\infty 
$ of $(x_n)_{n=1}^\infty $ such that for every $M$ infinite subset of $\mathbb{N%
}$ we can find $E\subset M$, $(\theta _n)_{n\in E}\subseteq \mathbb{R}^{+}$, $%
(y_n^{*})_{n\in E}\subseteq B_{X_A^{*}}$ such that 
\begin{equation}
\underset{n\in E}{\sum }\theta _n=1,\left| \left| \underset{n\in E}{\sum }%
\theta _ny_n\right| \right| <\frac 1{k2^k},  \label{o}
\end{equation}
\begin{equation}
\mathit{supp}y_n^{*}\subset \mathit{supp}y_n,y_n^{*}(y_n)>\frac 12
\label{oo}
\end{equation}
and 
\begin{equation}
||P||\leqslant C_{k2^k}\;\;\text{where\ }\;\;P(x)=\underset{n\in E}{\sum }%
y_n^{*}(x)y_n.  \label{ooo}
\end{equation}
\newline
Let $w_n\in W$ such that $||y_n-\lambda 2^kw_n||<a_k$ and set $E_n=\mathit{%
supp}y_n$.\newline
By proposition \ref{AE} there exists an infinite subset $M$ of $\mathbb{N}$ and
a family $(E_n)_{n\in M}$ of pairwise incomparable subsets of $\frak{D}$
such that for every $n\in \mathbb{N}$\newline
\begin{equation*}
D_n\subset F_n\;\;and\;\;||E_nw_n-F_nw_n||<\frac{a_k}{2^k\lambda }
\end{equation*}
Let $E\subset M$, $(\theta _n)_{n\in E}\subset \mathbb{R}^{+}$, $%
(y_n^{*})_{n\in E}\subset B_{X_A^{*}}$ to satisfy conditions (\ref{o}) , (%
\ref{oo}) , (\ref{ooo}) and set\newline
\begin{equation*}
z_n^{*}=y_n^{*}|_{D_n}\;\;\;,\;\;\;R(x)=\underset{n\in E}{\sum }z_n^{*}(x)y_n
\end{equation*}
\newline
It follows that 
\begin{equation}
\left| \left| R\right| \right| \leq \left| \left| P\right| \right| \leq
C_{k2^k}\;\;\;,\;\;\;\;|z_n^{*}(y_n)|>\frac 14.  \label{oooo}
\end{equation}
\newline
Set $Y_E=span[(y_n)_{n\in E}]$\newline
Then by (\ref{oe})\newline
$\;\;\;\;\;\;\;\;\;\;\;\;\;\;R(B_{Y_E})\subset \lambda 2^kR(W)+2\lambda
a_kR(B_{X_A})$\newline
and by (\ref{oooo}) and Lemma \ref{LA} we have that 
\begin{equation}
B_{Y_E}\subset 8\lambda 2^kR(W)  \label{99}
\end{equation}
Since the sets ($\mathit{supp}z_n^{*})_{n\in E}$ are pairwise incomparable
we get that 
\begin{equation}
R(W)\subset co[(y_i)_{i\in E}]  \label{66}
\end{equation}
\newline
Indeed,let $w=\underset{\gamma \in L}{\sum }\lambda _{\gamma _{}}x_{\gamma
_{}}\in W_0,$where $L$ is a set of branches of the tree such that $\left|
\left| x_{\gamma _{}}\right| \right| \leq 1$ for every $\gamma \in L$ and $%
\underset{\gamma \in L}{\sum }\left| \lambda _{\gamma _{}}\right| \leq 1.$We
set : 
\begin{equation*}
L_n=\left\{ \gamma \in L:\mathit{supp}z_n^{*}\cap \gamma \neq \emptyset
\right\} .
\end{equation*}
The sets $(L_n)_{n\in E}$ are pairwise disjoint and so 
\begin{equation*}
R(w)=\underset{n\in E}{\sum }z_n^{*}\left( \underset{\gamma \in L_n}{\sum }%
\lambda _{\gamma _{}}x_{\gamma _{}}\right) y_n,\text{ where }
\end{equation*}
\begin{equation*}
\underset{n\in E}{\sum }z_n^{*}\left( \underset{\gamma \in L_n}{\sum }%
\lambda _{\gamma _{}}x_{\gamma _{}}\right) \leq \underset{n\in E}{\sum }%
\left( \underset{\gamma \in L_n}{\sum }\left| \lambda _{\gamma _{}}\right|
\right) \leq \underset{\gamma \in L}{\sum }\left| \lambda _{\gamma
_{}}\right| \leq 1.
\end{equation*}
\newline
By (\ref{99}) and (\ref{66}) we get that 
\begin{equation*}
B_{Y_E}\subset 8\lambda 2^kco[(y_i)_{i\in E}]
\end{equation*}
But $k2^k\left( \sum_{n\in E}\theta _ny_n\right) \in B_{Y_E}$ and so $k\leq
8\lambda $ which is a contradiction since we selected $k>16\lambda .$%
\endproof%

\begin{remark}
\label{Rem3.2}If we assume that the sequence $(C_k)_{k=1}^\infty $ in the
above theorem is bounded then our arguments can show that the set $W$ is a
thin subset of $X_A$.
\end{remark}
\end{proof}

The following Lemma is well known (cf [LT] 1c8)

\begin{lemma}
\label{Lem3.10}Let $A$ be a reflexive Banach space with an unconditional
basis and $P:A\rightarrow Y$ a projection of $A$ onto a block subspace $Y$
of $A$ spanned by a normalized block sequence $(y_n)_{n=1}^\infty .$There
exists a projection $R:A\rightarrow Y$ of the form\newline
\begin{equation*}
R(x)=\underset{n=1}{\overset{\infty }{\sum }}y_n^{*}(x)y_n
\end{equation*}
\newline
where $(y_n^{*})_{n=1}^\infty $ is a uniformly bounded sequence of members
of $A^{*}$ with $supp(y_n^{*})\subset supp(y_n)$ for every $n\in \mathbb{N}%
\mathbf{.}$%
\endproof%
\end{lemma}

\begin{remark}
\label{Rem3.3}If a reflexive Banach space satisfies the property of the
above Lemma the same holds for the space $X_A.$
\end{remark}

The following result has been proved by N.Tomczak-Jaegermann who kindly
permit us to include it here.

\begin{proposition}
Let $A$ be a reflexive Banach space with an unconditional basis such that
every block subspace of $A$ has a block subspace complemented in $A.$ Then
the set $W$ is a thin subset of $X_A.$
\end{proposition}

\proof%
Suppose that the set $W$ is not a thin subset of $X_A.$Then there exists a
normalized block sequence $(x_n)_n$ of $X_A$ such that for every $%
\varepsilon >0$ there exists a $\lambda >0$ such that \newline
\begin{equation}
B_X\subset \lambda W+\varepsilon B_{X_A}  \label{T11}
\end{equation}
\newline
where $X=\overline{span}[(x_n)_n]$.From Lemma \ref{Lem3.10} and Remark \ref
{Rem3.3} we can suppose,by passing to a normalized block sequence of $%
(x_n)_n $ ,if necessary,that there exists a sequence $(x_n^{*})_n$ $\subset
X_A^{*}$ , uniformly bounded by a number $C>0$ ,such that for every $n\in 
\mathbb{N}$ ,$supp(x_n^{*})\subset range(x_n),$\newline
$x_n^{*}(x_n)=1$ and the operator\newline
\begin{equation}
P(x)=\underset{n=1}{\overset{\infty }{\sum }}x_n^{*}(x)x_n  \label{T12}
\end{equation}
\newline
is bounded.\newline
Let $\varepsilon =\dfrac 1{8\left\| P\right\| }$ and $\lambda >0$ such that
the relation (\ref{T11}) is satisfied.\newline
For every $n\in \mathbb{N}$ we select $w_n\in W_0$ such that $\left\|
x_n-\lambda E_nw_n\right\| <\varepsilon $\newline
where $E_n=range(x_n).$\newline
Applying Proposition \ref{AE} we can find an infinite subset $I$ of $\mathbb{N}$
and $(F_n)_{n\in I}$ pairwise incomparable subsets of $\frak{D}$ such that
for every $n\in I\;$\newline
$\,\,\,\,\,\,\,\,\,\,\,\,\,\,\,\,\,\,\,\,\,\,\,\,\,\,\,\,\,\,\,\,\,\,\,\,\,%
\,\,\,\,\,\,\,\,\,\,\,\,\,\,\,\,\,F_n\subset E_n$ and $\left\|
E_nw_n-F_nw_n\right\| <\frac 1{4\lambda }$\newline
We set \newline
\begin{equation}
y_n^{*}=F_nx_n^{*}\;\;\;\;\text{and\ \ \ \ \ }R(x)=\underset{n=1}{\overset{%
\infty }{\sum }}y_n^{*}(x)x_n  \label{T14}
\end{equation}
\newline
We observe that \newline
\begin{equation}
\left\| R\right\| \leq \left\| P\right\| \text{ \ \ \ \ \ and\ \ \ \ \ \ \ \ 
}\left| y_n^{*}(x_n)\right| >\frac 14  \label{T15}
\end{equation}
\newline
Let $\;\;\;\;\;\;\;\;\;\;\;\;\;\;\;\;\;\;\;\;\;\;Y=\overline{span}[%
(x_n)_{n\in I}]$\newline
By Lemma \ref{LA} , (\ref{T14}) and (\ref{T15})\newline
we get that $B_Y\subset 8\lambda \overline{R(W)}$\newline
and since the sets $(suppy_n^{*})_{n\in I}$ are pairwise incomparable we
have that\newline
$R(W)\subset C\overline{co}[(x_n)_{n\in I}]$ and so \newline
\begin{equation}
B_Y\subset 8\lambda C\overline{co}[(x_n)_{n\in I}]  \label{T16}
\end{equation}
\newline
This leads to a contradiction by a theorem of Linderstrauss and Phelps (\cite
{LP}) which asserts that a convex body in a reflexive Banach space cannot
have countable many extreme points.\newline
We can also derive a contradiction by selecting a combination $z=\underset{}{%
\underset{n\in E}{\sum }\theta _nx_n}\;\;,\;\;\underset{n\in E}{\sum }\left|
\theta _n\right| =1$ such that $\left\| z\right\| <\frac 1{16\lambda }$
which contradicts relation (\ref{T16})%
\endproof%

\begin{definition}
Let X be a Banach space with a basis $(e_n)_{n=1}^\infty $ and $1<p\leq
\infty .$ We say that $X$ {\bf has an upper }$p-${\bf estimate }if
there exists $C>0$ such that for every normalized block sequence $(x_n)_n,$
\end{definition}

\begin{center}
$\left\| \underset{i=1}{\overset{n}{\sum x_i}}\right\| \leq Cn^{\frac 1p}$
\end{center}

\begin{remark}
\label{3}If a space $X\,$has an upper-p estimate for some $p>1$ then it
satisfies property (P2) and hence property (P)(see Remark \ref{Rem3.1}
).Spaces with an upper $p-$estimate are $L^p,1<p<\infty $ and Banach spaces
where $\ell ^1$ is not uniformly block representable (\cite{MS}).It follows
that for every such reflexive space with an unconditional basis has the
property that the set $W$ in the space $X_A$ is $\mathbf{a-}$thin for an
appropriate sequence $\mathbf{a}$ depending on $p$.\newline
It can be shown for $A$ the Tsirelson's space T,its dual ,Schlumprecht's
space $S$ and other spaces the set $W$ in $X_A$ is also $\mathbf{a-}$thin .
\end{remark}

We pass now to show that certain classes of operators are thin operators .As
we have pointed out every thin operator is strictly singular but the
converse is not always true. However,if we restrict our interesting in the
classes $\mathcal{L}\left( \ell ^p,\ell ^r\right) $ $1\leq p,r<\infty $ then
the subsets of strictly singular operators and thin operators coincide. This
follows from the next lemmas:

\begin{lemma}
\label{3.11}Let $T$ be in $\mathcal{L}\left( X,\ell ^p\right) $ $1\leq
p<\infty $ which is not a thin operator. Then there exists $Z$ block
subspace of $\ell ^p$ and a projection $P:\ell ^p\rightarrow Z$ with $||P||=1
$ such that $P\circ T$ is a quotient map.
\end{lemma}

\proof%
Choose any $\varepsilon <1.$ then there exists $Z$ infinite dimensional
closed subspace of $\ell ^p$ and $\lambda \in \mathbb{R}^{+}$ such that $%
B_Z\subset \lambda T[B_X]+\varepsilon B_{\ell ^p}.$ Further we may assume
that $Z$ is a block subspace of $\ell ^p$ hence there exists $P:\ell
^p\rightarrow Z$ a projection with $||P||=1$. Therefore $B_Z\subset \lambda
P\left( T[B_X]\right) +\varepsilon B_Z$ and from Lemma \ref{LA} $B_Z\subset
\dfrac \lambda {1-\varepsilon }P\circ T[B_X]$ which implies that $P\circ T$
is a quotient map 
\endproof%

\begin{lemma}
\label{3.12}If $T:\ell ^p\rightarrow \ell ^p$ $1\leq p<\infty $ is a bounded
linear operator then there exists $Z$ an infinite dimensional closed
subspace of $\ell ^p$ such that $T|_Z$ is an isomorphism
\end{lemma}

\proof%
Since the operator $T$ is onto there exists $C\geq 1$ and a sequence $\left(
z_n\right) _{n\in \mathbb{N}}$ such that $Tz_n=e_n.$ Further by passing to a
subsequence we may assume that $\left( z_n\right) _n$ is a block basic
sequence in $\ell ^p$. Therefore 
\begin{eqnarray*}
\left( \sum_{i=1}^n|a_i|^p\right) ^{\dfrac 1p} &\leq &\left\|
\sum_{i=1}^na_ie_i\right\| =\left\| \sum_{i=1}^na_iTz_i\right\| \leq \\
&\leq &||T||\;\left\| \sum_{i=1}^na_iz_i\right\| \leq ||T||\;C\;\left(
\sum_{i=1}^n|a_i|^p\right) ^{\dfrac 1p}
\end{eqnarray*}
and the proof is complete%
\endproof%

\begin{corollary}
\label{3.13}(a) If 1$\leq p\neq r<\infty $ then every $T\in $ $\mathcal{L}%
\left( \ell ^p,\ell ^r\right) $ is a thin operator.\newline
(b) If 1$\leq p\neq r<\infty $ and $T\in $ $\mathcal{L}\left( \ell ^p,\ell
^p\right) $ is a strictly singular operator then $T$ is a thin operator.
\end{corollary}

\proof%
It follows from Lemmas \ref{3.11} and \ref{3.12} 
\endproof%

\begin{remark}
It follows readily that Lemma \ref{3.11} and Lemma \ref{3.12} remain true if
instead of $\ell ^p$ we consider the space $c_0$. Therefore we have similar
results of the above Corollary if either $\ell ^p$ or $\ell ^r$ is
substituted by the space $c_0\left( \mathbb{N}\right) $
\end{remark}

\section{{\bf THIN NORMING SETS II}}

\begin{center}
{\bf (general reflexive case)}
\end{center}

In this section we prove that every reflexive space $A$ with an
unconditional basis has a subspace $B$ such that the set $W$ in $X_B$ is an $%
\mathbf{a}-$thin set for an appropriate null sequence $\mathbf{a}$. The
proof heavily depends on the uniform control on a subspace of the rate of
convergence in norm to zero of certain convex combination. Therefore the
ordinal order of the complexity of weakly null sequences seems to be
necessary and we strongly use the results from \cite{AMT}.\newline
The main parts of this section are the following:\newline
The first one contains the definitions of Schreier families, initially
introduced in \cite{AA}, and R.A.-Hierarchy from \cite{AMT}. We also present
some results from \cite{AMT} that we will use here. The most important of
them are the dichotomy principle (Proposition \ref{Prop4.01}) and the large
families theorem (Theorem \ref{Th4.01}).\newline
In the second part for a countable ordinal $\xi $, a natural number $n$ and $%
B$ a subspace of $A$ we introduce a numerical quantity $\tau _{\xi ,n}(B)$
which is basic for our proof. We show that there exists a countable ordinal $%
\xi _0$ such that $\tau _{\xi ,n}(B)=0$ for every $\xi \geq \xi _0$ and $%
n\in \mathbb{N}$ and there exists $B$ subspace of $A$ on which $\tau _{\xi ,n}$
is stabilized for any further subspace. This subspace $B$ is the desired
subspace of $A$. For this, we need to show that for every $\varepsilon >0$
there exists $\left( \xi ,n\right) $ such that $0<\tau _{\xi
,n}(B)<\varepsilon $.\newline
Finally, in the last part, (Propositions \ref{Prop4.04} , \ref{Prop4.05} ,%
\ref{Th4.04}), we prove that the set $W$ in $X_B$ is an $\mathbf{a-}$thin
set.\newline
{\bf The Schreier Families}

For every countable ordinal $\xi $ we define a family $S_{\xi _{}}$ of
finite subsets of $\mathbb{N}$ as follows:

\begin{description}
\item[1.]  $S_0=\left\{ \{n\}:n\in \mathbb{N}\right\} \cup \{\emptyset \}$

\item[2.]  If $S_{\xi _{}}$ has been defined then 
\begin{equation*}
S_{\xi +1}=\left\{ \underset{i=1}{\overset{n}{\bigcup }}F_i:n\in \mathbb{N},%
\text{ }n\leq F_1<F_2<...<F_n,\text{ }F_i\in S_{\xi _{}}\right\} \cup
\{\emptyset \}
\end{equation*}

\item[3.]  If $\xi $ is a limit ordinal and $S_{\zeta _{}}$has been defined
for every $\zeta <\xi $ we fix a strictly increasing sequence $(\xi
_n)_{n\in \mathbb{N}}$ of non-limit ordinals with $\underset{n\in \mathbf{N}}{%
\sup }\xi _n=\xi $ and we define 
\begin{equation*}
S_{\xi _{}}=\{F:n\leq \min F\;\;\text{and\ \ }F\in S_{\xi _n\;}\;\text{for\
\ some\ \ \ }n\in \mathbb{N}\}
\end{equation*}
{\bf The Repeated Averages Hierarchy or RA-Hierarchy}
\end{description}

We denote by $S_{\ell _1}^{+}$ the positive part of the unit sphere of $\ell
_1(\mathbb{N})$, i.e. the set of all sequences $(a_n)_{n=1}^\infty $ of
non-negative numbers such that $\underset{n=1}{\overset{\infty }{\sum }}%
a_n=1 $, set also $(e_n)_{n=1}^\infty $ the usual basis of $\ell _1(\mathbb{N})$%
.\newline
For every countable ordinal $\xi $, every infinite subset $M$ of $\mathbb{N}$
there exists a sequence $(\xi _n^M)_{n\in \mathbb{N}}$ satisfying the following
properties:

\begin{itemize}
\item[1.]  $\xi _n^M\in S_{\ell _1}^{+}\;$for\ every$\ n\in \mathbb{N}$ and $%
\xi _n^M(m)\leq \frac 1{minM}$ for every $M\in [\mathbb{N]},n,m\in \mathbb{N},\xi
>0.$

\item[2.]  $\mathit{sup}p\xi _n^M<\mathit{sup}p\xi _{n+1}^M\;\;_{}$and$%
\;\;_{}\mathit{sup}p\xi _n^M\in S_{\xi _{}}$for every $n\in \mathbb{N}$
\end{itemize}

\begin{enumerate}
\item[3.]  $M=\underset{n=1}{\overset{\infty }{\bigcup }}\mathit{sup}p\xi
_n^M$
\end{enumerate}

\begin{itemize}
\item[4.]  For every $M$ infinite subset of $\mathbb{N}$ and $%
(n_k)_{k=1}^\infty $ a strictly increasing sequence if we set 
\begin{equation*}
M^{\prime }=\underset{k=1}{\overset{\infty }{\bigcup }}\mathit{sup}p\xi
_{n_k}^M\text{ then }\xi _k^{M^{\prime }}=\xi _{n_k}^M
\end{equation*}
\end{itemize}

We recall the definition of $(\xi _n^M)_{n\in \mathbb{N}}$ from \cite{AMT}. It
is given by induction

\begin{enumerate}
\item  For $\xi =0$, $M=\{m_1,m_2,...\}$, $m_1<m_2<...$ we set 
\begin{equation*}
\xi _n^M=e_{m_n}
\end{equation*}

\item  If $\xi =\zeta +1$ and $(\zeta _n^M)_{n\in \mathbb{N}}$ has been defined
for all $M\in [\mathbb{N}]$.We set 
\begin{equation*}
M=(m_n)_{n=1}^\infty \text{ and }\xi _1^M=\frac{\zeta _1^M+...+\zeta
_{m_{m_1}}^M}{m_{m_1}}
\end{equation*}
and then inductively we define 
\begin{equation*}
\xi _n^M=\xi _1^{M_n}
\end{equation*}
where 
\begin{equation*}
M_1=M,\;\;\;\;\;\;\;\;M_n=M\backslash \underset{i=1}{\overset{n-1}{\bigcup }}%
\mathit{sup}p\xi _i^M
\end{equation*}

\item  If $\xi $ is a limit ordinal we consider the sequence $(\xi _n)_{n\in 
\mathbb{N}}$ of successors ordinals that defines $S_{\xi _{}}$. For any $M\in [%
\mathbb{N}]$ $\;\;M=\{m_1,m_2...\}$we set$\;n_1=m_1$ 
\begin{equation*}
\xi _1^M=(\xi _{n_1})_1^M
\end{equation*}
and inductively $M_k=M\backslash $\textit{sup}$p[\xi
_{n_{k-1}}^{}]_1^{M_{k-1}},n_k=\min M_k$%
\begin{equation*}
\xi _k^M=(\xi _{n_k})_1^{M_k}
\end{equation*}
\newline
By standard induction one could verify the properties (1)-(4) listed above.
\end{enumerate}

The following dichotomy principle is important for our proofs ; it permit us
to stabilize functions acting on the summability methods $(\xi _n^M)_{n\in 
\mathbb{N}}$ without loosing the order $\xi .$It is a consequence of stability
property (4) of the repeated averages hierarchy and the following
combinatorial result

\begin{theorem}[Nash-Williams,Galvin and Prikry,Silver,Ellentuck ]
\label{Ra}Let $A$ be an analytic subset of $[\mathbb{N]}$Then for every $M\in [%
\mathbb{N]}$ there exists $L\in [M]$ such that either $[L]\subset A$ or $%
[L]\subset \mathbb{N}|\backslash A.$
\end{theorem}

(for a proof we refer to \cite{Ke} or \cite{T} )

\begin{proposition}[Dichotomy Principle]
{\bf \label{Prop4.01} }Let $\xi <\omega _1$, $M\in [\mathbb{N}]$ and $n\in 
\mathbb{N}.$ We set: 
\begin{equation*}
I_{M,n}=\{(\xi _{k_1}^N,...,\xi _{k_n}^N):k_1<...<k_n\;,N\in [M]\}.
\end{equation*}
Clearly $I_{M,n}\subset \left( S_{\ell _1}^{+}\right) ^n.$\newline
Then for every real function $\phi :I_{M,n}\rightarrow \mathbb{R}$ and $\tau
\in \mathbb{R}$ one of the following holds:
\end{proposition}

\begin{itemize}
\item[(i)]  \textit{There exists }$L\in [M]$\textit{\ such that }$\phi
(s)>\tau $\textit{\ for all }$s\in I_{L,n}$

\item[(ii)]  \textit{There exists }$L\in [M]$\textit{\ such that }$\phi
(s)\leq \tau $\textit{\ for all }$s\in I_{L,n}$
\end{itemize}

\proof%
Consider the following partition of $M.$\newline
$A=\{N\in [M]:\phi ((\xi _1^N,...,\xi _n^N))>\tau \},B=[M]\backslash A.$%
\newline
It is clear that $A$ is an open subset of $[M]$ hence from Theorem \ref{Ra}
we get that either there exists $L\in [M]$ such that $[L]\subset A$ or $_{}$%
there exists $L\in [M]$ such that $[L]\subset B.$\newline
Assume that $_{}[L]\subset A$ .Since (by property (4) of the
R.A.Hierarchy)for $N\in [L]$ and $k_1<...<k_n$ there exists $N^{\prime }\in
[L]$ such that $\xi _i^{N^{\prime }}=\xi _{k_i}^N$ for $i=1,...n$ we get
that for every $N\in [L]$ and $k_1<...<k_n$, $\phi ((\xi _{k_1}^N,...,\xi
_{k_n}^N))=\phi ((\xi _1^{N^{\prime }},...,\xi _n^{N^{\prime }}))>\tau $ and
case (i) is satisfied.\newline
If $_{}[L]\subset B$ we obtain case (ii).%
\endproof%

{\bf Notation} Let $M\in [\mathbb{N}]$ and $\xi $ a countable ordinal .We
introduce two more families related to $S_{\xi _{}}$ and the set $M$.\newline
If $M=(n_k)_{k\in \mathbb{N}}$ we set\newline
\begin{equation*}
S_{\xi _{}}^M=\{G\subset M:\text{there exists }F\in S_{\xi _{}}\;\text{such
that }G=\{n_\ell :\ell \in F\}\}
\end{equation*}
\newline
We also denote by 
\begin{equation*}
S_{\xi _{}}[M]=\{F\subset M:F\in S_{\xi _{}}\}
\end{equation*}
\begin{equation*}
\lbrack S_{\xi _{}}^M]^n=\left\{ \underset{i=1}{\overset{n}{\bigcup }}%
U_i:U_i\in S_{\xi _{}}^M,U_1<U_2<...<U_n\right\}
\end{equation*}
Let us observe that by the definition of the families $S_{\xi _{}}$, $S_{\xi
_{}}^M\subset S_{\xi _{}}[M]$ but the converse is not true.\newline
The following lemma is due to G. Androulakis and T. Odell and shows that by
passing to a subset $N$ of $M$ we have that the members of $S_{\xi _{}}[N]$, 
$S_{\xi _{}}^M$ are almost equal.

\begin{lemma}
\cite{AO}\label{Lem4.01}. For every $\xi <\omega _1$ and $M\in [\mathbb{N}]$
there exists $N\in [M]$ such that for every $F\subset N$ if $F\in S_{\xi
_{}}[N]$ then $F\backslash \{\min F\}\in S_{\xi _{}}^M.$%
\endproof%
\end{lemma}

(For a proof see also \cite{AMT}).

\begin{definition}
A family $\frak{F}$ of subsets of $\mathbb{N}$ is said to be
\end{definition}

\begin{description}
\item[(a)]  {\bf adequate}\textit{\ if }$F\in \frak{F}$\textit{\ and }$%
G\subset F$\textit{\ then }$G\in \frak{F}$

\item[(b)]  {\bf compact }\textit{if it is compact in the topology of
pointwise convergence as a subspace of }$\{0,1\}^{\mathbb{N}}$
\end{description}

\begin{definition}
A compact adequate family $\frak{F}$ of subsets of $\mathbb{N}$ is said to be
\end{definition}

\begin{description}
\item[(a)]  $(\xi ,M,\varepsilon )$\textit{-{\bf large}, where }$\xi
<\omega _1$\textit{, }$M\subset [\mathbb{N}]$\textit{\ and }$\varepsilon >0$%
\textit{\ if for every }$L\in [M]$\textit{\ and }$n\in \mathbb{N}$\textit{\
there exists }$F\in \frak{F}$\textit{\ such that\newline
} 
\begin{equation*}
<\xi _n^L,F>=\underset{k\in F}{\sum }\xi _n^L(k)>\varepsilon 
\end{equation*}

\item[(b)]  $(n,\zeta ,M,\varepsilon )-${\bf \ large}\textit{\ if for
every }$N\in [M]$\textit{, } 
\begin{equation*}
\underset{F\in \mathcal{F}}{sup}\;\underset{i}{inf}\{<\zeta
_i^N,F>\}_{i=1}^n>\varepsilon 
\end{equation*}
\end{description}

The following result is proved in \cite{AMT} and they will be fundamental
ingrendients for our proofs

\begin{theorem}
\label{Th4.01}Let $\frak{F}$ be an adequate compact family of finite subsets
of $\mathbb{N}$ such that for some countable ordinal $\xi $, $M$ infinite
subset of $\mathbb{N}$ and $\varepsilon >0$, $\frak{F}$ is $(\xi ,M,\varepsilon
)$-large($(n,\xi ,M,\varepsilon )$-large for some $n\in \mathbb{N}$).Then there
exists an infinite subset $L$ of $M$ such that $S_{\xi _{}}^L\subset \frak{F}
$ ( $[S_{\zeta _{}}^L]^n$ $\subset \frak{F}\mathcal{)}$.
\end{theorem}

\proof%
For $n=1$ the result follows from Proposition 2.32 and Theorem 2.26 of\cite
{AMT} ;for the general $n$ the result follows from Lemma 2.3.5 in \cite{AMT}
and the Remark after Corollary 2.2.8 in \cite{AMT}.%
\endproof%

From Lemma \ref{Lem4.01} and Theorem \ref{Th4.01} we get the following
result.

\begin{corollary}
\label{Cor4.01}Let $\frak{F}$ be a $(\xi ,M,\varepsilon )$-{\bf large
family;}for any $\delta >0$ there exists a $L\in [M]$ such that for every $%
N\in [L]$ and $n\in \mathbb{N}$ there exists $F\in \frak{F}$ with the property $%
\left\langle \xi _n^N,F\right\rangle >1-\delta .$
\end{corollary}

{\bf Notation} .Let $X$ be a Banach space with a basis. We denote by:

\begin{enumerate}
\item[1.]  $\Sigma (X)\;$the set of all normalized block sequences of basis
i.e.. 
\begin{equation*}
\Sigma (X)=\{(x_i)_{i=1}^\infty :x_i<x_{i+1},||x_i||=1\;\text{and\ }x_i\in
X\}
\end{equation*}

\item[2.]  
\begin{equation*}
\Sigma (\mathbb{N})=\{(F_i)_{i=1}^\infty :F_i<F_{i+1},F_i\subset \mathbb{N}\}
\end{equation*}

\item[3.]  For every $\mathcal{F}=(F_i)_{i=1}^\infty \in \Sigma (\mathbb{N})$
and $s=(y_i)_{i=1}^\infty \in \Sigma (X)\;$we set$\;\;$\newline
\begin{equation*}
Y_i(\mathcal{F},s)=\left\langle (y_j)_{j\in F_i}\right\rangle \newline
\newline
\end{equation*}
\begin{equation*}
\Sigma (s,\mathcal{F})=\{(z_i)_{i=1}^\infty :z_i\in Y_i(\mathcal{F}%
,s),\left\| z_i\right\| =1\text{ for every }i\}.\newline
\end{equation*}

\item[4.]  For $s=(x_k)_{k\in \mathbb{N}}\in \Sigma (X)\;$we denote by 
\begin{equation*}
\xi _n^M\cdot s=\underset{k=1}{\overset{\infty }{\sum }}\xi _n^M(k)x_k
\end{equation*}
\newline
and

\item[5.]  further,if $\xi $ is a countable ordinal, $M\in [\mathbb{N}],$ $n\in 
\mathbb{N}$ we denote by 
\begin{equation*}
||\xi _n^M(s,\mathcal{F})||=max\left\{ ||\xi _n^M\cdot t||,t\in \Sigma (s,%
\mathcal{F})\right\} .
\end{equation*}
and

\item[6.]  
\begin{equation*}
\left\| \frac{\xi _1^N+...+\xi _n^N}n(s,\mathcal{F})\right\| =max\left\{
\left\| \frac{\xi _1^N\cdot t+...+\xi _n^N\cdot t}n\right\| ,t\in \Sigma (s,%
\mathcal{F})\right\} 
\end{equation*}

\item[7.]  For $s_1,s_2\in \Sigma (X)$ we write $s_2\prec s_1$ if $s_2$ is a
normalized block sequence of $s$

\item[8.]  If $s_1=(x_i)_{i=1}^\infty $ and $s_2=(x_i)_{i=n}^\infty $ for
some $n\in \mathbb{N}$ we say that $s_{2\text{ }}$is a {\bf tail subsequence 
}of $s_{1\text{.}}$

\item[9.]  If $Y_1=\overline{span}[(x_i)_{i=1}^\infty ]$ and $Y_2=\overline{%
span}[(x_i)_{i=n}^\infty ]$ for some $n\in \mathbb{N}$ we say that $Y_{2\text{ }%
}$is a {\bf tail subspace} of $Y_{1\text{.}}$
\end{enumerate}

\begin{definition}
\label{Def4.01}For every $s\in \Sigma (X)$, $\xi <\omega _1$,we set
\end{definition}

\begin{itemize}
\item[(a)]  $\;\;\;\;\;\;\;\;\;\;\;\;\;\;\;\;\;\;\tau _{\xi _{}}(s)=%
\underset{\mathcal{F}\in \Sigma (\mathbf{N})}{\sup }\;\underset{M\in [%
\mathbf{N}]}{\sup }\;\underset{N\in [M]}{\inf }||\xi _1^N(s,\mathcal{F})||$

\item[(b)]  $\;\;\;\;\;\;\;\;\;\;\;\;\;\;\;\;\;\;\;\;\;\;\;\;\;\;\;\;\;\;\;%
\tau _{\xi _{}}(X)=\underset{s\in \Sigma (X)}{\inf }\tau _{\xi _{}}(s).$
\end{itemize}

This is a complicated definition and we will attempt to explain its
necessity after Lemma 4.6.

\begin{remark}
\begin{enumerate}
\item  We observe that 
\begin{equation*}
\underset{N\in [M]}{inf}||\xi _1^N(s,\mathcal{F})||=\underset{N\in [M]}{\inf 
}\;\underset{n\in \mathbf{N}}{\inf }||\xi _n^N(s,\mathcal{F})||
\end{equation*}
\newline
This follows from the stability property (4) of $(\xi _n^N)_{n\in \mathbb{N}}.$

\item  $0\leq \tau _{\xi _{}}(s)\leq 1$ for every $\xi <\omega _1$ and $s\in
\Sigma (X)$
\end{enumerate}
\end{remark}

\begin{lemma}
\label{Lem4.5}Suppose that for some $s\in \Sigma (X)$, $\xi <\omega _1$ we
have that $\tau _{\xi _{}}(s_{})=\delta >0.$Then the following properties
hold:
\end{lemma}

\begin{enumerate}
\item  \textit{For every }$\mathcal{F}=(F_i)_{i=1}^\infty \in \Sigma (\mathbb{N}%
)$\textit{\ and }$M\in [\mathbb{N}\mathbf{]}$\textit{\ there exists }$L\in [M]$%
\textit{\ satisfying:For every }$N\in [L]$\textit{\ and }$n\in \mathbb{N}%
\mathbf{\;}$%
\begin{equation*}
||\xi _n^N(s,\mathcal{F})||<2\delta 
\end{equation*}

\item  \textit{There exists }$\mathcal{F}=(F_i)_{i=1}^\infty \in \Sigma (%
\mathbb{N})$\textit{\ and }$M\in [\mathbb{N}\mathbf{]}$\textit{\ such that for all 
}$N\in [M\mathbf{]}$\textit{\ }$\mathit{,}n\in \mathbb{N}\mathbf{\;}$%
\begin{equation*}
||\xi _n^N(s,\mathcal{F})||>\frac 12\delta 
\end{equation*}

\item  \textit{There exists }$\mathcal{F}\in \Sigma (\mathbb{N}),M\in [\mathbb{N}%
\mathbf{]}$\textit{\ such that for all }$N\in [M\mathbf{]}$\textit{\ ,}$n\in 
\mathbb{N}$%
\begin{equation*}
\frac 12\delta <||\xi _n^N(s,\mathcal{F})||<2\delta 
\end{equation*}
\end{enumerate}

\proof%
Property (1)follows by Proposition \ref{Prop4.01}(the Dichotomy Principle)
and property (2) follows from the definition .Property (3) is a combination
of (1) and (2).%
\endproof%

To explain the nature of the quantity $\tau _{\xi _{}}\left( s\right) $ let
us assume that $\tau _{\xi _{}}\left( s\right) =\delta >0$ . Then from the
previous Lemma \ref{Lem4.5} (2) we have that there exists $\mathcal{F}\in
\Sigma \left( \mathbb{N}\right) $ and $M\in [\mathbb{N]}$ such that $||\xi
_n^N\left( s, \mathcal{F} \right) ||>\frac \delta 2$ for every $N\in [M]$ .It follows
readily for every $N\in [M]$ there exists $t^N\in \Sigma \left( s,\mathcal{F}%
\right) $ such that $||\xi _n^N\left( t^N\right) ||>\frac \delta 2$ for
every $n\in \mathbb{N}$ . The difficulty is that by changing from $N$ to $L$ we
do not know that $t^N$ remains equal to $t^L$ .This forces us to be more
careful in our proofs. The advantage of the above definition comes
essentially from part $\left( 1\right) $ of Lemma \ref{Lem4.5} . Indeed if
,for a given $\varepsilon >0$ , we know that $0<\tau _{\xi _{}}\left(
s\right) =\delta <\varepsilon $ then for every $\mathcal{F}\in \Sigma \left( 
\mathbb{N}\right) $ and every $M\in [\mathbb{N]}$ there exists $L\in [M]$ such
that $||\xi _n^N\left( s,\mathcal{F}\right) ||<2\delta <2\varepsilon $ for
all $N\in [L],n\in \mathbb{N}$ .This permits us to have a uniform $2\delta $%
-control of all convex combinations of a certain form and at the same time
we know by part two of Lemma \ref{Lem4.5} that there are convex combination
of the same form remaining greater than $\frac \delta 2$ .This key
observation is crucial for our approach.

\begin{lemma}
\label{Lem4.06}Let $\xi <\omega _1$
\end{lemma}

\begin{description}
\item[(a)]  \textit{If }$s_1,s_2\in \Sigma (X)$\textit{\ and }$s_2\prec s_1$%
\textit{\ then }$\tau _{\xi _{}}(s_1)\geq \tau _{\xi _{}}(s_2)$

\item[(b)]  \textit{If }$Y\prec Z$\textit{\ then }$\tau _{\xi _{}}(Y)\geq
\tau _{\xi _{}}(Z)$

\item[(c)]  \textit{If }$s_1,s_2\in \Sigma (X)$\textit{and }$s_2$\textit{\
is a tail subsequence of }$s_1$\textit{\ for some }$n\in \mathbb{N}$\textit{\
then }$\tau _{\xi _{}}(s_1)=\tau _{\xi _{}}(s_2)$

\item[(d)]  \textit{If }$Z$\textit{\ is a tail subspace of }$Y$\textit{\
then }$\tau _{\xi _{}}(Y)=\tau _{\xi _{}}(Z)$\textit{.}
\end{description}

\proof%
(a) It is enough to show that if $\tau (s_2)>\theta $ then $\tau
(s_1)>\theta .$\newline
Suppose that $\tau (s_2)>\theta $ for some $\theta \in (0,1)$ then there
exists $\mathcal{F=}(F_n)_{n=1}^\infty \in \Sigma (\mathbb{N})$and \newline
$M\in [\mathbb{N}\mathbf{]}$ such that for every $N\in [M]$ we can select $%
t^N\in \Sigma (s_2,\mathcal{F})$ with $||\xi _n^N\cdot t^N||>\theta .$%
\newline
Let $s_1=(y_n)_{n=1}^\infty $ and $s_2=(z_n)_{n=1}^\infty $with $%
z_n=\sum_{i\in E_n}\lambda _iy_i$ where $(E_1)_{n=1}^\infty $ is a sequence
of successive subsets of $\mathbb{N}$.\newline
We set 
\begin{equation*}
F_n^{\prime }=\bigcup\limits_{i\in F_n}E_i\;\;\;\;\;\mathcal{F}^{\prime
}=(F_n^{\prime })_{n=1}^\infty
\end{equation*}
then $\Sigma (s_2,\mathcal{F}^{\prime })\subset \Sigma (s_1,\mathcal{F}%
^{\prime })$ and so for every $N\in [M]$ we select the above sequence $%
t^N\in \Sigma (s_1,\mathcal{F}^{\prime })$ with $||\xi _n^N\cdot
t^N||>\theta .$\newline
The cases (b) ,(c),(d) are obvious.%
\endproof%

\begin{lemma}
\label{Lem4.07}For every $\xi <\omega _1$ there exists a block subspace $Y$
of $X\,$ such that $\tau _{\xi _{}}(s)=\tau _{\xi _{}}(Y)$ for every $s\in
\Sigma (Y).$
\end{lemma}

\proof%
We inductively define sequences $(s_i)_{i=1}^\infty ,(Y_i)_{i=1}^\infty $ of
normalized block sequences of $Z$ ,block subspaces of $Z,$ respectively ,as
follows:

For $i=1$ we set $s_1=(e_n)_{n=1}^\infty \;\;,\;\;\;\;Y_1=Z$

If $Y_k,s_k$ have been defined we select 
\begin{equation}
s_{k+1}\prec s_k\;\;\;such\ \ that\;\;\tau _{\xi _{}}(s_{k+1})-\tau _{\xi
_{}}(Y_k)<\dfrac 1{k+1}  \label{rel111}
\end{equation}
and we define $Y_{k+1}=\overline{span}[s_{k+1}]$\newline
Let $s$ be a diagonal block sequence of $(s_i)_{i=1}^\infty $ and set $Y=%
\overline{span}[s].$\newline
If $t\prec s$ for every $k$ there exists a $t_k\prec s_k$ which is a tail
subsequence of $t$ so by Lemma \ref{Lem4.06} (a) and (c) we get that 
\begin{equation*}
\tau _{\xi _{}}(Y_k)\leq \tau _{\xi _{}}(t_k)=\tau _{\xi _{}}(t)=\tau _{\xi
_{}}(t_{k+1})\leq \tau _{\xi _{}}(s_{k+1})
\end{equation*}
\newline
and therefore by the relation (\ref{rel111}) above ,we get that $\tau _{\xi
_{}}(s)=\tau _{\xi _{}}(t)$%
\endproof%

The following definition is the natural extension of Definition \ref{Def4.01}
for $n-$averages of $(\xi _k^N)_{k\in \mathbb{N}}$ which is going to use in a
certain stage of our proof.

\begin{definition}
\label{4.4}(a) For every $s\in \Sigma (X)$, $\xi <\omega _1$,$n\in \mathbb{N}$
we set:\newline
\begin{equation*}
\tau _{\xi ,n}(s)=\underset{\mathcal{F}\in \Sigma (\mathbf{N})}{\sup }\ 
\underset{M\in [\mathbf{N}]}{\sup }\ \underset{N\in [M]}{\inf }\;||\frac{\xi
_1^N+...+\xi _n^N}n(s,\mathcal{F})||
\end{equation*}
\end{definition}

(b) $\tau _{\xi ,n}(X)=\underset{s\in \Sigma (X)}{\inf }\tau _{\xi ,n}(s)$

\begin{remark}
\begin{itemize}
\item[(a)]  For $n=1$ the previous Definition coincides with Definition \ref
{Def4.01}

\item[(b)]  The corresponding results of Lemma \ref{Lem4.06} and \ref
{Lem4.07} stated of $\tau _{\xi ,n}$ instead of $\tau _{\xi _{}}$ also hold
and the proofs are completely analogous .In particular for a given $\xi
<\omega _1$ ,$n\in \mathbb{N}$ $,$ there exists a block subspace $Y$ of $X\,$
such that $\tau _{\xi ,n}(Y)=\tau _{\xi ,n}(s)$ for every $s\in \Sigma (Y).$
\end{itemize}
\end{remark}

\begin{definition}
For $\delta >0$, $s=(y_i)_{i=1}^\infty \in \Sigma (X)$ and $\mathcal{F}%
=(F_i)_{i=1}^\infty \in \Sigma (\mathbb{N})$ we define:\newline
$\frak{F}_{\delta _{}}(s,\mathcal{F})=\{G\in [\mathbb{N}]^{<\infty }:\exists
x_G^{*}\in S_{X^{*}}\;$such that $||x_G^{*}|_{Y_{F_l}}||>\delta \;\;\forall
l\in G\},$\newline
where $Y_{F_l}=<(y_i)_{i\in F_l}>$ .
\end{definition}

\begin{remark}
If $X$ is a reflexive Banach space with an unconditional basis then it is
easy to check that the set $\frak{F}_{\delta _{}}(s,\mathcal{F})$ is an
adequate and compact family of finite subsets of $\mathbb{N}$.
\end{remark}

In the sequel we will denote by $X$ a reflexive Banach space with an
unconditional basis.

\begin{lemma}
\label{Lem4.03}Suppose that for some $\xi <\omega _1$, $s=(y_i)_{i=1}^\infty
\in \Sigma (X)$, $\mathcal{F}=(F_i)_{i=1}^\infty \in \Sigma (\mathbb{N})$,
there exists $M\in [\mathbb{N}]$ and $\delta >0$ such that $||\xi _1^N(s,%
\mathcal{F})||>\delta $ for all $N\in [M].$Then there exists $L\in [M]$ such
that $S_{\xi _{}}^L\subset \frak{F}_{\delta /2}(s,\mathcal{F}),$i.e. for
every $G\in S_{\xi _{}}^L$ there exists $x_G^{*}\in S_{X^{*}}$ satisfying $%
||x_G^{*}|_{Y_{F_l}}||>\frac \delta 2$ for every $l\in G$ ,where, $%
Y_{F_l}=<(y_i)_{i\in F_l}>.$
\end{lemma}

\proof%
Observe as we have mentioned above ,that from the stability property (4) of
summability methods we get that for all $n\in \mathbb{N}$ ,$N\in [M],\left|
\left| \xi _n^N(s,\mathcal{F})\right| \right| >\delta .$Further if $t\in
\Sigma (s,\mathcal{F})$ such that $\left| \left| \xi _n^N\cdot t\right|
\right| >\delta $ and $x^{*}\in S_{X^{*}}$ such that $x^{*}(\xi _n^N\cdot
t)>\delta ${\bf \ }we get that \newline
\begin{equation*}
\sum \left\{ \xi _n^N(l):||x^{*}|_{Y_{F_l}}||>\dfrac \delta 2\right\} >\frac
\delta 2
\end{equation*}
Therefore the compact and adequate family $\frak{F}_{\delta /2}\mathbf{(s,}%
\mathcal{F}\mathbf{)\;}is\;(\xi ,M,\delta )-$large.\newline
The result follows from Theorem \ref{Th4.01} and the proof is complete 
\endproof%

\begin{lemma}
\label{Lem4.05}Suppose that for $s\in \Sigma (X),\mathcal{F}\in \Sigma (X)$
and $\delta >0$ there exists $M\in [\mathbb{N}\mathbf{]}$such that $S_{\xi
_{}}^M$ is a subset of $\frak{F}_{\delta /2}(s,\mathcal{F})$.Let $m_0=\min M.
$Then for every $\varepsilon >0,$ there exist $x^{*}\in
B_{X^{*}},suppx^{*}\subset \bigcup_{l\in F_{m_0}}suppy_l$ and $L\in
[M\backslash \left\{ m_0\right\} ]$ satisfying the following property :%
\newline
For every $G\in S_{\xi _{}}^M,G\subset L\bigcup \left\{ m_0\right\} $ and $%
\min G=m_0$ there exists $x_G^{*}\in B_{X^{*}}$ , $||x_G^{*}|_{Y_{F_l}}||>%
\dfrac \delta 2-\varepsilon $ for all $l\in G$ and $%
x_G^{*}|_{Y_{F_{m_0}}}=x^{*}$.
\end{lemma}

\proof%
Set $Y_{F_{m_0}}^{*}=span[e_i^{*}:i\in \bigcup\limits_{l\in F_{m_0}}suppy_l]$
which is a finite dimensional subspace of $X^{*}.$ Choose $%
\{x_1^{*},...,x_d^{*}\}$ be a $\frac \varepsilon 2$-net in $%
B_{Y_{F_{m_0}}^{*}}$ and for $1\leq i\leq d$ define the set\newline
$A_i=\left\{ 
\begin{array}{c}
N\in M\backslash \{m_0\}:\text{if }N^{\prime }=N\cup \{m_0\}\text{ and }%
G=supp\xi _1^{N^{\prime }} \\ 
\text{ then }\left\| x_G^{*}\left| _{Y_{F_{n_{m_0}}}}\right.
-x_i^{*}\right\| <\frac \varepsilon 2
\end{array}
\right\} $\newline
Where $x_G^{*}$ denotes the functional in $B_{X^{*}}$ witnessing\ the fact $%
G\in \frak{F}_{\delta /2}(s,\mathcal{F}).$\newline
Clearly $[M\backslash \{m_0\}]=\underset{i=1}{\overset{d}{\bigcup }}A_i$
.Then by Theorem \ref{Ra}, there exists $L\in [M\backslash \left\{
m_o\right\} $ such that $[L]\subset A_{i_0}$ for some $i_o\in \left\{
1,...d\right\} .$The set $L$ is the desired .Indeed, it follows from the
definition of $A_{i_0}$ that for every $G=supp\xi _1^N,minN=m_0,N\backslash
\left\{ m_0\right\} \in [L]$ we have the desired property. In the general
case we simply observe that $G$ a subset of some $G^{\prime }=supp\xi _1^N$
for an appropriate set $N$ with $\min N=m_0$ and the proof is complete.%
\endproof%

Before pass to the next result we make a brief introduction to the $c_0$%
-tree.\newline
A tree $\mathcal{T}$ is said to be {\bf well-founded} if every linearly
ordered subset of $\mathcal{T}$ is finite .For a well founded tree we define
as the {\bf derivative tree} $\mathcal{T}^{\prime }$ the subtree of $%
\mathcal{T}\,$consisting of all no maximal elements of $\mathcal{T}$;
inductively we define the $\mathbf{\xi }$-{\bf derivative }as $\mathcal{T}%
^{(\zeta +1)}=(\mathcal{T}^{(\zeta )})^{\prime }$ and $\mathcal{T}^{(\xi
_{})}=\bigcap_{\zeta _{}<\xi _{}}\mathcal{T}^{(\zeta )}$ if $\xi $ is a
limit ordinal. The {\bf order } $\mathbf{o}(\mathcal{T})$ of a well
founded tree $\mathcal{T}$ is the smaller ordinal $\xi $ such that the $\xi
- $derivative of $\mathcal{T}\,$is an empty set.

$c_0${\bf -trees:}Let $X$ be a separable Banach space. We define ,for $%
\delta >0,$ 
\begin{equation*}
\mathcal{T}(c_0,\delta )=\left\{ 
\begin{array}{c}
\left( x_1,...,x_n\right) :\left| \left| x_i\right| \right| =1\text{ and }%
\left( x_1,...,x_n\right) \text{is } \\ 
\frac 1\delta -\text{equivalent to the usual basis of }\ell _n^\infty
\end{array}
\right\}
\end{equation*}
Then $\mathcal{T}(c_0,\delta )$ with the natural order is a tree and further
if $c_0$ is not isomorphic to a subspace of $X$ then $\mathcal{T}(c_0,\delta
)$ is a well founded tree.\newline
The $c_0\mathbf{-}${\bf index }of the space $X$ is defined to be 
\begin{equation*}
\mathbf{o}(X)=\sup \left\{ \mathbf{o}(\mathcal{T}(c_0,\delta )):\delta
>0\right\}
\end{equation*}
\newline
The following result is due to J.Bourgain\cite{Bo}:

\begin{theorem}[J.Bourgain]
\label{Th4.03}Let $X$ be a separable reflexive Banach space Then the $c_0-$%
index $\mathbf{o}(X)$is a countable ordinal$.$%
\endproof%
\end{theorem}

\begin{proposition}
\label{Prop4.02}Let $X$ be a reflexive Banach space with an unconditional
basis and $s\in \Sigma (X),\xi <\omega _1,$ such that $\tau _{\xi _{}}(s)>0$
then the $c_0$-index of $X^{*}$is greater or equal to $\xi .$
\end{proposition}

\proof%
Since $\tau _{\xi _{}}>0$ there exists $\mathcal{F}\in \Sigma (\mathbb{N})$ and 
$M\in [\mathbb{N}\mathbf{]}$ such that for every $N\in [M]$ 
\begin{equation*}
\left| \left| \xi _1^N(s,\mathcal{F})\right| \right| >\frac{\tau _{\xi
_{}}(s)}2=\delta
\end{equation*}
\newline
Therefore there exists $L\in [M]$ such that $S_{\xi _{}}^L\subset \frak{F}%
_{\delta /2}(s,\mathcal{F})$(Lemma \ref{Lem4.03}).We consider for $G=\left\{
l_1,...l_d\right\} \in \frak{F}_{\delta /2}(s,\mathcal{F})$ the set $\left\{
x_G^{*}(1),...x_G^{*}(d)\right\} $ defined by $%
x_G^{*}(l)=x_G^{*}|_{Y_{F_l}}. $Clearly $\left\{
x_G^{*}(1),...x_G^{*}(d)\right\} $ is equivalent to $\ell _d^\infty -$basis
with constant $\frac 2\delta $.Notice that if $\varepsilon <\frac \delta 4$
and $\left\{ x_1^{*},...x_d^{*}\right\} $ is such that $%
||x_G^{*}(l)-x_l^{*}||<\varepsilon $ then $\left\{
x_1^{*},...x_d^{*}\right\} $ is equivalent to $\ell _d^\infty -$basis with
constant $\frac 2{\delta ^2}.$\newline
We will proved the desired statement by induction\newline
{\bf The inductive hypothesis:}\textit{If }$\zeta <\omega _1$\textit{\
and }$M\in [\mathbb{N}]$\textit{are such that }$S_{\zeta _{}}^M\subset \frak{F}%
_{\delta /2}(s,\mathcal{F})$\textit{\ then for every }$N\in [M],\varepsilon
>0$\textit{\ there exists a tree }\newline
$\mathcal{T}_{\zeta _{}}\subset \left\{ \left( x_1^{*},...x_d^{*}\right)
:d\in \mathbb{N}\text{ and }x_i^{*}\subset B_{X^{*}}\right\} $\textit{\ with }$%
\mathbf{o}(\mathcal{T}_{\zeta _{}})>\zeta $ \textit{satisfying the following
properties:\newline
(i) : For every }$\left( x_1^{*},...x_d^{*}\right) \in \mathcal{T}_{\zeta
_{}},suppx_1^{*}<...<suppx_d^{*}$\textit{\newline
(ii): For every }$\left( x_1^{*},...x_d^{*}\right) \in \mathcal{T}_{\zeta
_{}}$\textit{\ there exists }$G=\left\{ m_1,...m_d\right\} \in S_{\zeta
_{}}^M[N]$\textit{\ with }$||x_G^{*}(i)-x_i^{*}||<\varepsilon $\textit{\ for
all }$i=1,...,d.$\newline
{\bf Proof of the inductive hypothesis:}We will prove it for $N=M$ and
the general case is similar. Let $M=\left\{ m_1<m_2<....\right\} .$\newline
1.For $\zeta =1$ is obvious.{\bf \ }\newline
2.If $\zeta $ is a limit ordinal we consider the sequence $(\zeta
_n)_{n=1}^\infty $ which defines the family $S_{\zeta _{}}.$Then for every $%
n\in \mathbb{N}$ $\;\;S_{\zeta _n}[\mathbb{N}\backslash \left\{ 1,...n\right\}
]\subset S_{\zeta _{}}$ therefore $S_{\zeta _n}^M[\mathbb{N}\backslash \left\{
1,...n\right\} ]\subset S_{\zeta _{}}\subset \frak{F}_{\delta /2}(s,\mathcal{%
F})$ and by the inductive hypothesis there exists a tree$\mathcal{T}_{\zeta
_n}$ satisfying the inductive hypothesis for the ordinal $\zeta _n$ ,the
number $\varepsilon $ and the set $N=M\backslash \left\{ m_1,...m_n\right\} $%
.\newline
If $\mathcal{T}_{\zeta _{}}=\bigcup_n\mathcal{T}_{\zeta _n}$ it is easy to
see that $\mathcal{T}_{\zeta _{}}$ satisfies the requirement properties.%
\newline
3.If $\zeta =\frak{\eta }+1$ assume that $m_0=\min M>1$ and set $\varepsilon
^{\prime }=\varepsilon /4.$\newline
Then there exists $x^{*}\in B_{X^{*}}$ and $N\in \left[ M\backslash \left\{
m_0\right\} \right] $ satisfying the conditions of Lemma \ref{Lem4.05}.%
\newline
by the inductive hypothesis there exists a tree $\mathcal{T}_{\eta _{}}$
satisfying the properties (i) and (ii) for the set $N$ and the number $%
\varepsilon ^{\prime }.$\newline
We define 
\begin{equation*}
\mathcal{T}_{\zeta _{}}=\left\{ \left( x^{*},x_1^{*},...x_d^{*}\right)
:\left( x_1^{*},...x_d^{*}\right) \in \mathcal{T}_{\eta _{}}\right\}
\end{equation*}
It is easy to check that $\mathcal{T}_{\zeta _{}}$ is the requirement tree
and the proof of the inductive hypothesis is complete.\newline
To finish the proof we simply observe that for every $\zeta <\xi $ there
exists a natural number $n$ such that $S_{\zeta _{}}[\mathbb{N}\mathbf{%
\backslash }\left\{ 1,...n\right\} ]\subset S_{\xi _{}}.$(This can be proved
by an easy inductive argument.)Therefore if $S_{\zeta _{}}^M\subset \frak{F}%
_{\delta /2}(s,\mathcal{F})$ then for every $\eta <\zeta $ there exists a $%
k_{\eta _{}}\in \mathbb{N}$ such that $S_{\eta _{}}^M[\mathbb{N}\mathbf{\backslash 
}\left\{ 1,...k_{\eta _{}}\right\} ]\subset S_{\zeta _{}}^M$ and by the
inductive hypothesis we get the desired result.%
\endproof%

\begin{proposition}
\label{Prop4.03}Every reflexive Banach space $X$ with an unconditional basis
contains a subspace $Y$ with the following property:\newline
There exists a unique $\xi _Y<\omega _1$ such that :\newline
(1) $\tau _{\xi _Y}(Y)=0$ and $\tau _{\xi _{}}(Y)\neq 0$ for every $\xi <\xi
_Y$. \newline
(2) For every $\xi <\omega _1$, every $s\in \Sigma (Y)$ and $n\in \mathbb{N}$%
\newline
$\;\;\;\;\tau _{\xi _{},n}(s)=\tau _{\xi _{},n}(Y)$
\end{proposition}

\proof%
It follows from Lemma \ref{Lem4.07} that for every $\xi <\omega _1\,,Z\prec
X,$the subspace $Z$ contains a block subspace $Y_{\xi _{}}$ such that for
every $s\in \Sigma (Y_{\xi _{}}),\tau _{\xi _{},1}(s)=\tau _{\xi
_{}}(s)=\tau _{\xi _{}}(Y_{\xi _{}}).$ By a similar argument we can show
that for every $n\in \mathbb{N}$ and $Z\prec X$ there exists $Y_n\prec Z$ such
that $\tau _{\xi _{},n}(s)=\tau _{\xi _{},n}(Y)$ for every $s\in \Sigma
(Y_n) $ and so we can find a sequence $Y_1\succ Y_2\succ ...$ of block
subspaces such that $\tau _{\xi _{},n}(s)=\tau _{\xi _{},n}(Y_n)$ for every $%
s\in \Sigma (Y_n).$A diagonal block space $Y$ of the sequence $%
(Y_n)_{n=1}^\infty $ satisfies $\tau _{\xi _{},n}(s)=\tau _{\xi _{},n}(Y)$
for every $s\in \Sigma (Y)$ and $n\in \mathbb{N}$ .Therefore for every $\xi
<\omega _1$ and $Z\prec X$ there exists $Y_{\xi _{}}\prec Z$ such that $\tau
_{\xi _{},n}(s)=\tau _{\xi _{},n}(Y_{\xi _{}})$ for every $s\in \Sigma
(Y_{\xi _{}})$ and $n\in \mathbb{N}\mathbf{.}$ \newline
From Theorem \ref{Th4.03} and Proposition \ref{Prop4.02} we get that there
exists a countable ordinal $\xi _0$ such that $\tau _{\xi _{}}(s)=0$ for all 
$s\in \Sigma (X)$and $\xi \geq \xi _0.$\newline
Enumerate the set $\left\{ \xi :\xi \leq \xi _0\right\} $ as $\left( \xi
_n\right) _{n\in \mathbb{N}}$ and inductively choose block subspaces $Y_1\succ
Y_2\succ ...$ of $X$ such that $\tau _{\xi _i,n}(s)=\tau _{\xi _i,n}(Y_i)$
for every $s\in \Sigma (Y_i)$ and $n\in \mathbb{N}\mathbf{.}$A diagonal block
space $Y$ of the sequence $(Y_i)_{i=1}^\infty $ is the desired subspace of $%
X $.%
\endproof%

\begin{definition}
A Banach space $Y$ which satisfies the conditions of the previous
proposition is called a {\bf stabilized }space.\newline
We denote by $\xi _Y=\min \left\{ \xi :\tau _{\xi _{}}(Y)=0\right\} $
\end{definition}

\begin{lemma}
\label{Lem4.13}Let $s\in \Sigma (X)$,$\mathcal{F}\in \Sigma (\mathbb{N}),M\in [%
\mathbb{N}]$ and $\delta >0$ .Assume that for $m>n\left( \left[ \frac 2\delta
\right] +1\right) $ and $\zeta <\omega _1$we have that for all $N\in [M]$ $%
\left\| \frac{\zeta _1^N+...+\zeta _m^N}m(s,\mathcal{F})\right\| >\delta $.
Then there exists $L\in [M]$ such that $[S_{\zeta _{}}^L]^n\subset \frak{F}%
_{\delta /4}(s,\mathcal{F})$.
\end{lemma}

\proof%
We observe that if $x^{*}\in B_{X^{*}}$ , $t\in \Sigma (s,\mathcal{F})$
satisfy $x^{*}(\frac{\zeta _1^N+...+\zeta _m^N}m\cdot t)>\delta $ then there
exists $\{k_1<k_2<...<k_n\}\subset \{1,2,...,m\}$ such that $x^{*}(\zeta
_{k_i}^N\cdot t)>\dfrac \delta 2$ for all $i=1,...n.$\newline
Hence we can define a partition of \newline
\begin{equation*}
\lbrack M]=\bigcup \left\{ A_{(k_1,...k_n)}:1\leq k_1<...<k_n\leq m\right\}
\end{equation*}
,where $N\in A_{(k_1,...k_n)}$ if and only if there exists $x^{*}\in
B_{X^{*}}$ with $x^{*}(\zeta _{k_i}^N\cdot t)>\dfrac \delta 2$ for all $%
i=1,...,n.$\newline
From the dichotomy principle $\left( \text{Proposition\ref{Prop4.01}}\right) 
$there exists a $L\in [M]$ and $1\leq k_1<...<k_n\leq m$ such that $%
[L]\subset A_{(k_1,...k_n)}$ and from the stability property (4) of repeated
averages hierarchy ,for every $N\in $ $[L]$ there exists $N^{\prime }\subset
N$ such that $N^{\prime }\subset A_{(1,...,n)}$ .From this we can assume
that actually $[L]\subset A_{(1,...,n)}.$\newline
Let $N\in [L],$choose $x^{*}\in B_{X^{*}}$ such that $x^{*}(\zeta
_{k_i}^N\cdot t)>\dfrac \delta 2$ .Then for every $i=1,...,n$ we get that $%
\sum \left\{ \zeta _i^N(d):\left\| x^{*}\left| _{Y_{F_d}}\right. \right\|
>\dfrac \delta 4\right\} >\dfrac \delta 4.$Therefore $\frak{F}_{\delta /4}(s,%
\mathcal{F})$ is $(n,\zeta ,L,\dfrac \delta 4)-$large and the result follows
from Theorem \ref{Th4.01}

\begin{lemma}
Assume that for $\xi <\omega _1,s\in \Sigma (X),\;\mathcal{F}=\{F_n\}_{n\in 
\mathbb{N}}\in \Sigma (\mathbb{N})$ , $\delta >0,k\notin \mathbb{N}$ we have that $%
\tau _{\xi _{},m}(s,\mathcal{F})\geq \delta $ for some $m\geq k(\frac
2\delta +1).$Then there exists $\mathcal{G}\in \Sigma (\mathbb{N}\mathbf{)},%
\mathcal{G}\subset \mathcal{F}$ such that for every $\mathcal{G}^{\prime
}\subset \mathcal{G},M\in [\mathbb{N}]$ we have that 
\begin{equation*}
\left| \left| \frac{\zeta _1^M+...+\zeta _k^M}k(s,\mathcal{G}^{\prime
})\right| \right| >\frac \delta 4
\end{equation*}
\end{lemma}

\proof%
It follows from Lemma \ref{Lem4.13} that there exists $L\in [\mathbb{N}]$ such
that $[S_{\xi _{}}^L]^k\subset \frak{F}_{\delta /4}(s,\mathcal{F}).$We set $%
\mathcal{G}=\left\{ G_l\right\} _{l\in \mathbb{N}}=\left\{ F_{n_l}\right\}
_{l\in L}.$ Observe that for $\mathcal{G}^{\prime }\subset \mathcal{G},%
\mathcal{G}^{\prime }=\left\{ G_l^{\prime }\right\} _{l\in \mathbb{N}},$ we
have $[S_{\xi _{}}^L]^k\subset \frak{F}_{\delta /4}(s,\mathcal{G}^{\prime
}). $Therefore for every $M\in [N],\left| \left| \frac{\zeta _1^M+...+\zeta
_k^M}k(s,\mathcal{G}^{\prime })\right| \right| >\frac \delta 4$ and the
proof is complete.%
\endproof%

\begin{remark}
The content of the above Lemma is that at the moment that we know that \\
$\left| \left| \frac{\zeta _1^M+...+\zeta _k^M}k(s,\mathcal{G}^{\prime
})\right| \right| >\delta $ for some $\delta >0,$ the pair $(s,\mathcal{F})$
and all $N\in [L]$ then by going to a certain $\mathcal{G}\subset \mathcal{F}
$ we know the conclusion for $\frac \delta 4$ and all $L$ subsets of $\mathbb{N}
$.
\end{remark}

\begin{proposition}
\label{4.15}Let $X$ be a stabilized reflexive space with an unconditional
basis .Then the following hold:\newline
(a) If $\xi =\xi _X$ is a limit ordinal and $(\xi _n)_{n\in \mathbb{N}}$ the
increasing sequence that defines $S_{\xi _{}}$ then \underline{$lim$}$\tau
_{\xi _n}=0$\newline
(b) If $\xi =\xi _X$ is of the form $\xi =\zeta +1$ then \underline{$lim$}$%
\tau _{\zeta _{},n}=0$.
\end{proposition}

\proof%
We prove (a). The proof of (b) is similar. Assume that (a) fails. Then if $%
s=\left( e_k\right) _{k\in \mathbb{N}}$,the basis of $X$,since $X$ is
stabilized there exists a $\delta >0$ such that $\tau _{\xi _n}(s)>\delta $
for all $\delta >0$. For $\mathcal{F}_1=(F_l^1)$, $\mathcal{F}_2=(F_l^2)$ ,
we denote by $\mathcal{F}_2\prec \mathcal{F}_1$ if $\mathcal{F}_2$ is a
block subfamily of $\mathcal{F}_1$ i.e. every $F_l^2=\underset{d\in G_l}{%
\bigcup }F_d^1$ for $G_1<G_2<....$ subsets of $\mathbb{N}.$\newline
We inductively construct $\mathcal{F}_1\succ \mathcal{F}_2\succ ...\succ 
\mathcal{F}_n\succ ...$ such that for all $n\in \mathbb{N}$, $\mathcal{G}%
_n\subset \mathcal{F}_n$ , $M\in [\mathbb{N}]$, $l\in \mathbb{N}$. $||(\xi
_k)_1^M(s,\mathcal{G}_k)||>\frac \delta 4$.\newline
We show how we produce $\mathcal{F}_1,\mathcal{F}_2$ and in the same manner
we get the general inductive step. Since $\tau _{\xi _1}(s)>\delta $ there
exists $\mathcal{G}_1$ and $M\in [\mathbb{N}]$ such that for all $N\in [M]$, $%
n\in \mathbb{N}\mathbf{,}\;\;$ $||(\xi _n)_l^N\cdot (s,\mathcal{G}_1)||>\delta $%
.\newline
From the previous lemma there exists $\mathcal{F}_1\subset \mathcal{G}_1$
satisfying the desired property. Set $s_1=(e_{l_n})_{n\in \mathbb{N}}$
subsequence of $(e_n)_{n\in \mathbb{N}}$ such that $\mathcal{F}_1=(F_{l)_{l\in 
\mathbb{N}}}$ and $n_l\in F_l$. Since $\tau _{\xi _{n_2}}(s_1)>\delta $ there
exists $\mathcal{G}_2$ and $M\in [\mathbb{N}]$ such that for all $N\in [M]$, $%
l\in \mathbb{N}$ $||(\xi _{n_2})_l^N(s_1,\mathcal{G}_2)||>\delta $.\newline
Observe that from the choice of $s_2$ the family $\mathcal{G}_2$ can be
taken such that $\mathcal{G}_2\prec \mathcal{F}_1$. Choose, as before ,$%
\mathcal{F}_2\subset \mathcal{G}_2$ satisfying the conclusion and the
inductive definition of $(\mathcal{F}_k)_{k\in \mathbb{N}}$ is complete.\newline
Each $\mathcal{F}_k=(F_n^k)_{n\in \mathbb{N}}$ and we set $\mathcal{F}_\infty
=(F_k^k)_{k\in \mathbb{N}}$.\newline
{\bf Claim}$:${\bf \ \ \ \ \ }For\ every\ \ $M\in \left[ \mathbb{N}%
\right] ,\;\;\;\;||(\xi _{})_1^M\cdot (s,\mathcal{F}_\infty )||>\frac \delta
4.$\newline
Indeed suppose that $M=\{m_1,m_2,....\}$ then by the definition of $(\xi
_0)_1^M$ we have that $(\xi _0)_1^M=(\xi _{m_1})_1^M$. Since for all $k\geq
m_1$, $\mathcal{F}_k\prec \mathcal{F}_{m_1}$ it follows that the exists $%
(F_{l_n}^{m_1})_{k\geq m_1}$ such that $l_{m_1}=m_1$ and $%
F_{l_k}^{m_1}\subset F_k^k$ set $\mathcal{G}=\{F_1^{m_1},...,F_{m_1}^{m_1}\}%
\cup \{F_{l_k}^{m_1}\}_{k>m_1}$ then we have that for all $L\in [\mathbb{N}]$ $%
||(\xi _{n_{m_1}})_1^L\cdot (s,\mathcal{G})||>\frac \delta 8$.\newline
Observe now that if $\mathcal{G}=[G_k]_{k\in \mathbb{N}}$ then for $k\geq m_1$ $%
G_m\subset F_k^k$ hence for the given $M$ 
\begin{equation*}
||(\xi _{n_{m_1}})_1^M\cdot (s,\mathcal{G})||\leq ||(\xi
_{n_{m_1}})_1^M\cdot (s,\mathcal{F}_\infty )||
\end{equation*}
and the proof, is complete. \newline
Therefore $\tau _{\xi _{}}(s)>\frac \delta 8$ and this contradiction
finishes the proof for the limit ordinal case.\newline
The case of successor ordinal is proved by similar arguments.

\begin{remark}
With a more careful construction we can actually show that $lim\tau _{\xi
_{n_k}}=0$.
\end{remark}

\begin{proposition}
\label{Prop4.05}Let $X$ be a stabilized reflexive Banach with an
unconditional basis such that $\xi _X$ is a successive ordinal. Then $X$
satisfies the property (P2).
\end{proposition}

\proof%
Let $\xi _X=\zeta +1$, $\tau =\tau _{\zeta _{}}(X)>0$. By Proposition \ref
{4.15} \underline{$lim$}$\tau _{\zeta _{},n}(X)=0$ and so for every $k\in 
\mathbb{N}$ there exists $n_k\in \mathbb{N}$ such that $\tau _{\zeta
_{},n_k}^{_{}}(X)<\frac \tau {2k}$. \newline
We set $C_k=n_k$, for every $k\in \mathbb{N}$.We will prove that for every $%
(x_i)_{i=1}^\infty \prec X$ and $k\in \mathbb{N}$ there exists $%
t=(z_i)_{i=1}^\infty \prec (x_i)_{i=1}^\infty $ such that for every $%
M=\left\{ m_1,...m_n,...\right\} \subseteq [\mathbb{N}],\left\| \dfrac{%
z_{m_1}+...+z_{m_{n_k}}}{n_k}\right\| <\dfrac 1k,$ which implies the
property $P(2).$\newline
Let $s=(x_i)_{i=1}^\infty \prec X$ and $Y=<(x_i)_{i=1}^\infty >$. Since $%
\tau _{\zeta _{},n_k}^{_{}}(s)=\tau _{\zeta _{},n_k}^{_{}}(X)<\dfrac \tau
{2k}$ ,for every $\mathcal{F}\in \Sigma (\mathbb{N})$ and $M\in [\mathbb{N}\mathbf{%
]}$ there exists $N\in [M]$ such that if $t\in \Sigma (s,\mathcal{F})$ then 
\begin{equation}
\left| \left| \dfrac{\zeta _1^N+...+\zeta _{n_k}^N}{n_k}\cdot t\right|
\right| <\frac \tau {2k}  \label{r1}
\end{equation}
Since $\tau _{\zeta _{}}(s)=\tau $ there exists $\mathcal{F}_1\in \Sigma (%
\mathbb{N})$ and $M_1\in [\mathbb{N}]$ such that for every $N\in [M_1]$ and every $%
n\in \mathbb{N}$\ $\;||\zeta _n^N(s,\mathcal{F}_1)||>\frac \tau 2.$\newline
For every $n\in \mathbb{N}$ we choose $t_n=(y_i^{(n)})_i\in \Sigma (s,\mathcal{F%
}_1)$ such that 
\begin{equation}
||\zeta _n^{M_1}\cdot t_n||>\dfrac \tau 2  \label{r2}
\end{equation}
\newline
Let $k_n=maxsupp\zeta _n^{M_1}$ and\newline
$\tilde{s}=(\tilde{y}_i)_{i=1}^\infty
=(y_1^{(1)},...,y_{k_1}^{(1)},y_{k_1+1}^{(2)},...y_{k_2}^{(2)},y_{k_2+1}^{(3)},...,y_{k_3}^{(3)},y_{k_3+1}^{(4)}...). 
$\newline
Then for every $n\in \mathbb{N}$%
\begin{equation}
||\zeta _n^{M_1}\cdot \tilde{s}||>\dfrac \tau 2  \label{r3}
\end{equation}
We set 
\begin{equation}
z_n=\frac{\zeta _n^{M_1}\cdot \tilde{s}}{\left\| \zeta _n^{M_1}\tilde{s}%
\right\| }\;\;\;\;\;,\;\;\;\;t=(z_n)_{n=1}^\infty  \label{r5}
\end{equation}
\newline
Let $M\in [M_1]$, $M=\{m_1,m_2,...\}$and set $M_2=\overset{\infty }{%
\bigcup\limits_{k=1}}$\textit{supp}$\zeta _{m_k}^{M_1}$. \newline
By the property (4) of the repeated averages hierarchy 
\begin{equation}
||\zeta _{m_k}^{M_1}\cdot \tilde{s}||=||\zeta _k^{M_2}\cdot \tilde{s}||
\label{r4}
\end{equation}
\newline
By(\ref{r5}) ( \ref{r3}) ( \ref{r1})and (\ref{r4}) and the unconditionality
of the basis we get \newline
\begin{eqnarray*}
\left\| \dfrac{z_{m_1}+...+z_{m_{n_k}}}{n_k}\right\| &\leq &\dfrac 2\tau
\left\| \dfrac{\zeta _{m_1}^{M_1}\cdot \tilde{s}+...+\zeta
_{m_{n_k}}^{M_1}\cdot \tilde{s}}{n_k}\right\| = \\
&=&\dfrac 2\tau \left\| \dfrac{\zeta _{_1}^{M_2}\cdot \tilde{s}+...+\zeta
_{_{n_k}}^{M_2}\cdot \tilde{s}}{n_k}\right\| <\dfrac 1k
\end{eqnarray*}
\newline
and the proof is completed%
\endproof%

Given $A$ be a reflexive Banach space with an unconditional basis and we
denote by $X_A$ the space constructed in the previous section 3. We denote
\\ $\;\Sigma (X_A)=\left\{ (x_i)_{i=1}^\infty :\left\| x_i\right\| =1\mathbf{,}
\text{ and }range(x_i)<range\left( x_{i+1}\right) \right\} .$\newline
The indexes $\tau _{\xi _{}}(s),\tau _{\xi _{},n}(s)$ for $s\in \Sigma
(X_A)\,$ are as in Definitions \ref{Def4.01} and \ref{4.4}.We also define $%
\tau _{\xi _{},n}(X_A)=\inf \{\tau _{\xi _{},n}(s):s\in \Sigma (X_A)\}.$

\begin{lemma}
\label{4.17}If the space $A$ is stabilized then $X_A$ is stabilized and for $%
s\in \Sigma (A),s^{\prime }\in \Sigma (X_A)$ we have $\tau _{\xi
_{},n}(s)=\tau _{\xi _{},n}(s^{\prime })$ for every $\xi <\omega _1,n\in 
\mathbb{N}\mathbf{.}$
\end{lemma}

\proof%
It is similar to that of Proposition \ref{PP} that if the space $A$
satisfies the property (P) the same fact also holds for the space $X_A$%
\endproof%

\begin{proposition}
\label{Prop4.04}Let $A$ be a stabilized reflexive space such that $\xi _A$
is a limit ordinal. Then the set $W$ in $X_A$ is an $\mathbf{a}$-thin set
for a null sequence $\mathbf{a}=(a_n)_{n\in \mathbb{N}}$ of positive numbers.
\end{proposition}

\proof%
It follows from Proposition \ref{4.15} that for the sequence $(\xi _n)_{n\in 
\mathbb{N}}$ that defines $S_{\xi _A}$ we have \underline{$\lim $}$\tau _{\xi
_n}=0$. \newline
For every $k\in \mathbb{N}$ we choose $\xi _{n_k}$ such that $\tau _{\xi
n_k}<\dfrac 1{k2^k}$ and we define $a_k=\dfrac{\tau _{\xi n_k}}{32k}.$%
\newline
We will show that $W$ is an $(a_k)_{k\in \mathbb{N}}$ thin set.\newline
Assume that it is not. Then there exists $Z$ infinite dimensional subspace
of $X_A$ and $\lambda >0$ such that for all $n\in \mathbb{N}$%
\begin{equation*}
B_Z\subset \lambda (2^kW+a_kB_{X_A})\text{.}
\end{equation*}
\newline
Choose $k>1000\lambda $ and in the sequel we will denote by $\zeta $ the
ordinal $\xi _{n_k}$ and by $\tau $ the number $\tau _{\xi _{n_k}}$.\newline
We also may assume that $Z$ is a block subspace generated by a sequence $%
(z_l)_{l\in \mathbb{N}}$ such that $rangez_l<rangez_{l+1}$,where \newline
$range\left( z\right) =\left\{ \delta \in \frak{D}:\exists \alpha ,\beta \in 
\text{\textit{supp}}z:|\alpha |\leq |\delta |\leq |\beta |\right\} $.\newline
Since the space $A$ is stabilized from Lemma \ref{4.17} we get that for
every $s\in \Sigma (X_A)$ and $\xi <\omega _1$,$\tau _{\xi _{}}(s)=\tau
_{\xi _{}}(A).$ Therefore for the sequence $s=(z_l)_{l\in \mathbb{N}}$ there
exists an $\mathcal{F}=(F_m)_m$, satisfying the property:\newline
There exists $M\in [\mathbb{N}]$ such that $\forall N\in [M]$, 
\begin{equation}
\frac 12\tau <||\zeta _1^N(s,\mathcal{F})||<2\tau \text{.}  \label{rr0}
\end{equation}
\newline
Further we may assume that $\mathcal{F}$ satisfies the following stronger
property:\newline
\begin{eqnarray}
\text{For every }N &\in &[\mathbb{N}\mathbf{]}\text{ there exists }x^{*}\in
B_{X_A^{*}}\text{ with}  \label{*} \\
\newline
||x^{*}|Z_{F_m}|| &>&\tau /4\,\,\,\,\text{ }\forall m\in \text{\textit{sup}}%
p\zeta _1^N.  \notag
\end{eqnarray}
\newline
The property \ref{*} follows from Lemma \ref{Lem4.03}\newline
For every $m\in \mathbb{N}$ we define $Z_{F_m}=span[\left( z_i\right) _{i\in
F_m}]$ and \newline
$X_{F_m}=span[\left\{ e_{\delta _{}}:\delta \in \bigcup_{i\in
F_m}rangez_i\right\} ]$.\newline
Consider the following family of finite subsets of $\mathbb{N}$\newline
\begin{equation*}
\frak{F}=\left\{ 
\begin{array}{c}
\;\;G\subset \mathbb{N}:\exists x_G^{*}\in B_{X_A^{*}}\text{, and }\left(
y_m\right) _{m\in G}\text{ such that} \\ 
\text{(a) \ }y_m\in Z_{F_m}\text{ and }||y_m||\leq 1\text{ \ \ \ \ \ \ \ \ \
\ \ \ \ \ \ \ \ \ \ } \\ 
\text{(b) }\forall m\exists E_m^{\prime }\subset
rangey_m:x_G^{*}(E_m^{\prime }y_m)\geq \dfrac{\tau \,}{16}\,\text{and } \\ 
\;\;\text{the sets }\left( E_m^{\prime }\right) _{m\in G}\text{ are pairwise
incomparable} \\ 
\text{(c)\ \ \ \ \ }||x_G^{*}|_{X_{F_m}}||<8\tau \text{ \ \ \ \ \ \ \ \ \ \
\ \ \ \ \ \ \ \ \ \ \ \ \ \ \ \ \ }
\end{array}
\right\}
\end{equation*}
\newline
{\bf Claim: }There exists $\varepsilon >0$\ such that $\frak{F}$\ is $%
(\zeta ,M,\varepsilon )$\ large.\newline
Suppose that the claim has been proved. Then we derive to a contradiction in
the following manner.\newline
By Theorem \ref{Th4.01} and Lemma \ref{Lem4.01} (see also Corollary\ref
{Cor4.01} )there exists a $L\in [M]$ such that $N\in L$the set $G^N=$\textit{%
sup}$p\zeta _1^N\backslash \{\min $\textit{sup}$p\zeta _1^N\}$ belongs to $%
\frak{F}.$\newline
From Property (1) of the repeated averages hierarchy we may consider $N\in
[L]$such that 
\begin{equation}
\underset{m\in G^N}{\sum }\zeta _1^N(m)>\frac 12.  \label{rr1}
\end{equation}
\newline
Also consider a functional $x^{*}=x_{G^N}^{*}\in B_{X_A^{*}}$ and a sequence 
$\left( y_m\right) _{m\in G^N}$ to satisfy conditions (a),(b),(c) above in
order $G^N$ belongs to $\frak{F}.$\newline
We set 
\begin{equation*}
\tilde{y}=\underset{m\in G^N}{\sum }\zeta _1^N(m)y_m,y=\underset{m\in G^N}{%
\sum }\zeta _1^N(m)E_m^{\prime }y_m
\end{equation*}
Choose $\tilde{w}\in W_0$ such that $\left| \left| \;\;\;||\tilde{y}||^{-1}%
\tilde{y}-\lambda \cdot 2^k\tilde{w}\right| \right| <\lambda a_k$ and set%
\newline
$w_m=E_m^{\prime }\tilde{w},w=\sum_{m\in G^N}w_m,u=||\tilde{y}%
||^{-1}y-\lambda \cdot 2^kw.$Then 
\begin{equation}
\left\| u\right\| =\left\| ||\tilde{y}||^{-1}y-\lambda \cdot 2^kw\right\|
<\lambda a_k  \label{rr2}
\end{equation}
\newline
By (\ref{rr0})\thinspace we get that 
\begin{equation}
||\tilde{y}||<2\tau  \label{rr3}
\end{equation}
\newline
and by \ref{rr1} and condition(b) of the definition of $\frak{F}$ 
\begin{equation}
x^{*}\left( y\right) >\frac \tau {32}  \label{rr4}
\end{equation}
\newline
Since the sets $\left( E_m^{\prime }\right) _{m\in G^N}$ are pairwise
incomparable we get that $w$ is a convex combination of vectors $w_m^{\prime
}\in B_{X_{F_m}},m\in G^N$.Therefore by property $\left( \text{c}\right) $of
the definition of the family $\frak{F}$ we have that 
\begin{equation}
x^{*}\left( w\right) \leq 8\tau  \label{rr5}
\end{equation}
\newline
Combining the above relations (\ref{rr3}),(\ref{rr4}) and (\ref{rr5}) we get
that 
\begin{equation*}
x^{*}\left( u\right) >\frac 1{2\tau }\frac \tau {32}-2^k\lambda 8\tau
>\lambda a_k
\end{equation*}
\newline
which contradicts to (\ref{rr2}) \newline
Hence it remains to prove the claim.\newline
Given any $N\in [M]$ set $G=$\textit{sup}$p\zeta _n^N$ and choose $x^{*}\in
B_{X_A^{*}}$ such that $||x^{*}|_{Z_{F_m}}||>\dfrac \tau 4$ for all $m\in G$.%
\newline
We write $x^{*}=\left( x^{*}\right) ^{+}-\left( x^{*}\right) ^{-}$ where $%
\left( x^{*}\right) ^{+}\left( e_{\delta _{}}\right) =\left\{ 
\begin{array}{c}
x^{*}(e_{\delta _{}})\text{ \ \ if }x^{*}(e_{\delta _{}})>0 \\ 
0\text{ \ \ \ otherwise}
\end{array}
\right. $\newline
and $\left( x^{*}\right) ^{-}\left( e_{\delta _{}}\right) =\left\{ 
\begin{array}{c}
-x^{*}(e_{\delta _{}})\text{ \ \ if }x^{*}(e_{\delta _{}})<0 \\ 
0\text{ \ \ \ otherwise}
\end{array}
\right. .$\newline
We set $G_{+}=\left\{ m\in G:||(x^{*})^{+}|_{Z_{F_m}}||>\dfrac \tau
8\right\} $ and \newline
$G_{-}=\left\{ m\in G:||(x^{*})^{-}|_{Z_{F_m}}||>\dfrac
\tau 8\right\} $.\newline
Then either $\newline
\sum_{m\in G_{+}}\zeta _1^N(m)>\frac 12$ or $\sum_{m\in G_{-}}\zeta
_1^N(m)>\frac 12.$\newline
Suppose that 
\begin{equation}
\sum_{m\in G_{+}}\zeta _1^N(m)>\frac 12  \label{rr9}
\end{equation}
and set$\newline
$%
\begin{equation*}
z^{*}=(x^{*})^{+}\;\;\;\;\;\;\;\;\;,G^1=\{m\in
G_{+}:||z^{*}|_{X_{F_m}}||<8\tau \}
\end{equation*}
Since for each $t=(t_m)_{m\in G}$, $t_m\in X_{F_m}$, $||t_m||=1$ we have $%
||\zeta _1^N\cdot t||<2\tau $ we get by \ref{rr9} that 
\begin{equation}
\newline
\sum_{m\in G^1}\zeta _1^N(m)>\frac 14\text{.}  \label{rr6}
\end{equation}
Further choose 
\begin{equation}
y_m\in Z_{F_m}:||y_m||=1,x^{*}(y_m)>\frac \tau 4.  \label{rr7}
\end{equation}
For every $m\in G^1$ choose $w_m\in W$ such that\newline
\begin{equation}
||y_m-\lambda 2^kw_m||<\lambda a_k  \label{rr8}
\end{equation}
\newline
Set $\varepsilon _1=\dfrac 1{32}\dfrac \tau {2^k\lambda },E_m=range\left(
y_m\right) .$\newline
From Proposition \ref{AC},there exists a partition $D_1,...D_N$ of the set $%
G $ and $E_m^{\prime }\subset E_m$ such that \newline
$\left( i\right) $ $z^{*}\left( \left( E_m\backslash E_m^{\prime }\right)
w_m\right) <\varepsilon _1$\newline
$\left( ii\right) $The sets $\left( E_m^{\prime }\right) _{m\in D_i}$ are
pairwise incomparable for every $i=1,...,N$\newline
$\left( iii\right) N<\dfrac 5{\varepsilon _1^2}$\newline
By (\ref{rr6}) there exists $i_0\in \left\{ 1,...,N\right\} $such that 
\newline
\begin{equation*}
\sum_{m\in D_{i_0}}\zeta _1^N(m)>\frac 1{4N}>\frac{\varepsilon _1^2}{20}%
\text{.}
\end{equation*}
\newline
By (\ref{rr7}) ,(\ref{rr8})and $\left( i\right) :$ 
\begin{equation*}
z^{*}\left( E_m^{\prime }y_m\right) >\frac \tau 8-2^k\lambda \varepsilon
_1-\lambda a_k>\frac \tau 8-\frac 1{32}\frac \tau {2^k\lambda }-\lambda
\frac \tau {32k}>\frac \tau {16}
\end{equation*}
\newline
Therefore $D_{i_0}\in \frak{F}$ \newline
Hence for $\varepsilon =\dfrac{\varepsilon _1^2}{20}$ the set $\frak{F}$ is $%
(\zeta ,M,\varepsilon )$ large and the proof is complete.\newline
\endproof%

\begin{theorem}
\label{Th4.04}Let $A$ be a reflexive Banach space with an unconditional
basis .Then there exists a block subspace $B$ and a null sequence $\mathbf{a}
$ such that the set $W$ in $X_B$ is an{\bf \ }$\mathbf{a-}$thin set
\end{theorem}

\proof%
From Proposition\ref{Prop4.03} we can choose a block subspace $B$ of $A$
which is stabilized .Then if $\xi _B$ is a successor ordinal or $\xi _A=0$
the result follows from Proposition \ref{Prop4.05} and if $\xi _B$ is a
limit ordinal then it follows from Proposition \ref{Prop4.04}%
\endproof%

\section{{\bf THIN NORMING SETS III}}

\begin{center}
{\bf (}$c_0${\bf (}$\mathbb{N}${\bf ) case)}
\end{center}

In this section for $A=c_0\left( \mathbb{N}\right) $ we construct a Banach
space $X_A$ which is of the form $\left( \overset{}{\underset{n=1}{\overset{%
\infty }{\sum }}\bigoplus \ell ^1\left( 2^n\right) }\right) _0\,$\,\,\\and a set
$W $ which norms a subspace of $X_A^{*}$ isometric and $w^{*}$-homeomorphic
to $\left( c_0\left( \mathbb{N}\right) \right) ^{*}$.The set $W$ is defined
similarly to the reflexive spaces case. The proof that $W$ is a thin set
uses the structure of the $w^{*}$-closure of $W$ in the second dual of $X_A.$

\begin{theorem}
\label{5.1}There exists a Banach space $X$ with separable dual $X^{*}$ and a
closed, bounded convex symmetric subset $W$ of $X$ such that $W$ is a thin
subset of $X$, $c_0(\mathbb{N})^{*}$ is isometric and $w^{*}$-isomorphic to a
closed subspace $Y$ of $X^{*}$ and $W$ is a $\frac 12$-norming set for $Y$
\end{theorem}

\proof%
{\bf (i) The space }$\mathbf{X}$\newline
We denote by $\frak{D}$ the dyadic tree, i.e. every $\delta \in \frak{D}$
has exactly two immediate successors. The tree $\frak{D}$ naturally
coincides with the set of all finite sequences $\{(\varepsilon
_1,...,\varepsilon _n):n\in \mathbb{N}$ and $\varepsilon _i\in \{0,1\}$, for $%
i=1,2,...,n\}$.\newline
We will use the notation and the definitions about the tree $\frak{D}$
introduced in the reflexive case.\newline
In the vector space $c_{00}(\frak{D})$ we define the following norm

\begin{center}
$\left\| \underset{\delta \in \frak{D}}{\sum }\lambda _{\delta _{}}e_{\delta
_{}}\right\| =\underset{n\in \mathbf{N}}{\max }\left( \underset{|\delta |=n}{%
\sum }|\lambda _{\delta _{}}|\right) $
\end{center}

The space $X$ is the completion of $c_{00}(\frak{D})$ with the above defined
norm. It is clear that $X$ is isometric to the space $\left( \underset{n\in 
\mathbf{N}}{\sum }\oplus \ell ^1(2^n)\right) _0$. The space $X^{*}$ is
isometric to the space $\left( \underset{n\in \mathbf{N}}{\sum }\oplus \ell
^\infty (2^n)\right) _1$ and $X^{**}$ is isometric to $\left( \underset{n\in 
\mathbf{N}}{\sum }\oplus \ell ^1(2^n)\right) _\infty $. The last space is
not a $d$-product of $(\ell ^1(2^n))_{n\in \mathbb{N}}$.

$_{}$\newline
{\bf (ii) The space }$\mathbf{Y}$\newline
The space $Y$ is defined in a similar way as in the reflexive case. We set 
\newline
$y_n^{*}=\underset{\delta (n)=1}{\underset{|\delta |=n}{\sum }}e_{\delta
_{}}^{*}$, where $\delta =(\varepsilon _1,...,\varepsilon _n)$, $\delta
(n)=\varepsilon _n$, and

\begin{center}
$Y=\overline{<(y_n^{*})_{n\in \mathbb{N}}>}$.
\end{center}

It is easy to see that the space $Y$ is isometric to $\ell ^1$ and $w^{*}$%
-isomorphic to $c_0^{*}(\mathbb{N})$.

$_{}$\newline
{\bf (iii) The set }$\mathbf{W}$\newline
For any initial segment $s$ of $\frak{D}$ we denote

\begin{center}
$x_s=\underset{\delta \in s}{\sum }e_{\delta _{}}$
\end{center}

Clearly, $||x_s||=1$.\newline
We set $W=\overline{co}\{\pm x_s:s$ is an initial segment$\}$.\newline
For any infinite branch $\gamma $ of $\frak{D}$ we set

\begin{center}
$x_{\gamma _{}}=\underset{\delta \in \gamma }{\sum }e_{\delta _{}}=w^{*}-%
\underset{n\rightarrow \infty }{\lim }x_{\gamma |_n}$
\end{center}

where the series is taken in the $w^{*}$-topology and $x_{\gamma _{}}\in
X^{**}$.

\begin{proposition}
The set $W$ is $\frac 12$-norming for the space $Y$.
\end{proposition}

\begin{proof}
Let $y^{*}\in Y$ such that $||y^{*}||=1$. Then $y^{*}=\underset{n=1}{%
\overset{\infty }{\sum }}\lambda _ny_n^{*}$ and $\underset{n=1}{\overset{%
\infty }{\sum }}|\lambda _n|=1$. Set $M=\{j\in \mathbb{N}:\lambda _j>0\}$. We
can assume that $\underset{j\in M}{\sum }\lambda _j\geq \frac 12$ and let $%
\gamma =(\varepsilon _i)_{i\in \mathbb{N}}$ where $\varepsilon _i=\left\{ 
\begin{array}{c}
1\,\,\,\,\,\text{if\thinspace \thinspace \thinspace \thinspace \thinspace }%
\lambda _i>0 \\ 
0\,\,\,\,\,\text{if\thinspace \thinspace \thinspace \thinspace }\lambda
_i\leq 0
\end{array}
\right. $.\newline
Then $x_{\gamma _{}}(y^{*})\geq \frac 12$ and hence $\underset{n\rightarrow
\infty }{\lim }x_{\gamma |_n}(y^{*})\geq \frac 12$ which proves that $W$ $%
\frac 12$-norms $Y.$
\end{proof}

\begin{proposition}
The set $W$ is a thin subset of $X.$
\end{proposition}

\proof%
Assume on the contrary that $W$ is not a thin set. Then for every $%
\varepsilon >0$ there exists a normalized block sequence $(z_n)_{n\in \mathbb{N}%
}$ in $X$ and $\lambda >0$ such that

\begin{center}
$B_Z\subset \lambda W+\varepsilon B_X$
\end{center}

where $Z=\overline{<(z_n)_{n\in \mathbb{N}}>}$.\newline
We denote by $\tilde{W}$ the $w^{*}$-closure of $W$ in $X^{**}$ and by the
compactness of $\tilde{W}$ in the $w^{*}$-topology we get that

\begin{center}
$B_{Z^{**}}\subset \lambda \tilde{W}+\varepsilon B_{X^{**}}$.
\end{center}

We denote by $Br(\frak{D})$, $S(\frak{D})$ the sets of infinite branches,
initial segments of $\frak{D}$ respectively.\newline
Define

$K_c=\{x_{\gamma _{}}:\gamma \in Br(\frak{D})\}$

$K_d=\{x_s:s\in S(\frak{D})\}$

$K=K_c\cup K_d$\newline
It is easy to check that $\tilde{W}=\overline{co}^{w^{*}}(K\cup -K)$.\newline
Finally $\mathcal{M}(K_c)$ denotes the space of the regular finite Borel
measures on $K_c$, $\mathcal{M}_1(K_c)=\{\mu \in \mathcal{M}(K):||\mu ||\leq
1\}$ and $Q:X^{**}\rightarrow X^{**}/X$ the natural quotient map.

\begin{lemma}
If $R:\mathcal{M}_1(K_c)\rightarrow \tilde{W}$ is the natural affine map
from the unit ball of the Borel measures on $K_c$ onto $\tilde{W}$ then $QR$
is an isometry between $\mathcal{M}_1(K_c)$ and $Q(\tilde{W})$.
\end{lemma}

\begin{proof}
Notice, first that $(K_c,w^{*})$ is a compact metric space homeomorphic to
the Cantor set $\{0,1\}^{\mathbb{N}}$. Hence for any $\delta \in \frak{D}$ the
family $(\mathcal{N}_{\delta _{}})_{\delta \in \frak{D}}$, where $\mathcal{N}%
_{\delta _{}}=\{x_{\gamma _{}}:\delta \in \gamma \}$ is a basis for the
topology of $K_c$ and each $\mathcal{N}_{\delta _{}}$ is a clopen subset of $%
K_c$.\newline
Now for $\mu \in \mathcal{M}_1(K_c)$ $||\mu ||=\underset{n\rightarrow \infty 
}{\lim }\underset{|\delta |=n}{\sum }|\mu (\mathcal{N}_{\delta _{}})|=%
\underset{n\rightarrow \infty }{\lim }\underset{|\delta |=n}{\sum }|R\mu
(e_{\delta _{}})|=d(R\mu ,X)=||QR\mu ||$, and this completes the proof of
the lemma.
\end{proof}

\begin{lemma}
\label{5.5}Let $(\mu _\xi )_{\xi \in \Xi }$ be an uncountable family of
measures in $\mathcal{M}_1(K_c)$. Then for any $\delta >0$ there exists an
uncountable family $\left\{ \{\xi _i,\zeta _i\},i\in I\right\} $ of pairwise
disjoint two-points sets such that $\{\mu _{\xi _i}-\mu _{\zeta _i},i\in I\}$
are {\bf pairwise }$\mathbf{\delta }${\bf -singular}, i.e. for every $%
i\in I$ there exists a $v_i\in \mathcal{M}_1(K_c)$ such that $||v_i||<\delta 
$ and the elements of the family $\{\mu _{\xi _i}-\mu _{\zeta
_i}-v_i\}_{i\in I}$ are pairwise singular.
\end{lemma}

\proof%
: It is well known ([L]) that $\mathcal{M}_1(K_c)$ is isometric to\newline
$\left( \underset{\gamma <2^w}{\sum }\oplus L^1(\mu _{\gamma _{}})\right) _1$
where $\{\mu _{\gamma _{}}\}_{\gamma <2^w}$ is a maximal family of pairwise
singular probability measures in $\mathcal{M}_1(K_c)$.\newline
The isometry assigns to each $\mu \in \mathcal{M}(K_c)$ the vector $\left\{ 
\frac{d\mu }{d\mu _{\gamma _{}}}\right\} _{\gamma <2^w}$.\newline
Since $K_c$ is a compact metric space, each $L^1(\mu _{\gamma _{}})$ is a
separable Banach space.\newline
Given $\delta >0$, for each $\xi \in \Xi $ we choose a finite subset $F_{\xi
_{}}$ of $2^w$ such that

$\left\| \mu _{\xi _{}}-\underset{\gamma \in F_{\xi _{}}}{\sum }\frac{d\mu
_{\xi _{}}}{d\mu _{\gamma _{}}}\right\| <\frac \delta 4$.

Apply Erd\"{o}s-Rado theorem ([KM]) to get an uncountable $\Xi ^{\prime }$
subset of $\Xi $ and $F$ subset of $2^{\omega _{}}$ such that for $\xi
_1,\xi _2\in \Xi ^{\prime }$, $\xi _1\neq \xi _2$ we have $F_{\xi _1}\cap
F_{\xi _2}=F$. If $F=\emptyset $ then the measures $(\mu _{\xi _{}})_{\xi
\in \Xi ^{\prime }}$ are pairwise $\frac \delta 4$-singular, hence for any
family $\left\{ \{\xi _i,\zeta _i\},i\in I\right\} $ of pairwise disjoint
two-point subset of $\Xi ^{\prime }$, the measures $\{\mu _{\xi _i}-\mu
_{\zeta _i}\}_{i\in I}$ are pairwise $\frac \delta 2$-singular. Therefore we
assume that $F\neq \emptyset $.\newline
We set $\tau _{\xi _{}}=\underset{\gamma \in F}{\sum }\frac{d\mu _{\xi _{}}}{%
d\mu _{\gamma _{}}}$ for any $\xi \in \Xi ^{\prime }$.\newline
Since $\left( \underset{\gamma \in F}{\sum }\oplus L^1(\mu _{\gamma
_{}})\right) _1$ is a separable Banach space there exists an uncountable
family $\left\{ \{\xi _i,\zeta _i\}\right\} _{i\in I}$, of pairwise disjoint
two-point subsets of $\Xi ^{\prime }$ such that for each $i\in I$

\begin{center}
$||\tau _{\xi _i}-\tau _{\zeta _i}||<\frac \delta 4$.
\end{center}

It is easy to check that the members of the family $\{\mu _{\xi _i}-\mu
_{\zeta _i}:i\in I\}$ are pairwise $\delta $-singular%
\endproof%

In the sequel by a {\bf normalized block} sequence in the space $X$ we
understand a sequence $(z_n)_{n\in \mathbb{N}}$ of vectors of $X$ with finite
supports $(E_n)_{n\in \mathbb{N}}$ and

\begin{center}
$\max \{|\delta |:\delta \in E_n\}<\min \{|\delta |:\delta \in E_{n+1}\}$
\end{center}

\begin{lemma}
Let $(z_n)_{n\in \mathbb{N}}$ be a normalized block sequence in the space $X$
and $\{M_{\xi _{}}\}_{\xi \in \Xi }$ an uncountable almost disjoint family
of subsets of $\mathbb{N}$.For any $\xi \in \Xi $ set $z_{\xi _{}}^{**}=%
\underset{n\in M_{\xi _{}}}{\sum }z_n$, in the $w^{*}$-sense. Then the
family $\{Qz_{\xi _{}}^{**}\}_{\xi \in \Xi }$ is isometrically equivalent to
the usual basis of $c_0(\Xi )$.
\end{lemma}

\begin{proof}
Since $(M_{\xi _{}})_{\xi \in \Xi }$ are almost disjoint and $(z_n)_{n\in 
\mathbb{N}}$ is a normalized block sequence then\newline
$\left\| \underset{i=1}{\overset{n}{\sum }}\lambda _iQz_{\xi
_i}^{**}\right\| =\underset{n\rightarrow \infty }{\lim }\left\| \underset{i=1%
}{\overset{n}{\sum }}\lambda _iz_{\xi _i}^{**}|\{\delta :|\delta |\geq
n\}\right\| =\max \{|\lambda _i|:1\leq i\leq n\}$
\end{proof}

$_{}$

{\bf Completion of the proof of Theorem \ref{5.1}}\newline
For $\varepsilon =\frac 1{16}$ choose a normalized block sequence $%
(z_n)_{n\in \mathbb{N}}$ and a $\lambda >0$ such that

\begin{center}
$B_Z\subset \lambda W+\varepsilon B_X$, where $Z=\overline{<z_n)_{n\in \mathbb{N%
}}>}$.
\end{center}

Let $(N_{\xi _{}})_{\xi \in \Xi }$ be an uncountable almost disjoint family
of subsets of $\mathbb{N}$ and for any $\xi \in \Xi $ define

\begin{center}
$z_{\xi _{}}^{**}=\underset{n\in N_{\xi _{}}}{\sum }z_n$.
\end{center}

Since $(z_{\xi _{}}^{**})_{\xi \in \Xi }$ is a subset of $B_{Z^{**}}$ there
exists a family $\{\mu _{\xi _{}}\}_{\xi \in \Xi }$ in $\mathcal{M}_1(K_c)$
such that

\begin{center}
$||\lambda QR\mu _{\xi _{}}-Qz_{\xi _{}}^{**}||<\varepsilon $.
\end{center}

Since $||Qz_{\xi _{}}^{**}||=1$, it follows that

$1-\varepsilon <||\lambda \mu _{\xi _{}}||=||\lambda QR\mu _{\xi
_{}}||<1+\varepsilon $.

By lemma \ref{5.5} we get an uncountable family $\{\lambda \mu _{\xi
_i}-\lambda \mu _{\zeta _i}\}_{i\in I}$ which is $\frac 1{16}$-singular.
Since \\
$||Q(z_{\xi _i}^{**}-z_{\zeta _i}^{**})||=1$, we get that

\begin{center}
$1-2\varepsilon <||\lambda QL(\mu _{\xi _i}-\mu _{\zeta
_i})||<1+2\varepsilon $
\end{center}

and also we have that

\begin{center}
$||\lambda QL(\mu _{\xi _i}-\mu _{\zeta _i})-(z_{\xi _i}^{**}-z_{\zeta
_i}^{**})||<2\varepsilon $.
\end{center}

Hence

\begin{center}
$1=\left\| \underset{i=1}{\overset{n}{\sum }}Q(z_{\xi _i}^{**}-z_{\zeta
_i}^{**})\right\| \geq \left\| \underset{i=1}{\overset{n}{\sum }}\lambda
QR(\mu _{\xi _i}-\mu _{\zeta _i})\right\| -2\varepsilon \geq $\\[0pt]
$(1-2\varepsilon -2\delta )n-2\varepsilon \geq \left( 1-\frac 14\right)
n-\frac 24$
\end{center}

which leads to a contradiction for $n\geq 4$ and the proof is complete.%
\endproof%

\section{{\bf CONSTRUCTIONS OF BLOCK-H.I. SPACES}}

In this section we present a general method of constructing block-H.I.
Banach spaces; thus for any sequence $\left( X_n\right) _{n\in \mathbb{N}}$ of
separable Banach spaces we define a $d$-product such that the resulting
space is a block H.I. space (Definition \ref{Def2.3} ). They are mainly two
ways to define such a norm. The first is to follow Gowers-Maurey methods and
the second is the method developed in \cite{AD2} using the analysis of the
norming functionals. There are several variations that give more information
but we do not use them since the whole construction becomes more
complicated. Thus, for example, one can define a block-H.I. $d$-product norm
such that every block subspace has to be an asymptotic $\ell ^1$ space. 
\newline
{\bf Notation. }Let $[\mathbb{N}]$ be the space $\{0,1\}^{\mathbb{N}}$ of all
the subsets of $\mathbb{N}$ with the product topology, $[\mathbb{N}\mathbf{]}%
_{<\omega _{}}$ the set of all finite subsets of $\mathbb{N}$. A set $\mathcal{M%
}\subseteq N_{<\omega _{}}$ is {\bf compact}\textit{\ }if\textit{\ }it is
a closed subset of\textit{\ }$[\mathbb{N}]$ and {\bf adequate}\textit{\ }if
whenever\textit{\ }$B\subset A$ and $A\in \mathcal{M}$ then $B\in \mathcal{M}
$.

A sequence $(E_1,...,E_n)$, $n\in \mathbb{N}$, of finite subsets of $\mathbb{N}$
is called $\mathcal{M}$\textit{-}{\bf admissible} if there exists a set $%
\{m_1,...,m_n\}\in \mathcal{M}$ such that\newline
$m_1\leq E_1<m_2\leq E_2<...<m_n\leq E_n$.\newline
A sequence $(\tilde{x}_1,...,\tilde{x}_n)$, $n\in \mathbb{N}$ of vectors of $%
\Omega _{00}=\left( \underset{n=1}{\overset{\infty }{\prod }X_n}\right)
_{00} $ is called $\mathcal{M}$\textit{-}{\bf admissible} if the sequence 
$($\textit{supp}$\tilde{x}_1,...,$\textit{supp}$\tilde{x}_n)$ is $\mathcal{M}
$-admissible.

Let $\theta \in \mathbb{R}$ and $\mathcal{M}\subseteq [\mathbb{N}]_{<w}$. A subset 
$K$ of $\Omega _{00}$ is called $(\mathcal{M},\theta )$\textit{-}{\bf %
closed} if for any $\mathcal{M}$-admissible sequence $(\tilde{x}_1,..,\tilde{%
x}_n)$ of elements of $K,\theta \underset{i=1}{\overset{n}{\sum }}\tilde{x}%
_i\in K$. If $A\subseteq \Omega _{00}$ and $(\mathcal{M}_n,\theta _n)_{n\in 
\mathbb{N}}$ is a sequence such that $\mathcal{M}_n\subseteq [\mathbb{N}]_{<\omega
_{}}$ and $\theta _n\in \mathbb{R}$ for every $n\in \mathbb{N}$, $K(A,(\mathcal{M}%
_n,\theta _n)_{n\in \mathbb{N}})$ is the smallest subset of $\Omega _{00}$
which contains $A$ and it is $(\mathcal{M}_n,\theta _n)$-closed for every $%
n\in \mathbb{N}$.\newline
In the sequel we consider families $\mathcal{M}_n\subseteq [\mathbb{N}%
]_{<\omega _{}}$ which are compact, adequate and $\theta _n\searrow 0$.

Let $(X_n,||.||)_{n\in \mathbb{N}}$ be a sequence of separable Banach spaces.
For each $n\in \mathbb{N}$ let $\mathcal{F}_n\subseteq X_n^{*}$ such that:%
\newline
1) $\mathcal{F}_n$ is a countable symmetric subset of $B_{X_n^{*}}$\newline
2) For every $x\in X_n$ $||x||=\underset{f\in \mathcal{F}_n}{\sup }f(x)$.%
\newline
Every $f=(f_n)_{n\in \mathbb{N}}\in \left( \underset{n=1}{\overset{\infty }{%
\prod }}X_n^{*}\right) _{00}$ is a functional $f:\underset{n=1}{\overset{%
\infty }{\prod }}X_n\rightarrow \mathbb{R}$ defined by the rule $f(\tilde{x})=%
\underset{n=1}{\overset{\infty }{\sum }}f_n(x_n)$, where $\tilde{x}%
=(x_n)_{n\in \mathbb{N}}\in \underset{n=1}{\overset{\infty }{\prod }}X_n$.
Every $f\in X_n^{*}$ is coincides with the element $(0,0,...,f,0,...)$ of $%
\left( \underset{n=1}{\overset{\infty }{\prod }}X_n^{*}\right) _{00}$.

\begin{definition}
The $(\mathcal{M}_n,\theta _n)_{n\in \mathbb{N}}$\textit{-d-product} of a
sequence $(X_n,||.||_n)_{n\in \mathbb{N}}$ of separable spaces is the
completion of $\left( \left( \underset{n=1}{\overset{\infty }{\prod }}%
X_n\right) _{00},||.||_K\right) $, where $||\tilde{x}||_K=\underset{f\in K}{%
\sup }f(\tilde{x})$ and $K=K\left( \underset{n\in \mathbf{N}}{\bigcup }%
\mathcal{F}_n,(\mathcal{M}_n,\theta _n)_{n\in \mathbb{N}}\right) $.
\end{definition}

If for every $n\in \mathbb{N}$ $\dim X_n=1$ then the $(\mathcal{M}_n,\theta
_n)_{n\in \mathbb{N}}$-$d$-product coincides with the mixed-Tsirelson space $%
T\left( (\mathcal{M}_n,\theta _n)_{n\mathbf{\in }}\right) $ introduced in
[AD2], where $\mathcal{F}_n=\{e_n^{*},-e_n^{*}\}$, $e_n^{*}(e_n)=1$ for $%
e_n\in X_n$ such that $||e_n||_n=1$.\newline
The sequence $(e_n)_{n=1}^\infty $ is a 1-unconditional basis of $T\left( (%
\mathcal{M}_n,\theta _n)_{n\mathbf{\in }\mathbb{N}}\right) $.

\begin{remark}
It is easy to see that the set $K$ in the above definition is also defined
inductively in the following manner
\end{remark}

\begin{center}
$K=\underset{s=0}{\overset{\infty }{\bigcup }}K^s$ where $K^0=\underset{n\in 
\mathbf{N}}{\bigcup }\mathcal{F}_n$ and if $s>1$ $K^s=K^{s-1}\bigcup \left( 
\underset{i=1}{\overset{\infty }{\bigcup }}K_i^n\right) $ where\\[0pt]
$K_i^n=\left\{ \theta _i\underset{k=1}{\overset{d}{\sum }}f_k:(f_1,...,f_d)%
\text{ is }\mathcal{M}_i\text{-admissible and }f_k\in K^{s-1}\right\} $.
\end{center}

We denote by

$A_0^{*}=K^0$ and if $n>0$ $A_n^{*}=\underset{s=0}{\overset{\infty }{\bigcup 
}}K_n^s$. For $n>0$ the set $A_n^{*}$ is described as.

\begin{center}
$\left\{ \theta _n\underset{k=1}{\overset{d}{\sum }}f_k:d\in \mathbb{N},\text{ }%
(f_1,...,f_d)\text{ is }\mathcal{M}_n\text{-admissible and }f_k\in K\right\} 
$.
\end{center}

For every $f\in K$ we define

\begin{center}
$j(f)=\max \{n:f\in A_n^{*}\}$
\end{center}

As in the case of previously defined H.I-spaces (\cite{GM}, \cite{AD2}) we
need an appropriate auxiliary function that we define below\newline
We denote by \newline
$\mathcal{K}=\{(f_1,...,f_d):d\in \mathbb{N},$ $f_1<...<f_d,$ $f_i\in K,$ $%
i=1,...,d\}$.\newline
Notice that the set $\mathcal{K}$ is countable.\newline
There exists $\sigma :\mathcal{K}\rightarrow \{2j:j\in \mathbb{N}\}$ such that%
\newline
a) $\sigma (f)>j(f)$ for any $f\in K$\newline
b) If $(f_1,...,f_n,f_{n+1})\in \mathcal{K}$ then $\sigma
(f_1,...,f_n)<\sigma (f_1,...,f_n,f_{n+1})$\newline
c) $\sigma $ is a one to one function.\newline
The injection $\sigma $ defines a subset $L$ of $K$ as follows:

$L=\underset{s=0}{\overset{\infty }{\bigcup }}L^s$ where\newline
(i) $L^0=K^0=\underset{n\in \mathbf{N}}{\bigcup }\mathcal{F}_n$, $%
L_j^0=\emptyset $ for every $j\in \mathbb{N}$\newline
(ii) If $s>0$ $L^s=\underset{s=0}{\overset{\infty }{\bigcup }}L_i^s$ where

$L_{2j}^s=\pm L_{2j}^{s-1}\bigcup \{\theta _{2j}(f_1+...+f_d):d\in \mathbb{N},$ 
$(f_1,...,f_d)$ is a $\mathcal{M}_{2j}$-admissible sequence of elements of $%
L^{s-1}\}$

$L_{2j+1}^{^{\prime }s}=\pm L_{2j+1}^{s-1}\bigcup \{\theta
_{2j+1}(f_1+...+f_d):d\in \mathbb{N},$ $(f_1,...,f_d)$ is $\mathcal{M}_{2j+1}$%
-admissible sequence of elements of $L^{s-1}$, $f_1\in L_{2k}^{s-1}$ with $%
k>n_{2j+1}$ and $f_{i+1}\in L_{\sigma (f_1,...,f_i)}^{s-1}$ for $i\geq 1\}$,
where $(n_j)_{j=1}^\infty $ will be a strictly increasing sequence of
positive integers.

$L_{2j+1}^s=\{f|_{[k,+\infty )}:f\in L_{2j+1}^{^{\prime }s}$ and $k\in \mathbb{N%
}\}$.\newline
We set $B_i^{*}=\underset{s=1}{\overset{\infty }{\bigcup }}L_i^s$.

\begin{definition}
The $(\mathcal{M}_n,\theta _n)_{n\in \mathbb{N}}^{\sigma _{}}${\bf -product}
of a sequence $(X_n,||.||_n)_{n\in \mathbb{N}}$ of separable spaces is the
completion of $\left( \left( \underset{n=1}{\overset{\infty }{\prod }}%
X_n\right) _{00},||.||_L\right) $.
\end{definition}

\begin{remark}
The norms $||.||_K,||.||_L$ defined above by the sets $K,L$ satisfy the
implicit relations:\newline
$||\tilde{x}||_K=\max \left\{ ||\tilde{x}||_\infty ,\underset{K\in \mathbf{N}%
}{\sup }\left\{ \sup \theta _k\underset{i=1}{\overset{d}{\sum }}||E_i\tilde{x%
}||_K\right\} \right\} $\newline
where the second ``$\sup $'' is over all the $\mathcal{M}_k$-admissible
sequences $(E_1,...,E_d),$ $d\in \mathbb{N}$ of intervals and\newline
$||\tilde{x}||_L=\max \left\{ ||\tilde{x}||_\infty ,\underset{K\in \mathbf{N}%
}{\sup }\left\{ \sup \theta _{2k}\underset{i=1}{\overset{d}{\sum }}||E_i%
\tilde{x}||_L\right\} ,\sup (f(x))\right\} $,\newline
where $f\in \underset{s=1}{\overset{\infty }{\bigcup }}B_{2s-1}^{*}$ $%
(E_1<...<E_d)\mathcal{M}_{2k}$-admissible sequence of intervals.
\end{remark}

\begin{remark}
Let $i(\mathcal{M})$ be the Cantor-Bendixson index of $\mathcal{M}$. If $i(%
\mathcal{M})\geq \omega $ we set $\frac 1{i(\mathcal{M})}=0$. A modification
of the arguments of Proposition 1.1 in \cite{AD2} proves the following:
\end{remark}

\begin{proposition}
\label{6.1}Let $(X_n,||.||_n)_{n\in \mathbb{N}}$ be a sequence of separable
Banach spaces\newline
(a) The norm $(\mathcal{M}_n,\theta _n)_{n\in \mathbb{N}}^{\sigma _{}}$ on $%
\underset{n=1}{\overset{\infty }{\prod }}X_n$ is shrinking.\newline
(b) If there exists $n\in \mathbb{N}$ such that $\theta _n>\frac 1{i(\mathcal{M}%
_n)}$ then this norm is also boundedly complete.\newline
(c) If (b) is satisfied then the $(\mathcal{M}_n,\theta _n)^{\sigma _{}}$%
-product is a reflexive space if and only if every $X_n$ is reflexive.%
\endproof%
\end{proposition}

Our goal is to show that for appropriate choice of sequences $(\mathcal{M}%
_n)_{n\in \mathbb{N}}$, $(\theta _n)_{n\in \mathbb{N}}$ we get that the $(\mathcal{%
M}_n,\theta _n)_{n\in \mathbb{N}}^{\sigma _{}}$-product norm on $\underset{n=1}{%
\overset{\infty }{\prod }}X_n$ is a block H.I-space for every sequence $%
(X_n,||.||_n)_{n\in \mathbb{N}}$ of separable Banach spaces.\newline
There are several choices of such sequences. We will use the following two
sequences that enable us to give a relative simple proof.\newline
We denote by $\mathcal{A}_n$ the compact family \newline
$\mathcal{A}_n=\{F\subset \mathbb{N}:|F|\leq n\}$.\newline
We select two sequences of natural numbers in the following manner.\newline
Set $m_1=2$ and inductively $m_{i+1}\geq m_i^5$. Next we set $n_i=3$, $%
n_{i+1}=n_i^{s_i}$, where $s_i$ is an integer such that $2^{s_i}>m_{i+1}^3$.%
\newline

\begin{theorem}
\label{6.2}Let $(X_n,||.||_n)_{n\in \mathbb{N}}$ be a sequence of separable
Banach spaces. Then the $\left( \mathcal{A}_{n_i},\frac 1{m_i}\right)
^{\sigma _{}}$-$d$-product of the sequence $(X_n,||.||_n)_{n\in \mathbb{N}}$ is
a block H.I-space.
\end{theorem}

{\bf Notation.} We denote $(\tilde{X}_j,|.|_j)$ the norm defined on $%
\underset{n=1}{\overset{\infty }{\prod }}X_n$ by the set of functionals $%
K_j=K\left( \underset{n\in \mathbf{N}}{\bigcup }\mathcal{F}_n,(\mathcal{A}%
_{n_i},\frac 1{m_i})_{i=1}^j\right) $ we first prove the following
proposition that gives some upper $p$-estimates for the $|.|_j$ norm.

\begin{proposition}
\label{6.3}For every $j=1,2,...$ and $\tilde{x}\in \tilde{X}_j$, $|\tilde{x}%
|_j\leq ||\tilde{x}||_p$ where $p=\min \{p_i:1\leq i\leq j\}$ and $p_i=\frac
1{1-\log _{n_i}m_i}$.
\end{proposition}

\begin{proof}
It is enough to show that $|f(\tilde{x})|\leq ||\tilde{x}||_p$ for every $%
f\in K_j=\overset{\infty }{\underset{s=0}{\bigcup }}K_j^s$.\newline
We will prove it by induction.\newline
For $s=0$ and $f\in K_j^0$ we have that $|f(\tilde{x})|\leq ||\tilde{x}%
||_\infty \leq ||\tilde{x}||_p$.Assume that the statement has been proved
for $0,1,...,s$ and let $f\in K_j^{s+1}\backslash K_j^s$. Then $f=\frac
1{m_k}(f_1+...+f_d)$ for some $1\leq k\leq j$ and $f_1,...,f_d$ successive
elements of $K_j^s,$ $d\leq n_k$.\newline
If $\frac 1{q_k}=1-\frac 1{p_k}$ then $\frac 1{q_k}=\log _{n_k}m_k$ so $%
m_k=n_k^{1/q_k}$. Let $E_i=ranf_i$ for $i=1,...,d$. Then\newline
$|f(\tilde{x})|\leq \frac 1{m_k}\underset{i=1}{\overset{d}{\sum }}|f_i(%
\tilde{x})|\leq \frac 1{m_k}\underset{i=1}{\overset{d}{\sum }}||E_i\tilde{x}%
||_p\leq \frac 1{m_k}d^{1/q_k}\left( \underset{i=1}{\overset{d}{\sum }}||E_i%
\tilde{x}||_p^{p_k}\right) ^{1/p_k}\leq $\newline
$\leq \left( \frac d{n_k}\right) ^{1/q_k}\left( \underset{i=1}{\overset{d}{%
\sum }}||E_i\tilde{x}||_p^p\right) ^{1/p}\leq ||\tilde{x}||_p$.
\end{proof}

\begin{definition}
Let $K=K\left( \underset{n\in \mathbf{N}}{\bigcup }\mathcal{F}_n,(\mathcal{M}%
_n,\theta _n)_{n\in \mathbb{N}}\right) $ and $\phi \in K$. An {\bf analysis}
of $\phi $ is a sequence $(K^s(\phi ))_{s=0}^m$ such that\newline
(1) $K^s(\phi )=\{f_1,...,f_{d_s}\}$ where $f_1<...<f_{d_s}$ and $f_i\in K^s$
for $i=1,...,d_s$\newline
(2) $\underset{f\in K^s(\phi )}{\bigcup }$\textit{supp}$f=$\textit{supp}$%
\phi $\newline
(3) If $s>0$ and $f\in K^s(\phi )\backslash K^{s-1}(\phi )$ there exist $%
n\in \mathbb{N}$ and a $\mathcal{M}_n$-admissible sequence $(f_1,...,f_d)$ of
elements of $K^{s-1}(\phi )$ such that $f=\theta _n\underset{i=1}{\overset{d%
}{\sum }}f_i$.\newline
(4) $K^m(\phi )=\{\phi \}$.\newline
If $\phi =\frac 1{m_q}(f_1+...+f_d)$, where $f_1,...,f_d\in K^{m-1}(\phi )$
we define $w(\phi )=\frac 1{m_q}$.
\end{definition}

\begin{remark}
(a) It follows from the inductive definition of the set $K$ that every $\phi
\in K$ has an analysis. Further if $L$ is the subset of $K$ that defines the 
$(\mathcal{M}_n,\theta _n)_{n\in \mathbb{N}}^{\sigma _{}}$-product norm on $%
\underset{n=1}{\overset{\infty }{\prod }}X_n$ then every $\phi \in L$ has an
analysis in elements that also belong to $L$.\newline
(b) The analysis of a functional $\phi $ is important for the proof that the
space is block-H.I. Essentially it will permit us to estimate the action of
a functional $\phi \in L$ on a certain block average by using a functional $%
\psi $ acting on an average of the unit vectors of the basis in the
following space $T\left[ \left( \mathcal{A}_{2n_i},\frac 1{m_i}\right)
_{i=1}^\infty \right] $ defined below.
\end{remark}

\begin{definition}
For $(m_i)_{i\in \mathbb{N}}$, $(n_i)_{i\in \mathbb{N}}$ we denote by $T\left[
\left( \mathcal{A}_{2n_i},\frac 1{m_i}\right) _{i=1}^\infty \right] $ the
mixed Tsirelson space defined by these sequences.
\end{definition}

\begin{remark}
Us we mentioned before mixed Tsirelson spaces are special cases of the $(%
\mathcal{M}_n,\theta _n)_{n\in \mathbb{N}}$-$d$-product norm of $\underset{n=1}{%
\overset{\infty }{\prod }}X_n$ where $\dim X_n=1$.\newline
In this special case we denote $B_i^{*}=\underset{i=1}{\overset{\infty }{%
\bigcup }}K_i^s$ which is equal to the set\newline
$\left\{ \frac 1{m_i}(f_1+...+f_d):d\leq 2n_i,f_1,...,f_d\in K\text{ and }%
f_1<...<f_d\right\} $.\newline
Further for $f\in B_i^{*}$ we denote $w(f)=\frac 1{m_i}$ to be the {\bf %
weight} of the functional $f$.
\end{remark}

\begin{lemma}
\label{6.4}For the sequences $(m_i)_{i\in \mathbb{N}}$, $(n_i)_{i\in \mathbb{N}}$
defined above we have that the sequences $p_i^{\prime }=\frac 1{1-\log
_{n_i}m_i}$, $p_i=\frac 1{1-\log _{2n_i}m_i}$, are strictly decreasing to $1$
hence the $T\left[ \left( \mathcal{A}_{2n_i},\frac 1{m_i}\right)
_{i=1}^j\right] $ is isomorphic to $\ell _{p_j}$.
\end{lemma}

\begin{proof}
The first part follows easily follows easily from the choice of $(m_i)_{i\in 
\mathbb{N}}$, $(n_i)_{i\in \mathbb{N}}$ and for the second we refer to \cite{AD1}.
\end{proof}

\begin{lemma}
\label{6.5}Let $k\in \mathbb{N}$ and $x=\frac 1{n_k}\underset{i=1}{\overset{n_k%
}{\sum }}e_i$. Then\newline
(a) For every $q=1,2,...$ and $\phi \in B_q^{*}$\newline
$|\phi (x)|\leq \left\{ 
\begin{array}{c}
\frac 2{m_qm_k}\text{ \thinspace \thinspace \thinspace \thinspace \thinspace
if }q<k \\ 
\frac 1{m_q}\text{ \thinspace \thinspace \thinspace \thinspace \thinspace
\thinspace \thinspace \thinspace \thinspace \thinspace \thinspace \thinspace
if }q\geq k
\end{array}
\right. $\newline
(b) If $\phi \in B_q^{*}$ and there exists an analysis $(K^s(\phi ))_{s=0}^m$
of $\phi $ such that for every $f\in \underset{s=0}{\overset{\infty }{%
\bigcup }}K^s(\phi )$, $f\notin B_k^{*}$ then $|\phi (x)|\leq \frac 1{m_k^3}$%
.
\end{lemma}

\begin{proof}
(a): If $q\geq j$ the estimate is obvious.\newline
Let $q<j$, $\phi =\frac 1{m_q}(f_1+...+f_d)$, $d\leq 2n_q,$ $f_1<...<f_d\in
K $ and $(K^s(\phi ))_{s=0}^m$ be an analysis of $\phi $.\newline
We set\newline
$I=\left\{ 
\begin{array}{c}
i:1\leq i\leq n_k\text{ and there exists }f\in \underset{s=0}{\overset{m}{%
\bigcup }}K^s(\phi ) \\ 
\text{ such that }w(f)<\frac 1{m_k}\text{ and }i\in \mathit{supp}f
\end{array}
\right\} $ and\newline
$J=\{1,...,n_k\}\backslash I$\newline
$x_1=\frac 1{n_k}\underset{i\in I}{\sum }e_i,$ $x_2=\frac 1{n_k}\underset{%
i\in J}{\sum }e_i$.\newline
Clearly $|\phi (x_1)|\leq \frac 1{m_qm_k}$\newline
Also for estimating we observe that\newline
$n_k^{\frac 1{q_{k-1}}}=n_{k-1}^{\frac{s_{k-1}}{q_{k-1}}}=\frac{(2n_{k-1})^{%
\frac{s_{k-1}}{q_{k-1}}}}{2^{\frac{s_{k-1}}{q_{k-1}}}}\geq \frac{%
m_{k-1}^{s_{k-1}}}{2^{s_{k-1}}}\geq \frac{2^{5s_{k-1}}}{2^{s_{k-1}}}\geq
2^{4s_{k-1}}\geq m_k^{12}$\newline
and by lemma \ref{6.4} and Proposition \ref{6.3}\newline
$|\phi (x_2)|\leq |x_2|_{k-1}\leq \frac 1{n_k}\left| \underset{i=1}{\overset{%
n_k}{\sum }}e_i\right| _{k-1}\leq \frac 1{n_k}n_k^{\frac 1{p_{k-1}}}=\frac
1{n_k^{\frac 1{q_{k-1}}}}\leq \frac 1{m_k^3}$.\newline
So $|\phi (x)|\leq |\phi (x_1)|+|\phi (x_2)|\leq \frac 1{m_qm_k}+\frac
1{m_k^3}\leq \frac 2{m_qm_k}$\newline
(b) If $q>k$ then $|\phi (x)|\leq \frac 1{m_{k+1}}<\frac 1{m_k^3}$\newline
\thinspace \thinspace \thinspace \thinspace \thinspace \thinspace \thinspace
\thinspace \thinspace \thinspace \thinspace If $q<k$ then we observe that $%
\phi |_J$ where $J$ is defined as in the case $K$ belongs to the norming set
of $T\left[ \left( \mathcal{A}_{2n_i},\frac 1{m_i}\right)
_{i=1}^{k-1}\right] $ hence\newline
$\left| \phi |_J\left( \frac 1{n_k}\underset{i=1}{\overset{n_k}{\sum }}%
e_i\right) \right| \leq \left| \frac 1{n_k}\underset{i=1}{\overset{n_k}{\sum 
}}e_i\right| _{k-1}\leq \frac 1{m_k^{12}}$
\end{proof}

\begin{definition}
A vector $\tilde{x}\in \left( \underset{n=1}{\overset{\infty }{\prod }}%
X_n\right) _{00}$ is called a $n_j${\bf -average} if there is a sequence $%
\tilde{x}_1<...<\tilde{x}_{n_j}$ of normalized vectors of $\left( \underset{%
n=1}{\overset{\infty }{\prod }}X_n\right) _{00}$ such that $\tilde{x}=\frac
1{n_j}(\tilde{x}_1+...+\tilde{x}_{n_j})$.\newline
A $n_j$-average $\tilde{x}$ is called {\bf seminormalized }if $||\tilde{x}%
||\geq \frac 12$.
\end{definition}

\begin{definition}
1) A sequence $(x_j)_{j=1}^{n_k}$ of successive elements of $\left( 
\underset{n=1}{\overset{\infty }{\prod }}X_n\right) _{00}$ is called a $n_k$-%
{\bf rapidly increasing sequence} or a $n_k$ {\bf RIS} if\newline
(i) For every $j=1,...,n$, $\tilde{x}_j$ is a seminormalized $n_{r_j}$%
-average where $k<r_1<r_2<...<r_k$.\newline
(ii) $||\tilde{x}_j||_{\ell _1}\leq \frac{m_{r_{j+1}}}{m_{r_{j+1}-1}}$ for
every $j=1,2,...$ .\newline
For every $k\in \mathbb{N}${\bf \ }the vector $\tilde{x}=\frac 1{n_k}(\tilde{%
x}_1+...+\tilde{x}_{n_k})$ is called a $n_k${\bf -RIS average}.\newline
2) A sequence $(x_j)_{j=1}^\infty $ of successive elements of $\left( 
\underset{n=1}{\overset{\infty }{\prod }}X_n\right) _{00}$ is called a 
{\bf rapidly increasing sequence} or a {\bf RIS} if\newline
(i) For every $j=1,...,n$, $\tilde{x}_j$ is a seminormalized $n_{r_j}$%
-average where $r_1<r_2<...$.\newline
(ii) $||\tilde{x}_j||_{\ell _1}\leq \frac{m_{r_{j+1}}}{m_{r_{j+1}-1}}$ for
every $j=1,2,...$
\end{definition}

\begin{remark}
The above definitions are similar to the corresponding in Gowers-Maurey's
example \cite{GM}. The existence of seminormalized $n_j$-averages for every
block sequence is established in the following proposition and the proof
follows their method.
\end{remark}

\begin{proposition}
\label{6.6}Let $(\tilde{y}_i)_{i=1}^\infty $ be a normalized block sequence
of $\left( \underset{n=1}{\overset{\infty }{\prod }}X_n\right) _{00}$. For
every $k>$ there exist $\tilde{x}_1,...,\tilde{x}_{n_k}$ normalized blocks
of .$(\tilde{y}_i)_{i=1}^\infty $ such that
\end{proposition}

\begin{center}
$\left\| \frac 1{n_k}\underset{i=1}{\overset{n_k}{\sum }}\tilde{x}_i\right\|
\geq \frac 12$.
\end{center}

\begin{proof}
{\bf Case 1: }$k=2j+1$ for some $j\in \mathbb{N}$.\newline
We define $n_{2j+1}$-averages $\tilde{x}_k^{(t)}$ in $\tilde{Y}=<(\tilde{y}%
_i)_{i=1}^\infty >$ as follows:\newline
For every $k\in \mathbb{N}$ we set $\tilde{x}_k^{(0)}=\tilde{y}_k^{(0)}=\tilde{y%
}_k$ and inductively\newline
$\tilde{x}_k^{(t+1)}=\frac 1{n_{2j+1}}\underset{i=(k-1)n_{2j+1}+1}{\overset{%
kn_{2j+1}}{\sum }}\tilde{y}_i^{(t)}$ and $\tilde{y}_k^{(t+1)}=\frac{\tilde{x}%
_k^{(t+1)}}{||\tilde{x}_k^{(t+1)}||}$.\newline
We claim that for some $t$, $1\leq t\leq s_{2j+1}-1$ and $%
k=1,...,n_{2j+1}^{s_{2j+1}-1}$, $||\tilde{x}_k^{(t)}||\geq \frac 12$.
Indeed, if this not true then\newline
$||\tilde{x}_1^{(s_{2j+1})}||\geq 2^{s_{2j+1}-1}\frac
1{n_{2j+1}^{s_{2j+1}}}||\underset{i=1}{\overset{n_{2j+1}^{s_{2j+1}}}{\sum }}%
\tilde{y}_i||=\frac{2^{s_{2j+1}-1}}{n_{2j+2}}||\underset{i=1}{\overset{%
n_{2j+2}}{\sum }}\tilde{y}_i||\geq \frac{2^{s_{2j+1}-1}}{m_{2j+2}}>1$, which
is a contradiction.\newline
{\bf Case 2:} $k=2j$ for some $j\in \mathbb{N}$.\newline
We define $n_{2j}$-averages $\tilde{x}_k^{(t)}$ in $\tilde{Y}=<(\tilde{y}%
_i)_{i=1}^\infty >$ as in previous case and then for some $1\leq t\leq
s_{2j}s_{2j+1}-1$ and $k=1,...,n_{2j}^{s_{2j}s_{2j+1}-1}$, $||\tilde{x}%
_k^{(t)}||\geq \frac 12$ because if this is not true then\newline
$||\tilde{x}_1^{s_{2j}s_{2j+1}}||\geq 2^{s_{2j}s_{2j+1}-1}\frac 1{n_{2j+2}}||%
\underset{i=1}{\overset{n_{2j+2}}{\sum }\tilde{y}_i||\geq \frac{%
2^{s_{2j}s_{2j+1}-1}}{m_{2j+2}}}>1$\newline
which is a contradiction.
\end{proof}

\begin{remark}
(a) From the above Proposition \ref{6.6} we get that for every $(\tilde{y}%
_i)_{i=1}^\infty $ a normalized block sequence of vectors of $\left( 
\underset{n=1}{\overset{\infty }{\prod }}X_n\right) _{00}$ there exists a
sequence $(\tilde{x}_i)_{i=1}^{n_k}$ of successive blocks of $(\tilde{y}%
_i)_{i=1}^\infty $ which is a $n_k$-RIS.\newline
(b) Our strategy in order to prove that the space $\tilde{X}$ is a block
H.I-space is the following.\newline
In the first step we will show that every $n_k$-RIS average behaves like the 
$n_k$-average of the vectors of the basis in the space $T\left[ \left( 
\mathcal{A}_{2n_i},\frac 1{m_i}\right) _{i=1}^\infty \right] $. Then we
shall introduce the dependent sequences which will give us the desired
result. To estimate the norms on the averages of the dependent sequences we
will use again the space $T\left[ \left( \mathcal{A}_{2n_i},\frac
1{m_i}\right) _{i=1}^\infty \right] $. The main tool for passing the
estimates from one space to the other is the analysis of the functionals
that we introduced above.
\end{remark}

\begin{lemma}
\label{6.7}Let $j\geq 2$ and $\tilde{x}\in \left( \underset{n=1}{\overset{%
\infty }{\prod }}X_n\right) _{00}$ be a $n_j$-average. Then for every $n\leq
n_{j-1}$ and $E_1<...<E_n$ successive intervals of integers $\underset{i=1}{%
\overset{n}{\sum }}||E_i\tilde{x}||\leq 2$.
\end{lemma}

\begin{proof}
Let $\tilde{x}=\frac 1{n_j}\underset{k=1}{\overset{n_j}{\sum }}\tilde{x}_k$
where $\tilde{x}_1<...<\tilde{x}_{n_j}$ and $||\tilde{x}_k||=1$ for every $%
k=1,...,n_j$. For every $i=1,...,n$ define

\begin{center}
$A_i=\{k:1\leq k\leq n_j$ and \textit{supp}$\tilde{x}_k\subset E_i\}$

$B_i=\{k:1\leq k\leq n_j$ and \textit{supp}$\tilde{x}_k\cap E_i\neq
\emptyset \}$
\end{center}

The sets $A_i$ are pairwise disjoint and $|B_i|\leq |A_i|+2$\newline
Since the norm is bimonotone, for every $i=1,...,n$

\begin{center}
$||E_i\tilde{x}||\leq \frac 1{n_j}||\underset{k\in B_i}{\sum }x_k||\leq
\frac 1{n_j}|B_i|\leq \frac 1{n_j}(|A_i|+2)$
\end{center}

Therefore

\begin{center}
$||\underset{i=1}{\overset{n}{\sum }}E_i\tilde{x}||\leq \frac 1{n_j}\left( 
\underset{i=1}{\overset{n}{\sum }}|A_i|+2n\right) \leq \frac
1{n_j}(n_j+2n_{j-1})\leq 2$.
\end{center}

The proof is complete.
\end{proof}

\begin{proposition}
\label{6.8}Let $k\geq 2$ and $(x_j)_{j=1}^{n_k}$ be a $n_k$-RIS with $x_j$
be a $n_{r_j}$-seminormalized average and $\phi \in K$ such that $w(\phi
)=\frac 1{m_r}$. Then there exists a functional $\psi \in K_T$ such that for
every $j=1,...,n_k$ $|\phi (\tilde{x}_j)|\leq 2\psi (e_j)$.\newline
Moreover if $r<r_1$ then $\psi \in B_r^{*}$ and if $r_t\leq r<r_{t+1}$ then $%
\psi =\frac 12(\psi _1+e_t^{*})$ where $\psi _1\in B_{r-1}^{*}$ and $t\notin 
$\textit{supp}$\psi _1$.
\end{proposition}

\begin{proof}
Let $(K^s(\phi ))_{s=0}^m$ be an analysis of $\phi $. For every $f\in 
\underset{s=0}{\overset{m}{\bigcup }}K^s(\phi )$ we set\newline
\begin{equation*}
D_f=\{j:\mathit{supp}\phi \cap \mathit{supp}\tilde{x}_j=\mathit{supp}f\cap 
\mathit{supp}\tilde{x}_j\neq \emptyset \}.
\end{equation*}
\newline
Using induction on $s=0,...,m$ we shall define for every $f\in \underset{s=0%
}{\overset{m}{\bigcup }}K^s(\phi )$ a functional $g_f\in K_T$ with the
following properties.\newline
(a) \textit{supp}$g_f\subset D_f$\newline
(b) For every $j\in D_f$, $|f(\tilde{x}_j)|\leq 2g_f(e_j)$\newline
(c) $g_f\in K_T$\newline
Moreover $g_f\in (B_r^{*})$ if $r<r_1$ and if $r_t\leq r<r_{t+1}$ $%
^{}g_f=\frac 12(g_f^1+e_t^{*})$ and $t\notin $\textit{supp}$g_f^1$.\newline
If $f\in K^0(\phi )$ and $D_f\neq \emptyset $ then $D_f=\{k\}$ for some $k$, 
$1\leq k\leq n$. In this case we set $g_f=e_k^{*}$.\newline
Let $s>0$ and suppose that $g_f$ has been defined for every $f\in \underset{%
i=0}{\overset{s-1}{\bigcup }}K^i(\phi )$ and let $f\in K^s(\phi )\backslash
K^{s-1}(\phi )$. Then $f=\frac 1{m_q}(f_1+...+f_d)$, $d\leq n_q$, $%
f_1,...,f_d$ successive functionals of $K$.\newline
{\bf Case 1:} $q<r_1$\newline
Define$^{}$ $I=\{i:1\leq i\leq d$ and $D_{f_i}\neq \emptyset \}$ and $%
^{}T=D_f\backslash \underset{i\in I}{\bigcup }D_{f_i}$\newline
We set 
\begin{equation*}
g_f=\frac 1{m_q}\left( \underset{i\in I}{\sum }g_{f_i}+\underset{j\in T}{%
\sum }e_j^{*}\right)
\end{equation*}
\newline
Property (a) is obvious. Since $|I|\leq n_q$ and $|T|\leq n_q$. We get that $%
g_f\in K_T$ hence property (c).\newline
To prove (b) let $j\in D_f$. If $j\in D_{f_i}$ for some $i=1,...,d$ then%
\newline
$|f(\tilde{x}_j)|=\frac 1{m_q}|f_i(\tilde{x}_j)|\leq \frac
2{m_q}g_{f_i}(e_j)=2g_f(e_j)$.\newline
If $j\in T$, let $E_i=ranf_i$ for $i=1,...,d$. Then by lemma \ref{6.7}%
\newline
$|f(\tilde{x}_j)|\leq \frac 1{m_q}\underset{i=1}{\overset{d}{\sum }}||E_i%
\tilde{x}_j||\leq \frac 2{m_q}=\frac 2{m_q}e_j^{*}(e_j)=2g_f(e_j)$\newline
{\bf Case 2:} $q\geq r_1$\newline
Let $t$ the integer with the property $r_t\leq q<r_{t+1}$ for $t=1,...,n_k$.%
\newline
{\bf Subcase 2a:} $t\notin D_f$ or $t\in \underset{i\in I}{\bigcup }%
D_{f_i}$.\newline
We set 
\begin{equation*}
g_f=\frac 1{m_{q-1}}\left( \underset{i\in I}{\sum }g_{f_i}+\underset{j\in T}{%
\sum }e_j^{*}\right) \newline
\end{equation*}
$g_f$ satisfies properties (a) and (c) since $|I|+|T|\leq n_k<n_{r_1-1}\leq
n_{q-1}$.\newline
To prove (b), let $j\in D_f$. If $j\in D_{f_i}$ for some $i=1,...,d$ then 
\begin{equation*}
|f(\tilde{x}_j)|=\frac 1{m_q}|f_i(\tilde{x}_j)|\leq \frac
2{m_q}g_{f_i}(e_j)<\frac 1{m_{q-1}}g_{f_i}(e_j)=g_f(e_j).
\end{equation*}
If $j\in T$ and $j<t$ then $r_{j+1}\leq r_t\leq q$ and by the property (ii)
of RIS we get\newline
\begin{eqnarray*}
|f(\tilde{x}_j)| &\leq &\frac 1{m_q}\underset{i=1}{\overset{d}{\sum }}|f_i(%
\tilde{x}_j)|\leq \frac 1{m_q}||\tilde{x}_j||_{\ell _1}\leq \frac 1{m_q}%
\frac{m_{r_{j+1}}}{m_{r_{j+1}-1}}\leq \frac 1{m_{q-1}}= \\
&=&\frac 1{m_{q-1}}e_j^{*}(e_j)=g_f(e_j)
\end{eqnarray*}
\newline
If $j\in T$ and $t<j$ then by Lemma \ref{6.7} we get\newline
$|f(\tilde{x}_j)|=\frac 1{m_q}|\underset{i=1}{\overset{d}{\sum }}f_i(\tilde{x%
}_j)|\leq \frac 1{m_q}\underset{i=1}{\overset{d}{\sum }}||E_i\tilde{x}%
_j||\leq \frac 1{m_q}<\frac 1{m_{q-1}}e_j^{*}(e_j)=g_f(e_j)$\newline
{\bf Subcase 2b:} $t\in T$. Then we set\newline
\begin{equation*}
g_f^1=\frac 1{m_{q-1}}\left( \underset{i\in I}{\sum }g_{f_i}+\underset{j\in
T\backslash \{t\}}{\sum }e_j^{*}\right)
\end{equation*}
\newline
as in subcase 2a, $|f(\tilde{x}_j)|\leq g_f^1(e_j)$ for every $j\in
D_f\backslash \{t\}$.\newline
For $j\in T$, $j\neq t$ by the previous subcase we get that $|f(\tilde{x}%
_j)|\leq 2g_f(e_j)$.\newline
For $j\in D_{f_i}$ then $|f(\tilde{x}_j)|\leq \frac 1{m_q}|f_i(\tilde{x}%
_j)|\leq g_{f_i}(e_j)\leq \frac 2{m_{q-1}}g_f(e_j)$.\newline
We set $g_f=\frac 12(g_f^1+e_t^{*})$. Then $|f(\tilde{x}_t)|\leq 1=2\frac
12e_t^{*}(e_t)=2g_f(e_t)$. This completes the proof.%
\endproof%
\end{proof}

\begin{proposition}
\label{6.9}Let $(\tilde{x}_j)_{j=1}^{n_k}$ be a $n_k$-RIS of elements of $%
\left( \underset{n=1}{\overset{\infty }{\prod }}X_n\right) _{00}$, $\phi \in
K$ such that $w(\phi )=\frac 1{m_r}$ and $\tilde{x}=\frac 1{n_k}(\tilde{x}%
_1+...+\tilde{x}_{n_k})$. Then\newline
(a) $|\phi (\tilde{x})|\leq \left\{ 
\begin{array}{c}
\frac
4{m_km_r}\,\,\,\,\,\,\,\,\,\,\,\,\,\,\,\,\,\,\,\,\,\,\,\,\,\,\,\,\,\,\,\,\,%
\,\,\,\,\, \\ 
\frac{\overset{}{2}}{m_r}\,\,\,\,\,\,\,\,\,\,\,\,\,\,\,\,\,\,\,\,\,\,\,\,\,%
\,\,\,\,\,\,\,\,\,\,\,\,\,\, \\ 
\frac{\overset{}{2}}{m_{r-1}}+\frac 2{n_k}\leq \frac 4{m_k}
\end{array}
\right. \left. 
\begin{array}{c}
\underset{}{\text{if}}\text{ \thinspace \thinspace \thinspace \thinspace
\thinspace \thinspace \thinspace \thinspace \thinspace \thinspace \thinspace
\thinspace \thinspace \thinspace \thinspace \thinspace \thinspace \thinspace 
}r<k \\ 
\text{if }k\leq r<r_1 \\ 
\text{if \thinspace \thinspace \thinspace \thinspace \thinspace \thinspace
\thinspace \thinspace \thinspace \thinspace \thinspace \thinspace \thinspace
\thinspace }r_1\leq r
\end{array}
\right. $\newline
(b) If $k=2k^{\prime }$ for some $k^{\prime }\in \mathbb{N}$
\end{proposition}

\begin{center}
$\frac 1{2m_k}\leq ||\tilde{x}||\leq \frac 4{m_k}$
\end{center}

\begin{proof}
(a) it follows from the above Proposition and Lemma \ref{6.5} (a)\newline
(b) The fact $||\tilde{x}||\leq \frac 4{m_k}$ follows from (a) If $k$ is an
even there exist successive functionals $f_j\in L$ such that $f_j(x_j)\geq
\frac 12$ for every $j=1,...,n_k$ hence the functional $f=\frac
1{m_k}(f_1+...+f_{n_k})\in L$. So $||\tilde{x}||\geq f(\tilde{x})\geq \frac
1{2m_k}$.
\end{proof}

We pass now to the final stage of the proof. In the following we introduce
the $(2j+1)$-dependent sequences that we will use to prove that the space is
block-H.I.

\begin{definition}
\label{D6.7}Let $j\geq 2$. A sequence $(\tilde{y}_k)_{k=1}^{n_{2j+1}}$ of
successive vectors of $\left( \underset{n=1}{\overset{\infty }{\prod }}%
X_n\right) _{00}$ is called a $\mathbf{(2j+1)}${\bf -dependent sequence}
if there exists a sequence $(y_k^{*})_{k=1}^{n_{2j+1}}$ of successive
elements of $L$ such that for every $k=1,...,n_{2j+1}$\newline
(a) $y_k^{*}\in B_{2r_k}$\newline
(b) $n_{2j+1}<r_1<...<r_{n_{2j+1}}$\newline
(c) $2r_k=\sigma (y_1^{*},...,y_{k-1}^{*})$ for $k=2,...,n_{2j+1}$\newline
(d) Each $\tilde{y}_k$ is a $n_{2r_k}$-RIS average\newline
(e) \textit{supp}$\tilde{y}_k^{*}\subset [\min $\textit{supp}$\tilde{y}_k$, $%
\max $\textit{supp}$\tilde{y}_k]=range(\tilde{y}_k)$ and
\end{definition}

\begin{center}
$\frac 1{2m_{2r_k}}\leq \tilde{y}_k^{*}(\tilde{y}_k)\leq \frac 4{m_{2r_k}}$.
\end{center}

\textit{The sequence} $(y_k^{*})_{k=1}^{n_{2j+1}}$ \textit{is called the 
{\bf dual} sequence of} $(\tilde{y}_k)_{k=1}^{n_{2j+1}}$.\newline
\textit{Given a }$(2j+1)$-\textit{dependent sequence} $(\tilde{y}%
_k)_{k=1}^{n_{2j+1}}$ \textit{for every} $k$ \textit{we set }$\theta _k$ 
\textit{to be} \textit{the positive number such that} $y_k^{*}(\theta
_km_{2r_k}\tilde{y}_k)=1$. \textit{Clearly} $\frac 12\leq \theta _k\leq $5.

\begin{remark}
For every $j>2$ and for every normalized block sequence $(\tilde{x}%
_i)_{i=1}^\infty $ there exists a $(2j+1)$-dependent sequence $(\tilde{y}%
_k)_{k=1}^{n_{2j+1}}$ such that each $\tilde{y}_k$ is a block vector of $(%
\tilde{x}_i)_{i=1}^\infty $.
\end{remark}

\begin{lemma}
\label{6.10}Let $(\tilde{y}_k)_{k=1}^{n_{2j+1}}$ be a $(2j+1)$-depended
sequence and $(m_{2r_k})_{k=1}^{n_{2j+1}}$, $(\theta _k)_{k=1}^{n_{2j+1}}$
as in definition 5.7. For every sequence $(z_s^{*})_{s=1}^d$ of elements of $%
L$, $d\leq n_{2j+1}$ such that $z_s^{*}\in B_{2t_s}^{*}$ $s=1,...,d$ and $%
n_{2j+2}<t_1<...<t_s$, and for every $k=1,...,n_{2j+1}$ such that $r_k\neq
t_s$ for every $s=1,...,d$
\end{lemma}

\begin{center}
$\left| \left( \underset{s=1}{\overset{d}{\sum }}z_s^{*}\right) (m_{2r_k}%
\tilde{y}_k)\right| \leq \frac 1{n_{2j+1}^2}$
\end{center}

\begin{proof}
$\tilde{y}_k=\frac 1{n_{2r_k}}\underset{j=1}{\overset{n_{2r_k}}{\sum }}%
\tilde{x}_j$ where $(x_j)_{j=1}^{n_{2r_k}}$ is a $n_{2r_k}$-RIS.\newline
Let $s_1$ be the maximal integer with the property $t_{s_1}<r_k$. If $s\leq
s_1$ then by Proposition \ref{6.9} we get $|z_s^{*}(\tilde{y}_k)|\leq \frac
4{m_{2t_s}m_{2r_k}}$ and therefore

\begin{center}
$\left| \left( \underset{s=1}{\overset{d}{\sum }}z_s^{*}\right) (m_{2r_k}%
\tilde{y}_k)\right| \leq 4\underset{s=1}{\overset{s_1}{\sum }}\frac
1{m_{2t_s}}\leq \frac 8{m_{2t_1}}\leq \frac 1{2n_{2j+1}^2}$.
\end{center}

For every $s\geq s_1+1$ set

$D_s=\{j:$\textit{supp}$\tilde{x}_j\cap $\textit{supp}$z_s^{*}=$\textit{supp}%
$\tilde{x}_j\cap $\textit{supp}$\left( \underset{s=s_1+1}{\overset{d}{\sum }}%
z_s^{*}\right) \neq \emptyset \}$.\newline
The sets $D_s$ are pairwise disjoint. We define also\newline
$I=\{s\geq s_1+1:D_s\neq \emptyset \}$ and\newline
$\left| \underset{s=s_1+1}{\overset{d}{\sum }}z_s^{*}(\tilde{y}_k)\right|
\leq \underset{s\in I}{\sum }\left| z_s^{*}\left( \frac 1{n_{2r_k}}\underset{%
j\in D_s}{\sum }\tilde{x}_j\right) \right| +\underset{j\in T}{\sum }\frac
1{n_{2r_k}}\left| \left( \underset{s=s_1+1}{\overset{d}{\sum }}%
z_s^{*}\right) (\tilde{x}_j)\right| $.\newline
By Proposition \ref{6.9} for every $s\in I$

\begin{center}
$\left| z_s^{*}\left( \frac 1{n_{2r_k}}\underset{j\in D_s}{\sum }\tilde{x}%
_j\right) \right| \leq \frac 2{m_{2t_s-1}}+\frac 2{n_{2r_k}}$.
\end{center}

So

\begin{center}
$\underset{s\in I}{\sum }\left| z_s^{*}\left( \frac 1{n_{2r_k}}\underset{%
j\in D_s}{\sum }\tilde{x}_j\right) \right| \leq \underset{s=s_1+1}{\overset{d%
}{\sum }}\frac 2{m_{2t_s-1}}+\frac{2n_{2j+1}}{n_{2r_k}}$.
\end{center}

If $j\in T$, $\tilde{x}_j$ is a $n_{l_j}$-average and $l_j>2j+1$ and by
Lemma \ref{6.7}\newline
$\left| \left( \underset{s=s_1+1}{\overset{d}{\sum }}z_s^{*}\right) (\tilde{x%
}_j)\right| \leq 2$. So

\begin{center}
$\underset{j\in T}{\sum }\frac 1{n_{2r_k}}\left| \left( \underset{s=s_1+1}{%
\overset{d}{\sum }}z_s^{*}\right) (\tilde{x}_j)\right| \leq \frac{2n_{2j+1}}{%
n_{2r_k}}$.
\end{center}

Putting all this estimates together

\begin{center}
$\left| \underset{s=1}{\overset{d}{\sum }}z_s^{*}(\tilde{y}_k)\right| \leq
\frac 1{2n_{2j+1}^2}+\underset{s=s_1+1}{\overset{d}{\sum }}\frac
2{m_{2t_s-1}}+\frac{4n_{2j+1}}{n_{2r_k}}\leq $

$\leq \frac 1{2n_{2j+1}^2}+\frac 4{m_{2t_1-1}}+\frac{4n_{2j+1}}{n_{2r_k}}%
\leq \frac 1{n_{2j+1}^2}$
\end{center}
\end{proof}

\begin{lemma}
\label{6.11}Let $(\tilde{y}_k)_{k=1}^{n_{2j+1}}$ be a $(2j+1)$-dependent
sequence, $(y_k^{*})_{k=1}^{n_{2j+1}}$ its dual and $(\theta
_k)_{k=1}^{n_{2j+1}}$, $(m_{2r_k})_{k=1}^{n_{2j+1}}$ as above\newline
Then for every $\phi \in B_{2j+1}^{*}$
\end{lemma}

\begin{center}
$\phi \left( \underset{k=1}{\overset{n_{2j+1}}{\sum }}(-1)^k\theta _km_{2r_k}%
\tilde{y}_k\right) \leq 1$.
\end{center}

\begin{proof}
Let $\phi \in B_{2j+1}^{*}$. Then there exists a $\sigma $-dependent
sequence $(z_1^{*},...,z_d^{*})$, $d\leq n_{2j+1}$ and an interval $%
E=[l,+\infty )$ such that $\phi =\frac 1{m_{2j+1}}E(z_1^{*}+...+z_d^{*})$.
We set $k_1=\min \{k=1,...,d:Ez_k^{*}\neq \emptyset \}$ and $k_2=\max
\{k=1,...,d:(z_1^{*},...,z_d^{*})=(y_1^{*},...,y_d^{*})\}$ otherwise we set $%
k_2=0$.\newline
Then $\phi =\frac
1{m_{2j+1}}(x_{k_1}^{*}+y_{k_1+1}^{*}+...+y_{k_2}^{*}+z_{k_2+1}^{*}+...+z_d^{*}) 
$ where\newline
$x_{k_1}^{*}=Ey_{k_1}^{*},z_{k_2+1}^{*}\neq y_{k_2+1}^{*}$.\newline
We set $\tilde{z}_k=\theta _km_{2r_k}\tilde{y}_k$ and note that $||\tilde{z}%
_k||\leq 16$.\newline
We get the following estimates:\newline
(a) $\left| \phi \left( \underset{k=k_1+1}{\overset{k_2-1}{\sum }}(-1)^k%
\tilde{z}_k\right) \right| \leq \frac 1{m_{2j+1}}\left| \underset{k=k_1+1}{%
\overset{k_2-1}{\sum }}(-1)^ky_k^{*}(\theta _km_{2r_k}\tilde{y}_k)\right| =$%
\newline
$=\frac 1{m_{2j+1}}\left| \underset{k=k_1+1}{\overset{k_2-1}{\sum }}%
(-1)^k\right| \leq \frac 1{m_{2j+1}}$\newline
(b) $\left| \phi (\tilde{z}_{k_1})\right| =\frac 1{m_{2j+1}}\left|
x_{k_1}^{*}(\tilde{z}_{k_1})\right| \leq \frac 1{m_{2j+1}}||\tilde{z}%
_{k_1}||\leq \frac{16}{m_{2j+1}}$\newline
(c) $\left| \phi (\tilde{z}_{k_2})\right| \leq \frac 1{m_{2j+1}}\left|
y_{k_2}^{*}(\tilde{z}_{k_2})\right| +\frac 4{m_{2j+1}}\left| \left( 
\underset{k=k_2+1}{\overset{d}{\sum }}z_k^{*}\right) (\tilde{z}%
_{k_2})\right| $\newline
If $k\geq k_2+1$ then $z_k^{*}\in B_{2t_k}^{*}$,\newline
$2t_k=\sigma (y_1^{*},...,y_{k_1}^{*},z_{k-1}^{*})>\sigma
(y_1^{*})>2n_{2j+1} $\newline
$2t_k=\sigma (y_1^{*},...,y_{k_1}^{*},z_{k-1}^{*})\neq \sigma
(y_1^{*},...,y_{k_{2-1}}^{*})=2r_{k_2}$ and we apply Lemma \ref{6.10} to get%
\newline
$\frac 4{n_{2j+1}}\left| \left( \underset{k=k_2+1}{\overset{d}{\sum }}%
z_k^{*}\right) (\tilde{z}_{k_2})\right| \leq \frac 4{m_{2j+1}}\frac
1{n_{2j+1}^2}$ and so $\left| \phi (\tilde{z}_{k_2})\right| \leq \frac
2{m_{2j+1}}$\newline
(d) $\left| \phi (\tilde{z}_{k_2+1})\right| \leq \frac 1{m_{2j+1}}\left|
z_{k_2+1}^{*}(\tilde{z}_{k_2+1})\right| +\frac 1{m_{2j+1}}\left| \left( 
\underset{k>k_2+1}{\sum }z_k^{*}\right) (\tilde{z}_{k_2+1})\right| \leq 
\frac{17}{m_{2j+1}}$\newline
(by Lemma \ref{6.10})\newline
(e) If $k<k_1$ then $\phi (z_k)=0$\newline
If $k>k_2+1$ then $\left| \phi (\tilde{z}_k)\right| =\frac 1{m_{2j+1}}\left|
\left( \underset{l=k_2+1}{\overset{d}{\sum }}z_l^{*}\right) (\tilde{z}%
_k)\right| $.\newline
Since $l\geq k_2+1$ and $k>k_2+1$ $z_l^{*}\in B_{2t_l}^{*}$ where\newline
$2t_l=\sigma (y_1^{*},...,y_{k_2}^{*},..,z_{l-1}^{*})\neq \sigma
(y_1^{*},...,y_{k_2}^{*},...,y_{k-1}^{*})=2r_k$ and by Lemma \ref{6.10}

\begin{center}
$\left| \phi (\tilde{z}_k)\right| \leq \frac 1{m_{2j+1}}\frac{\theta _k}{%
n_{2j+1}^2}<\frac 1{n_{2j+1}^2}$.
\end{center}

By (a), (b), (c), (d), (e)\newline
$\left| \phi \left( \overset{n_{2j+1}}{\underset{k=1}{\sum }}(-1)^k\tilde{z}%
_k\right) \right| \leq \frac 1{m_{2j+1}}+\frac{16}{m_{2j+1}}+\frac
2{m_{2j+1}}+\frac{17}{m_{2j+1}}+\frac{n_{2j+1}}{n_{2j+1}^2}\leq 1$.
\end{proof}

\begin{lemma}
\label{6.12}Let $\phi \in L$, $(\tilde{y}_k)_{k=1}^{n_{2j+1}}$, $%
(m_{2r_k})_{k=1}^{n_{2j+1}}$, $(\theta _k)_{k=1}^{n_{2j+1}}$ as in
Definition \ref{D6.7} and
\\ $\tilde{x}=\frac{1}{n_{2j+1}}\underset{k=1}{
 \overset{n_{2j+1}}{\sum }}(-1)^k\theta _km_{2r_k}\tilde{y}_k$, then there
exist a $\psi \in K_T$ and an analysis $(K^s(\psi ))_{s=0}^m$ of $\psi $
such that $w(g)\neq \frac 1{m_{2j+1}}$ for every $g\in \underset{s=0}{%
\overset{m}{\bigcup }}K^s(\psi )$ and $(\phi (\tilde{x}%
)|\leq 2\psi \left( \frac 1{n_{2j+1}}\underset{k=1}{\overset{n_{2j+1}}{\sum }%
}e_k\right) $
\end{lemma}

\begin{proof}
Let $\phi \in L$ and $\{K^s(\phi )\}_{s=0}^m$ be an analysis of $\phi $ such
that \newline
$\underset{s=0}{\overset{m}{\bigcup }}K^s(\phi )\subset L$. Let $F_1,...,F_l$
be the maximal elements of $\underset{s=0}{\overset{m}{\bigcup }}K^s(\phi )$
which belong to $B_{2j+1}^{*}$, i.e. there is no $h\in \underset{s=0}{%
\overset{m}{\bigcup }}K^s(\phi )$, such that $w(h)=\frac 1{m_{2j+1}}$ and 
\textit{supp}$F_i\subsetneqq $\textit{supp}$h$ for some $i=1,...,l$. The
functionals $F_1,...,F_l$ have disjoint supports.\newline
For every $f\in \underset{s=0}{\overset{m}{\bigcup }}K^s(\phi )$ we define $%
D_f=\emptyset $ if $f\varsubsetneq F_i$ for some $i=1,...,l$.\newline
Otherwise we define\newline
$D_f=\{k=1,...,n_{2j+1}:$\textit{supp}$\tilde{y}_k\cap $\textit{supp}$f=$%
\textit{supp}$\tilde{y}_k\cap $\textit{supp}$\phi \}$.\newline
For every $f\in \underset{s=0}{\overset{m}{\bigcup }}K^s(\phi )$ we shall
define a functional $g_f\in K_Y$ with the following properties:\newline
(a) supp$g_f\subseteq D_f$\newline
(b) $\left| f\left( \frac 1{n_{2j+1}}\underset{i\in D_f}{\sum }(-1)^i\tilde{z%
}_i\right) \right| \leq 2g_f\left( \frac 1{n_{2j+1}}\underset{k=1}{\overset{%
n_{2j+1}}{\sum }}e_k\right) $\newline
If $f\in K^0(\phi )$ and $D_f=\{k\}$ then we define $g_f=e_k^{*}$. Let $f\in
K^{s+1}(\phi )\backslash K^s(\phi )$ and suppose that $g_f$ has been defined
for every $f\in K^s(\phi )$.\newline
If $w(f)=\frac 1{m_{2j+1}}$ and $D_f\neq \emptyset $ then $f=F_i$ for some $%
i=1,...,l$. Let $k_f=\min \{k:k\in D_f\}$ we define $g_f=e_{k_f}^{*}$.%
\newline
By Lemma \ref{6.11}\newline
$|f(\tilde{x})|\leq \frac 2{n_{2j+1}}=2e_{k_ff}^{*}\left( \frac 1{n_{2j+1}}%
\underset{k=1}{\overset{n_{2j+1}}{\sum }}e_k\right) =2g_f\left( \frac
1{n_{2j+1}}\underset{k=1}{\overset{n_{2j+1}}{\sum }}e_k\right) $.\newline
If $w(f)\neq \frac 1{m_{2j+1}}$ and $D_f\neq \emptyset $ we define $g_f$ as
in the Proposition \ref{6.8}.\newline
Finally $\psi =g_{\phi _{}}$. It follows from the inductive construction of $%
\psi $ that the family $\{g_f:f\in \cup K^s(\phi )$, $D_f\neq \emptyset \}$
defines an analysis of $\psi $ in $K_T$ hence $w(g)\neq \frac 1{m_{2j+1}}$
for every $g\in \cup K^s(\phi )$.
\end{proof}

\begin{proposition}
\label{6.13}Let $(\tilde{y}_k)_{k=1}^{n_{2j+1}}$ be a $(2j+1)$-dependent
sequence, $(y_k^{*})_{k=1}^{n_{2j+1}}$ its dual, $(r_k)_{k=1}^{n_{2j+1}}$, $%
(\theta _k)_{k=1}^{n_{2j+1}}$ as in Definition \ref{D6.7}.Then
\end{proposition}

$^{}$\newline
(a) $\left\| \underset{k=1}{\overset{n_{2j+1}}{\sum }}\frac
1{n_{2j+1}}\theta _km_{2r_k}\tilde{y}_k\right\| \geq \frac 1{m_{2j+1}}$%
\newline
$^{}$\newline
(b) $\left\| \underset{k=1}{\overset{n_{2j+1}}{\sum }}\frac{(-1)^k}{n_{2j+1}}%
\theta _km_{2r_k}\tilde{y}_k\right\| \leq \frac 2{m_{2j+1}^3}$

\begin{proof}
(a) Since $f=\frac 1{m_{2j+1}}(y_1^{*}+...+y_{n_{2j+1}}^{*})\in L$\newline
$\left\| \underset{k=1}{\overset{n_{2j+1}}{\sum }}\frac 1{n_{2j+1}}\theta
_km_{2r_k}\tilde{y}_k\right\| \geq \left| f\left( \underset{k=1}{\overset{%
n_{2j+1}}{\sum }}\frac 1{n_{2j+1}}\theta _km_{2r_k}\tilde{y}_k\right)
\right| =\frac 1{m_{2j+1}}$\newline
(b) It is a consequence of Lemma \ref{6.12} and Lemma \ref{6.5}(b).
\end{proof}

{\bf Proof of the Theorem \ref{6.2}:}

Let $\tilde{Y}$, $\tilde{Z}$ be two block subspaces of $\tilde{X}$ and $%
\varepsilon >0$. Let $j$ such that $\frac 2{m_{2j+1}^2}<\varepsilon $ and $%
r_1>2j+2$. By Proposition \ref{6.6} there is a $\tilde{y}_1\in \tilde{Y}$
which is a $n_{2r_1}$-RIS average and $y_1^{*}\in B_{2r_1}^{*}$ such that $%
y_1^{*}(\tilde{y}_1)\geq \frac 1{2m_{2r_1}}$.\newline
Let $\tilde{y}_2\in \tilde{Z}$ such that \textit{supp}$\tilde{y}%
_2>rangy_1^{*}$ which is a $n_{2r_2}$-RIS average, where $2r_2=\sigma
(y_1^{*})$ and $y_2^{*}\in B_{2r_2}^{*}$ such that $y_2^{*}(\tilde{y}_2)\geq
\frac 1{2m_{2r_2}}$.\newline
In this manner we obtain a $(2j+1)$-dependent sequence $(\tilde{y}%
_k)_{k=1}^{n_{2j+1}}$ such that $\tilde{y}_k\in \tilde{Y}$ if $k$ is an odd
number and $\tilde{y}_k\in \tilde{Z}$ if $k$ is an even number.\newline
Let $(y_k^{*})_{k=1}^{n_{2j+1}}$ its dual sequence and $(\theta
_k)_{k=1}^{n_{2j+1}}$ such that $y_k^{*}(\theta _km_{2r_k}\tilde{y}_k)=1$.%
\newline
If $\tilde{y}=\underset{k\text{ }odd}{\underset{k=1}{\overset{n_{2j+1}}{\sum 
}}}\theta _km_{2r_k}\tilde{y}_k$\newline
$\tilde{z}=\underset{k\text{ }even}{\underset{k=1}{\overset{n_{2j+1}}{\sum }}%
}\theta _km_{2r_k}\tilde{y}_k$ then $\tilde{y}\in \tilde{Y}$, $\tilde{z}\in 
\tilde{Z}$ and by Proposition \ref{6.13}.\newline
$||\tilde{y}+\tilde{z}||\geq \frac 1{m_{2j+1}}$\newline
$||\tilde{y}-\tilde{z}||\leq \frac 2{m_{2j+1}^3}$\newline
so\newline
$||\tilde{y}-\tilde{z}||\leq \frac 2{m_{2j+1}^2}||\tilde{y}+\tilde{z}||\leq
\varepsilon ||\tilde{y}+\tilde{z}||$. and the proof is complete.%
\endproof%

\begin{corollary}
\label{6.14}Let $X$ be a separable Banach space and $W$ a bounded convex
closed symmetric thin subset of $X$. If $||.||_n$ is the equivalent norm on $%
X$ defined by the set $2^nW+2^{-n}B_X$ and $\tilde{X}$ the $\left( \mathcal{A%
}_{n_i},\frac 1{m_i}\right) ^{\sigma _{}}$-product of the sequence $%
(X,||.||_n)_{n\in \mathbb{N}}$ then the diagonal space $\Delta \tilde{X}$ is a
H.I. space.
\end{corollary}

\begin{proof}
It follows from Theorem \ref{6.2} and Proposition \ref{Prop2.3}.
\end{proof}

In the sequel by $T_3$ we denote the mixed Tsirelson space $T\left( \mathcal{%
A}_{3n_i},\frac 1{m_i}\right) _{i\in \mathbb{N}}$, where $(m_i)_i$, $(n_i)_i$
are the sequences used in the definition of $\tilde{X}$

\begin{proposition}
\label{6.15}Let $(\tilde{x}_i)_{i=1}^\infty $ be a RIS in $\tilde{X}$ of $%
(n_{r_i})_i-$averages. Then for every $\phi \in K$ there exists $\psi \in
K_{T_3}$ such that for every $(\lambda _i)_{i=1}^\infty \subset \mathbb{R}$ 
\begin{equation*}
\phi \left( \underset{i=1}{\overset{\infty }{\sum }}\lambda _i\tilde{x}%
_i\right) \leq 2\psi \left( \underset{i=1}{\overset{\infty }{\sum }}|\lambda
_i|e_i\right) +2\max \{|\lambda _i|:i\in \mathbb{N}\}
\end{equation*}
\end{proposition}

\begin{proof}
For $\phi \in K$, we choose $(K^s(\phi ))_{s=0}^m$ an analysis of $\phi $.
For each $i\in \mathbb{N}$ we denote by $f^i$ the unique element of $\underset{%
i=0}{\overset{m}{\bigcup }}K^s(\phi )$ such that\newline
(1) \textit{supp}$f^i\cap $\textit{supp}$\tilde{x}_i=$\textit{supp}$\phi
\cap $\textit{supp}$\tilde{x}_i$\newline
(2) $f^i\in K^{s_i}(\phi )$ and for $s<s_i$ there is no $f\in K^s(\phi )$
satisfying (1) for the vector $\tilde{x}_i$.\newline
We set $M=\left\{ i\in \mathbb{N}:w(f^i)\frac 1{m_{r_{i+1}}}\right\} $ and $I=%
\mathbb{N}\backslash M$ and for every $f\in \underset{s=0}{\overset{m}{\bigcup }%
}K^s(\phi )$ we set $D_f=\{i\in I:$\textit{supp}$f\cap $\textit{supp}$\tilde{%
x}_i=$\textit{supp}$\phi \cap $\textit{supp}$\tilde{x}_i\}$. We inductively
define for every $f\in \underset{s=0}{\overset{m}{\bigcup }}K^s(\phi )$ a
functional $\psi _f\in 2K_{T_3}$ such that, for every $i\in I$\newline
(a) \textit{supp}$\psi _f\subset D_f$\newline
(b) $f(\tilde{x}_i)\leq 2\psi _f(e_i)$ and either\newline
(c$_1$) $\psi _f\in K_{T_3}$ and $w(\psi _f^{})=w(f)$ or \newline
(c$_2)$There exists $\psi _f^1\in K_{T_3}$ and $i\in I$, $i<$\textit{supp}$%
\psi _f^1$ such that $w(\psi _f^1)=w(f)$ and $\psi _f=e_i^{*}+\psi _f^1$.%
\newline
For $s=0$ the construction is obvious. Suppose that for all $f\in \underset{%
t=0}{\overset{s-1}{\bigcup }}K^t(\phi )$, $\psi _f$ has been defined and let 
$f\in K^s(\phi )$, $f=\frac 1{m_q}(f_1+...+f_d)$, $f_1<...<f_d$ are
successive elements of $K^{s-1}(\phi )$. We set $J=D_f\backslash \underset{%
j=1}{\overset{d}{\bigcup }}D_i$.\newline
{\bf Case 1.} For every $i\in I$, $q\notin [r_i,r_{i+1})$. In this case
we simply observe from the definition of $I$ that $q<r_i$ for all $i\in
D_f\backslash \underset{j=1}{\overset{d}{\bigcup }}D_{f_j}$ Hence $\left| 
\underset{j=1}{\overset{d}{\sum }}f_j(\tilde{x}_i)\right| \leq 2$. We set $%
F_1=\{j\leq d:\psi _{f_j}\in K_{T_3}\}$. Then define $\psi _f=\frac
1{m_q}\left[ \underset{j\in F_1}{\sum }\psi _{f_j}+\underset{j\in F_1^c}{%
\sum }\psi _j^1+\underset{j\in F_1^c}{\sum }e_{i_j}^{*}+\underset{i\in J}{%
\sum }e_i^{*}\right] $. It can be easily checked that $\psi _f\in K_{T_3}$.%
\newline
{\bf Case 2.} There exists (a unique) $i_0\in J$ such that $r_i\leq
q<r_{i+1}$. In this case we set\newline
$\psi _f^1=\frac 1{m_q}\left[ \underset{j\in F_1}{\sum }\psi _{f_j}+%
\underset{j\in F_1^c}{\sum }\psi _{f_J}^1+\underset{j\in F_1^c}{\sum }%
e_{i_j}^{*}+\underset{i\in J\backslash \{i_0\}}{\sum }e_i^{*}\right] $ and $%
\psi _f=\psi _f^1+e_{i_0}^{*}$ and the inductive construction is complete.%
\newline
It follows that there exists $\psi =\psi _{\phi _{}}$ such that either $\psi
_{\phi _{}}\in K_{T_3}$ and $w(\psi _{\phi _{}})=w(\phi )$ or $\psi
=e_{i_0}^{*}+\psi _\phi ^1$ with $w(\psi _{\phi _{}}^1)=w(\phi )$ and $\phi (%
\tilde{x}_i)\leq 2\psi (e_i)$. The desired result follows easily from this
and the proof is complete.
\end{proof}

\begin{corollary}
\label{6.16}Let $(\tilde{x}_i)_{i=1}^\infty $ be a RIS in $\tilde{X}$, $Z=<(%
\tilde{x}_i)_{i=1}^\infty >$ and $S:T_3\rightarrow Z$ the linear operator
defined by $S(e_i)=x_i$ for every $i\in \mathbb{N}$. Then $S$ is a bounded
non-compact operator.
\end{corollary}

\begin{remark}
\label{R6.2}The space $T_3$ has property (P2). Indeed if $(\tilde{x}%
_i)_{i=1}^\infty $ is any block sequence in $T_3$ by Proposition \ref{6.6}
there exists a RI block sequence of $(\tilde{x}_i)_{i=1}^\infty $. Since
very subsequence of a RIS is RIS we easily get that $T_3$ satisfies (P2)
with $(C_k)_k=(n_k)_k$.
\end{remark}

\section{{\bf FINAL RESULTS-PROBLEMS}}

This final section contains the main results of the paper. These are
theorems concerning quotients of H.I. spaces and factorization of $\mathbf{a-%
}$thin operators through H.I. spaces. In the final part we prove similar
results for quotients of $\ell ^p$-saturated Banach spaces.

\begin{theorem}
\label{Th7.1}Let $A\,$ be a reflexive Banach space with an unconditional
basis satisfying property (P). Then there exists a reflexive H.I. Banach
space $X$ such that $A$ is a quotient of $X$.
\end{theorem}

\proof%
Recall that Property (P) is defined in Definition \ref{DP}. By Theorem \ref
{Th3.9} it follows that the set $W$ in $X_A$ is an $\mathbf{a-}$thin set for
an appropriate positive null sequence $\mathbf{a=}\left( a_n\right) _{n\in 
\mathbb{N}}$ ,that norms a subspace $Y$ of $X_A^{*}$ isometric to $A^{*}$.%
\newline
If $||\;||_n$ is the equivalent norm defined by Minkowski's gauge $%
2^nW+a_nB_{X_A}$ then by Proposition \ref{6.1} the $\left( A_{n_i},\frac
1{m_i}\right) ^{\sigma _{}}$-$d-$product of the sequence $\left( \left(
X_A,||\;||_n\right) \right) _{n\in \mathbb{N}}$ is reflexive.\newline
Since $W$ is a thin subset of $X_A$ by Theorem \ref{6.2} and Proposition \ref
{Prop2.3} the diagonal space $D$ is a H.I.space.The space $D$ as a subspace
of a reflexive space is in itself reflexive and further $A^{*}$ is
isomorphic to a subspace of $D^{*}$ (Proposition \ref{Prop2.4}) . It follows
that $D$ is the desired space $X$ and the proof is complete%
\endproof%

\begin{theorem}
\label{Th7.2}Suppose that $A$ is either $c_0$ or a reflexive space with an
unconditional basis such that every block subspace $B$ contains a further
block subspace $Z$ complemented in $A$ . Then $A$\ is a quotient of a H.I.
space $X.$In the late case $X$ can be chosen to be reflexive..
\end{theorem}

\proof%
The proof is completely analogous to the previous one. For the space $X_A$
and the set $W$, as defined in sections 5 for the space $c_0\left( \mathbb{N}%
\right) $ and section 3 for the reflexive space $A$, we have in both cases
that $W$ is a thin norming set.Then we continue the proof in the same manner
as before.%
\endproof%

Next, we list some classical spaces that are quotients of an H.I space.

\begin{corollary}
Every $L^p\left( \lambda \right) $, $1<p<\infty $,$\;\ell ^p\left( \mathbb{N}%
\right) $, $1<p<\infty $, is a quotient of a reflexive H.I. space. Further
on, Tsirelson's space $T$ together with its dual and Schlumprecht's space $S$
are quotients of a reflexive H.I.space.
\end{corollary}

\proof%
It was shown in section 3 (Remark\ref{Rem3.1}, Remark\ref{3}) that $%
L^p\left( \lambda \right) $, $1<p<\infty $,$\;\ell ^p\left( \mathbb{N}\right) $%
, $1<p<\infty $, satisfy the property (P) and as it is well known, they have
an unconditional basis. Therefore the result follows from Theorem \ref{Th7.1}
.\newline
For Tsirelson's space $T$, its dual and Schlumprecht's space $S$ it is known
that they satisfy the complemented subspace condition of Theorem \ref{Th7.2}
hence there also exists the corresponding reflexive H.I. space mapping onto
them%
\endproof%
.

\begin{remark}
(i) It follows by standard duality arguments that $L^p\left( \lambda \right) 
$, $1<p<\infty $,$\;\ell ^p\left( \mathbb{N}\right) $, $1<p<\infty $, the
spaces $T$, $T^{*}$, $S$, $S^{*}$ are isomorphic to a subspace of $X^{*}$
for a reflexive H.I. space $X$ .\newline
In particular, the imbedding of $\ell ^1\left( \mathbb{N}\right) $ into $X^{*}$
shows that the dual of a H.I. space is not necessarily arbitrarily
distortable. Recall that in \cite{To} it has been shown that every H.I.
space is arbitrarily distortable.\newline
(ii) As we mentioned in the Introduction it is not possible $\ell ^1\left( 
\mathbb{N}\right) $ to be a quotient of a H.I. space. What's more, it is not
possible $c_0\left( \mathbb{N}\right) ,L^1\left( \lambda \right) $ to be
isomorphic to a subspace of $X^{*}$ for $X$ a H.I. space.\newline
Indeed, it is well known that if $X^{*}$ contains some of $c_0\left( \mathbb{N}%
\right) $, $L^1\left( \lambda \right) $ then $\ell ^1\left( \mathbb{N}\right) $
is isomorphic to a subspace of $X$ which is impossible if $X$ is a H.I.
space.
\end{remark}

We pass now to prove the dichotomy that we mentioned in the Introduction
concerning the quotients of H.I. spaces.

\begin{theorem}
Let $A$ be a Banach space. Then either $\ell ^1\left( \mathbb{N}\right) $ is
isomorphic to a subspace of $A$ or there exists $B$ infinite dimensional
closed subspace of $A$ which is a quotient of a H.I. space.
\end{theorem}

\proof%
Suppose that $\ell ^1\left( \mathbb{N}\right) $ is not isomorphic to a subspace
of $A$.Then by Gowers' Dichotomy Theorem (\cite{G}) there exists $Z$
subspace of $A\,$ which is either a H.I. space or it has an unconditional
basis.\newline
If the first case occurs then we have finished.\newline
In the second case, by James Theorem (\cite{J}), either $Z$ is a reflexive
space or it contains a subspace isomorphic to $c_0\left( \mathbb{N}\right) $.%
\newline
If $c_0\left( \mathbb{N}\right) $ is a subspace of $A$ then as we showed in
Theorem \ref{5.1} $c_0\left( \mathbb{N}\right) $ is a quotient of a H.I. space
and therefore we get the desired result.\newline
The remaining case is that of the space $Z$ being reflexive.\newline
It follows from Theorem \ref{Th4.04} that $Z$ has a subspace $B$ such that
the set $W$ in the space $X_B$ is an $\mathbf{a}$-thin subset of $X_B$ and
norming for a subspace of $X_B^{*}$ isometric to $B^{*}$. As we have seen
before (Theorem \ref{Th7.1}) this implies that $B$ is a quotient of a H.I.
space so the proof is complete%
\endproof%
.

The next result is related to the factorization of linear operators between
Banach spaces through H.I. spaces.\newline
We recall that a bounded linear operator $T:X\rightarrow Y$ is $\mathbf{a}$%
-thin if $T[B_X]$ is an $\mathbf{a}$-thin subset of $Y$, for a positive null
sequence $\mathbf{a}=\left( a_n\right) _{n\in \mathbb{N}}$.

\begin{theorem}
\label{Th7.5}If $X$, $Y$ are Banach spaces and $T:X\rightarrow Y$ is an $%
\mathbf{a}$-thin operator then there exists $D$ a H.I. space such that $T$
is factorized through $D$.
\end{theorem}

\proof%
Set $W=T[B_X]$, which is by our assumption an $\mathbf{a}$-thin subset of $Y$%
, where $\mathbf{a=}\left( a_n\right) _{n\in \mathbb{N}}$ is a positive null
sequence. Denote by $||\;||_n$ the equivalent norm on $Y$ defined by the
Minkowski gauge of the set $2^nW+a_nB_Y$. We assume the $\left( \mathcal{A}%
_{n_i},\frac 1{m_i}\right) ^{\sigma _{}}$-product of the sequence $\left(
\left( Y,||\;||_n\right) \right) _{n\in \mathbb{N}}$ which is by Theorem \ref
{6.2} a block H.I. space. The diagonal space $D$ of that $d$-product is by
Proposition \ref{Prop2.3} a H.I. space and further on the ball of $D$
contains the set $W$. Therefore the operator $Q:X\rightarrow D$, defined by $%
Q\left( x\right) =T\left( x\right) $, is a bounded linear operator and $%
T=j\circ Q$, where $j:D\rightarrow Y$, is the natural inclusion map.%
\endproof%

\begin{theorem}
(a) For every $r\neq p$, $1\leq r,p<\infty $ and $T\in \mathcal{L}\left(
\ell ^r,\ell ^p\right) $, $T$ is factorized through a H.I. Banach space.%
\newline
This remains valid if one of $\ell ^r\left( \mathbb{N}\right) $, $\ell ^p\left( 
\mathbb{N}\right) $ is substituted by $c_0\left( \mathbb{N}\right) $.\newline
(b) Every $T\in \mathcal{L}\left( \ell ^r\left( \mathbb{N}\right) ,\ell
^p\left( \mathbb{N}\right) \right) $ is factorized through a H.I. space.\newline
(c) The identity map $I:L^\infty \left( \lambda \right) \rightarrow
L^1\left( \lambda \right) $ is factorized through a H.I. space.
\end{theorem}

\proof%
(a) and (b) follow from the above Theorem \ref{Th7.5} and Corollary \ref
{3.13}.\newline
(c) follows from the theorem above and Proposition \ref{M} that asserts that 
$B_{L^\infty }$ is an $\mathbf{a}$-thin set%
\endproof%
.

The next result concerns the structure of $\mathcal{L}\left( X,X\right) $
for $X$ a H.I. space. It was shown in \cite{GM} that any such operator is of
the form $\lambda I+S$, where $I$ is the identity map and $S$ is a strictly
singular operator. It is an open question if there exists a H.I. space such
that every bounded linear operator is of the form $\lambda I+K$ where $K$ is
a compact operator .W.Gowers has constructed a subspace $Y$ of Gowers-Maurey
space $X$ and an operator $T:Y\rightarrow X\ $which is strictly singular but
not compact.\newline
In the next Theorem we show that there exists H.I. spaces $X$ with ``many '' 
$T:X\rightarrow X$ which are strictly singular and not compact.

\begin{theorem}
There exists a H-I space $X$ such that for every infinite dimensional closed
subspace $Z$ of $X$ there exists a strictly singular operator $%
T:X\rightarrow Z$ which is not compact.
\end{theorem}

\begin{proof}
Let $T_3$ denote the space $T\left[ \left( \mathcal{A}_{3n_i},\frac
1{m_i}\right) _{i=1}^\infty \right] $ then by Remark \ref{R6.2} the set $W$
in $X_{T_3}$ is a $\mathbf{a}$-thin set hence the diagonal space $\Delta $
defined by the set $W$ and the connecting norm defined by the family $\left( 
\mathcal{A}_{n_i},\frac 1{m_i}\right) _{i\in \mathbb{N}}^{\sigma _{}}$ is an
H-I space that has as quotient the space $T_3$. Given any $Z$ closed
subspace of $\Delta $. Then there exists a normalized basic sequence in $Z$
equivalent to a block subsequence $(y_i)_{i=1}^\infty $ of $\tilde{X}$. It
follows for Corollary \ref{6.16} that there exists a non compact operator $%
S:T_3\rightarrow Z$. If $Q$ denotes the surjection of $\Delta $ onto $T_3$
then $S\circ Q$ is a non compact strictly singular operator from $\Delta $
to $Z$ and the proof is complete.
\end{proof}

A by-product of our method is the following result related to $\ell _p(\mathbb{N%
})$ $\left( 1<p<\infty \right) $ and $c_0\left( \mathbb{N}\right) $-saturated
Banach spaces. First we recall their definition.

\begin{definition}
A Banach space $X$ is said to be $\ell _p(\mathbb{N})-${\bf saturated }($%
c_0\left( \mathbb{N}\right) ${\bf -saturated} ) if every subspace $Y$ of $X$
contains a further subspace isomorphic to $\ell _p(\mathbb{N})$ ($c_0\left( 
\mathbb{N}\right) $ ).
\end{definition}

\begin{theorem}
Let $A$ be a reflexive Banach space with an unconditional basis .Then there
exists $B$ subspace of $A$ such that for every $p\in \left( 1,\infty \right) 
$ there exists $X_p$ $\ell _p(\mathbb{N})-$saturated reflexive space which has
as quotient the space $B$.Further there exists $c_0\left( \mathbb{N}\right) $%
-saturated space with the same property.
\end{theorem}

\proof%
It follows from Theorem \ref{Th4.04} that there exists $B$ subspace of $A$
with an unconditional basis such that the closed bounded symmetric set $W$
in $X_{A\text{ }}$ is an $\mathbf{a-}$ thin set ,norming a subspace $Y$ of $%
X_A^{*}$ isometric to $B^{*}.$ We denote by $||\;||_n$ the equivalent norm
defined on $X_A$ by the set $2^nW+a_nB_{X_A}$ and by $X_n$ the space $\left(
X,||\;||_n\right) .$ Then each $X_n$ is a reflexive space and we have the
same property for the space $\left( \sum_{n=1}^\infty \bigoplus X_n\right)
_p,1<p<\infty $ . Since $W$ is a thin set ,as follows from Proposition \ref
{2.3} ,the diagonal space is $\ell _p(\mathbb{N})-$saturated and if we consider
the space $\left( \sum_{n=1}^\infty \bigoplus X_n\right) _0$ then the
diagonal space is also $c_0$-saturated .Also $B^{*}$ is isomorphic to a
subspace of $D^{*}.$ Therefore $B$ is a quotient of the space $D$ and the
proof is complete%
\endproof%

\begin{remark}
(a) It follows from the above Theorem that for every $p,r\in \left( 1,\infty
\right) $ ,$\ell ^p\left( \mathbb{N}\right) $ is a quotient of a reflexive $%
\ell ^p\left( \mathbb{N}\right) $-saturated space. The result that $\ell
^2\left( \mathbb{N}\right) $ is a quotient of a $c_0\left( \mathbb{N}\right) $%
-saturated space has been proved with a different method by D.Leung (\cite
{Le}).\newline
(b) It can be shown that also $L^p\left( \lambda \right) ,1<p<\infty ,$
Tsirelson's space $T$ ,Schlumprecht's $S$ are quotient of reflexive $\ell
^p\left( \mathbb{N}\right) $-saturated or $c_0\left( \mathbb{N}\right) $-saturated
spaces
\end{remark}

We close this section by a list of questions related to our results.

\begin{center}
{\bf Problems}
\end{center}

\begin{enumerate}
\item  Given $X$ a separable Banach space not containing a subspace
isomorphic to $\ell ^1\left( \mathbb{N}\right) $. Is it true that $X$ is a
quotient of a H.I. space ?

\item  If $X$ is an H.I space does there exist a $Y$ subspace of $X^{*}$
which is a H.I. space ?

\item  If $X$ is an H.I space does there exist a $Y$ subspace of $X$ such
that $Y^{*}$ is a H.I. space ?

\item  Does there exists a bounded linear operator $T:X\rightarrow Y$ which
is strictly singular but not $\mathbf{a}$-thin for every positive null
sequence $\mathbf{a}$ ?

\item  Can the existence of quotient maps from H.I. spaces onto reflexive
spaces be used for the solution of the distortion problem for general
reflexive spaces ?
\end{enumerate}

S.A. Argyros and V. Felouzis ,
Department of Mathematics \\     
University of Athens Greece ,
e-mail sargyros@atlas.uoa.gr

\begin{thebibliography}{ADKM}
\bibitem[AA]{AA}  D.Alspach and S.Argyros.\textit{Complexity of weakly null
sequences,}Dissertationes Mathematicae CCCXXI (1992)

\bibitem[{AD1}]{AD1}  S. Argyros and I. Deliyanni. \textit{Banach spaces of
the type of Tsirelson}, preprint 1992.

\bibitem[{AD2}]{AD2}  S. Argyros and I. Deliyanni. \textit{Examples of
asymptotic} $\ell _1$ \textit{Banach spaces}. Trans. AMS 349 (1997) 973-995.

\bibitem[{ADKM}]{ADKM}  S. Argyros, I. Deliyanni, D. Kutzarova and A.
Manoussakis, \textit{Modified mixed Tsirelson spaces }preprint (1997).

\bibitem[AMT]{AMT}  S.Argyros,S.Merkourakis and A.Tsarpalias \textit{Convex
unconditionality and summability of weakly null sequences} (to appear in
Isr.J.Math.)

\bibitem[AO]{AO}  G.Androulakis and E.Odell \textit{Distorting mixed
Tsirelson spaces(}Preprint)

\bibitem[{B}]{B}  S.F. Bellenot. \textit{Tsirelson} \textit{superspaces and} 
$\ell _p$. Journ. of Funct. Analysis. 69 No2, 1986, 207-228.

\bibitem[Bo1]{Bo}  J.Bourgain.\textit{Convergent sequences of continuous
functions}, Bull .Soc . Math . Belg . Ser.B 32(1980),235-249.

\bibitem[Bo2]{Bo2}  J.Bourgain.\textit{La propri\'{e}t\'{e} de Radon-Nicodym 
}Math.Univ.Pierre et Marie Curie 36 (1979)

\bibitem[D]{D}  J.Diestel.\textit{Sequences and series in Banach spaces,}%
Graduate texts in Math.92,Springer-Verlag 1984

\bibitem[{DFJP}]{DFJP}  W.J. Davis, T. Figiel, W.B. Johnson and A.
Pelczynski. \textit{Factoring Weakly Compact Operators}. Journ. of Funct.
Analysis 17 (1974) 311-327.

\bibitem[{F1}]{F1}  V. Ferenczi. \textit{Quotient Hereditarily Idecomposable
spaces} preprint

\bibitem[{F2}]{F2}  V. Ferenczi. \textit{A uniformly convex heredirarily
indecomposable Banach space} preprint.

\bibitem[J]{J}  R.C.James.\textit{Bases and reflexivity of Banach spaces.}%
Ann.of.Math.52\newline
(1950),518-527

\bibitem[{G}]{G}  W.T. Gowers. \textit{A new dichotomy for Banach spaces. }%
GAFA 6(1996) 1083-1093.

\bibitem[Gr]{Gr}  A.Grothendieck.\textit{Crit\`{e}res de compacit\'{e} dans
les espaces fonctionelles g\'{e}n\'{e}raux,}Amer.J.Math.74 (1952) 168-186

\bibitem[{GM}]{GM}  W.T. Gowers and B. Maurey. \textit{The unconditional
basic sequence problem} Journal of AMS 6, 1993, 851-874

\bibitem[Ke]{Ke}  A.Kechris \textit{Classical Descriptive Set Theory }%
Springer-Verlag, Berlin, (1994).

\bibitem[KM]{KM}  K. Kuratowski and A. Mostowski. \textit{Set Theory}.
Amsterdam (1968).

\bibitem[{KW}]{KW}  N. Kalton and A. Wilanski. \textit{Tauberian operators
on Banach spaces} Proc. A.M.S. 57 (1967) 251-255.

\bibitem[{L}]{L}  H.E. Lacey. \textit{The Isometric Theory of Classical
Banach Spaces, }Springer-Verlag, Berlin, 1974.

\bibitem[{Le}]{Le}  D.H. Leung. \textit{On }$c_0$-\textit{saturated spaces. }%
Illinois J. Math. 39 (1995) 15-29.

\bibitem[{LT}]{LT}  J. Lindenstrauss and L. Tzafriri. \textit{Classical
Banach spaces I} Springer Verlag 92, 1977.

\bibitem[{LP}]{LP}  J. Lindenstrauss and R.R. Phelps. \textit{Extreme point
properties of convex bodies in Reflexive Banach spaces}.Israel J. Math., 6,
(1968) 39-48.

\bibitem[MR]{MR}  B.Maurey and H.Rosenthal .\textit{Normalized weakly null
sequences with no unconditional subsequence }Studia Math.61 (1977) 77-98

\bibitem[MS]{MS}  V.Milman and G.Schechtman .\textit{Asymptotic Theory of
Finite Dimensional Normed Spaces.}Lecture Notes in Math. 1200 Springer-Verlag

\bibitem[{N}]{N}  R. Neidinger. \textit{Factoring Operators through
hereditarily-}$\ell _p$ \textit{spaces} Lecture Notes in Math. 1166 (1985).

\bibitem[N2]{N2}  R. Neidinger. \textit{Properties of Tauberian Operators.}
Dissertation. University of Texas at Austin 1984.

\bibitem[{NR}]{NR}  R. Neidinger and H.P. Rosenthal. \textit{%
Norm-attainement of linear functional on subspaces and characterizations of
Tauberian operators}. Pac. J. Math. 118 (1985) 215-228.

\bibitem[{To}]{To}  N. Tomczak-Jaegermann. \textit{Banach spaces of type }$p$
\textit{have arbitrarily distortable subspaces }GAFA 6 (1996) 1075-1082.

\bibitem[T]{T}  A.Tsarpalias \textit{A note on Ramsey property },Proc.AMS
(to appear)

\bibitem[Ts]{Ts}  B.S. Tsireslon. \textit{Not every Banach space contains} $%
\ell _p$ \textit{or} $c_0$ Funct. Anal. Appl. 8 (1974) 138-141.
\end{thebibliography}
\end{document}